\documentclass[preprint,sort&compress,11pt,3p]{elsarticle}

\usepackage{moreverb}
\usepackage{amsmath}
\newtheorem{remark}{Remark}
\usepackage{algorithmic}
\usepackage{graphics}
\usepackage{graphicx}
\usepackage{subfigure}
\usepackage{enumitem}
\usepackage{color}
\usepackage{framed}
\usepackage{wrapfig}
\usepackage{bm}
\usepackage{mathrsfs}
\usepackage{amsfonts}
\usepackage{multirow}
\usepackage{longtable}
\usepackage{subfigure}
\usepackage{hyperref}
\usepackage{paralist}
\usepackage{indentfirst}
\usepackage{relsize}
\usepackage{extarrows}
\usepackage{lineno,xcolor}
\usepackage{upgreek}
\usepackage{bm}
\usepackage{diagbox}
\usepackage{xcolor}
\usepackage{booktabs}
\usepackage{tabularx}
\usepackage{multirow}
\usepackage[flushleft]{threeparttable}
\usepackage[linesnumbered,lined,ruled]{algorithm2e}

\newcommand{\eref}[1]{(\ref{#1})}

\hypersetup{
bookmarks=true,
bookmarksopen=true,
bookmarksnumbered=true,
unicode=false,
pdftoolbar=true,
pdfmenubar=true,
pdffitwindow=false,
pdfstartview={FitH},
pdftitle={My title},
pdfauthor={Author},
pdfsubject={Subject},
pdfcreator={Creator},
pdfproducer={Producer},
pdfkeywords={keywords},
pdfnewwindow=true,
colorlinks=true,
linkcolor=blue,
citecolor=blue,
filecolor=magenta,
urlcolor=blue
}

\journal{Journal}

\begin{document}
	
\title{\LARGE Physics informed deep learning for computational elastodynamics without labeled data}

\author[NU1]{Chengping~Rao}
\author[NU2,MIT]{Hao~Sun}
\author[NU1]{Yang~Liu\corref{cor}}
\ead{yang1.liu@northeastern.edu}

\cortext[cor]{Corresponding author. Tel: +1 617-373-3888}

\address[NU1]{Department of Mechanical and Industrial Engineering, Northeastern University, Boston, MA 02115, USA}
\address[NU2]{Department of Civil and Environmental Engineering, Northeastern University, Boston, MA 02115, USA}
\address[MIT]{Department of Civil and Environmental Engineering, MIT, Cambridge, MA 02139, USA}

\begin{abstract}
	\small
Numerical methods such as finite element have been flourishing in the past decades for modeling solid mechanics problems via solving governing partial differential equations (PDEs). A salient aspect that distinguishes these numerical methods is how they approximate the physical fields of interest. Physics-informed deep learning (PIDL) is a novel approach developed in recent years for modeling PDE solutions and shows promise to solve computational mechanics problems without using any labeled data (e.g., measurement data is unavailable). The philosophy behind it is to approximate the quantity of interest (e.g., PDE solution variables) by a deep neural network (DNN) and embed the physical law to regularize the network. To this end, training the network is equivalent to minimization of a well-designed loss function that contains the residuals of the governing PDEs as well as initial/boundary conditions (I/BCs). In this paper, we present a physics-informed neural network (PINN) with mixed-variable output to model elastodynamics problems without resort to the labeled data, in which the I/BCs are hardly imposed. In particular, both the displacement and stress components are taken as the DNN output, inspired by the hybrid finite element analysis, which largely improves the accuracy and the trainability of the network. Since the conventional PINN framework augments all the residual loss components in a ``soft'' manner with Lagrange multipliers, the weakly imposed I/BCs cannot not be well satisfied especially when complex I/BCs are present. To overcome this issue, a composite scheme of DNNs is established based on multiple single DNNs such that the I/BCs can be satisfied forcibly in a ``hard'' manner. The propose PINN framework is demonstrated on several numerical elasticity examples with different I/BCs, including both static and dynamic problems as well as wave propagation in truncated domains. Results show the promise of PINN in the context of computational mechanics applications.

\end{abstract}


\maketitle

\section{Introduction}\label{sintro}

The principled modeling of physical systems (e.g., materials) plays a critical role in scientific computing, where numerical methods are commonly used. The general idea behind numerical methods is to establish an approximate solution through a finite set of basis functions and unknown parameters (or variables). The unknown parameters can be solved through a system of algebraic equations after the discretization of the domain. The selection of the basis function leads us to different numerical methods. The classical finite element method (FEM) employs piecewise polynomials to approximate the solution with the nodal displacement and/or stress components as the unknown variables. By contrast, in the framework of isogeometric analysis (IGA) \cite{hughes2005isogeometric}, the non-uniform rational basis splines (NURBS) are used for the solution approximation with the variables of interest at control points acting as the unknown to be solved. Chebyshev and Fourier series are mostly employed in the spectral method as the basis functions \cite{gopalakrishnan2007spectral}.

In recent years, advances in deep learning have attracted drastic attention to the field of computational modeling and simulation of physical systems, thanks to the rich representations of deep neural networks (DNNs) for learning complex nonlinear functions. Latest studies that leverage DNNs for physical modeling branch into two streams: (1) the use of experimental or computationally-generated data to create coarse-graining reduced-fidelity or surrogate models \cite{WANG2018337,wang2019meta, HAN2019112603, rao2020three,vlassis2020geometric}, and (2) physics-informed neural network (PINN) for modeling the solution of partial differential equations (PDEs) that govern the behavior of physical systems \cite{sirignano2018dgm,raissi2019physics,raissi2020hidden}. The former requires rich and sufficient data to learn a reliable generative model and typically fails to satisfy physical constraints, whereas the latter relies only on small or even zero labeled datasets and enables data-scarce, physics-constrained learning. The embedded physics is expected to provide constraints to the trainable parameters, alleviate overfitting issues, reduce the need of big training datasets, and thus improve the robustness of the trained model for reliable prediction. In fact, the idea of using neural networks to solve PDEs is not new and can date back to the last century \cite{lee1990neural, psichogios1992hybrid, lagaris:1998, lagaris:2000}. These early works rely on the function approximation capabilities of a feedforward fully-connected neural networks to solve initial/boundary value problems. The solution to the system of equations can be obtained through minimization of the network's loss function, which typically consists of the residual error of the governing equations along with initial/boundary values. More recently, Raissi \emph{et al.} \cite{raissi2018deep,raissi2019physics,raissi2019deepJFM,raissi2020hidden} has inherited and extended this concept, leveraged the strong expressibility of DNNs, and developed the general PINN framework to solve the forward and inverse problems involving the system of nonlinear PDEs with small datasets or even without any labeled data. The PINN has found vast applications, within a short time, in a wide range of physical problems including modeling fluid flows and Navier-Stokes equations \cite{raissi2020hidden,sun2020surrogate,rao2020physics,jin2020nsfnets,mao2020physics,gao2020phygeonet}, solving stochastic PDEs \cite{karumuri2020simulator}, flows transport in porous media \cite{zhu2018bayesian, tripathy2018deep,he2020physics2}, cardiovascular systems \cite{sun2020physics,kissas2020machine,sahli2020physics}, design of metamaterials \cite{fang2019deep,liu2019multi,chen2020physics}, metamodeling of nonlinear structural systems \cite{zhang2020physics,ZhangCMAME2020}, and discovery of physical laws \cite{BERG2019239,Both2019,chen2020deep}, among others. To further improve the learning performance, the PINN framework has also been extended via incorporating variational/energy formulations of the residual physics loss function \cite{weinan2018deep,kharazmi2019variational,samaniego2020energy,goswami2020transfer}, distributed learning using domain decomposition \cite{kharazmi2020hp,jagtap2020conservative}, and uncertainty quantification via variational/Bayesian inference \cite{yang2019adversarial,sun2020physics,yang2020b,zhang2019quantifying}. It is important to note that a few recent attempts show the promise and power of PINN for addressing computational mechanics relevant challenges such as solving mechanical problems \cite{samaniego2020energy, he2020physics}, modeling fracture in materials \cite{goswami2020transfer, goswami2020adaptive}, and detecting cracks via ultrasound nondestructive testing \cite{shukla2020physics}.


The main contribution of this paper is to develop the PINN framework for modeling elastodynamics problems (e.g., wave propagation in bounded or truncated domains) in the absence of labeled data. In particular, a feedforward DNN is used as the a global approximator of the concerned physical quantities such as displacement and stress field. In this sense, the PINN shares a salient feature with the spectral method since both are global methods. In such as way, the DNN can be also viewed as a mapping from the independent spatiotemporal variable $\mathbf{X}=(\mathbf{x},t)$ to the determined variable $\mathbf{Y}=(\mathbf{u},\bm{\sigma})$, denoted by $\mathcal{N}(\mathbf{W},\mathbf{b}):\mathbf{X}\mapsto\mathbf{Y}$, where $\mathcal{N}$ denotes the DNN with trainable weights $\mathbf{W}$ and biases $\mathbf{b}$. The reason why DNN is well qualified as the solution approximator includes: (1) its extraordinary approximation capability proven by the universal approximation theorem \cite{cybenko1989approximation}, (2) its infinite continuity property \cite{lagaris:2000} with proper activation function, (3) the derivatives that appears in the PDEs can be calculated exactly via automatic differentiation \cite{baydin2017automatic}, and (4) a great variety of deep learning frameworks, such as TensorFlow \cite{abadi2016tensorflow} and PyTorch \cite{NEURIPS2019_9015}, make the implementation and parallelization efficient and convenient. 

The remaining of this paper is organized as follows. In Section \ref{dnn_ad}, the basic knowledge of DNN and automatic differentiation, crucial to formulation of the loss function, is introduced briefly. In Section \ref{pinnassolver}, we elaborate how PINN can be employed as a general PDE solver. The formulation of a constrained minimization problem equivalent to solving the PDEs is presented. The elasticity theory is presented in Section \ref{elas_theory}, altogether with the framework of PINN and the loss function for solving the elastodynamics problem. We propose a composite scheme of DNNs for the construction of a synergy solution to the elastodynamics problem. One of the most significant benefits of the constructed synergy solution is that the initial/boundary conditions (I/BCs) will be satisfied forcibly in a ``hard'' manner. In addition, the mixed-variable output of PINN is proven to be crucial for the training. Several numerical examples, including the defected plate under cyclic uni-axial tension and the elastic wave propagation in confined and truncated (e.g., infinite and semi-infinite) domains, are given in Section \ref{results} to illustrate the capability of the proposed PINN framework for modeling elastodynamics problems. Section \ref{conc} is dedicated to the conclusion of the current work and the outlook of our future work.

\section{Method}\label{method}

In this section, the proposed framework of PINN for solving elastodynamics problems is presented. The basic concepts of DNN and automatic differentiation, which are the prerequisites for designing a PINN, are introduced briefly. The construction of the loss function for the training of the PINN, as well as the discretization of the problem domain, is introduced subsequently. The employed mixed-variable formula and the construction of a synergy solution are also elaborated. 

\subsection{Deep neural network and automatic differentiation}\label{dnn_ad}

In recent years, DNN has led to many successful applications such as image recognition and natural language processing thanks to its exceptional expressibility. A feedforward fully-connected neural network can be assumed to be the stack of the input layer, multiple hidden layers and the output layer. The connection between two adjacent layers, say from $(i-1)$th to $i$th layer, can be expressed concisely in the form of tensor, as follows
\begin{equation} 
\label{feedforwardnn} 
\begin{aligned}
 \mathbf{z}_{i}={\sigma}\left(\mathbf{b}_{i}+{\mathbf{W}_{i}}\mathbf{z}_{i-1}\right ),~\text{for }1\leq i\leq n+1
\end{aligned}
\end{equation}
\noindent where $n$ is the total number of hidden layers, ${\rm \sigma(\cdot)}$ denotes activation function acting element-wise, $\mathbf{z}_0$ and $\mathbf{z}_{n+1}$ denotes the input and output tensors respectively, and $\mathbf{W}_{i}$ and $\mathbf{b}_{i}$ are the trainable weight matrix and bias vector in the $i$ th layer. In this work, we utilize DNN as the parameterized approximate solution to the elastodynamics problems in which the spatiotemporal location $\mathbf{X}=(\mathbf{x},t)$ denotes the independent variables. To find a set of trainable parameters that achieve good approximation to the solution, the loss function with the physical law embedded needs to be designed as the training target. In contrast with the Galerkin method widely used in computational mechanics, we concentrate on the strong form of the governing PDEs throughout the paper. We will show the benefits brought by handling the strong-form equation in an elastic wave propagation problem. The design of the loss function is presented in details in Section \ref{pinnassolver}.
 
\begin{figure}[t!]
	\centering
	\includegraphics[width=0.7\textwidth]{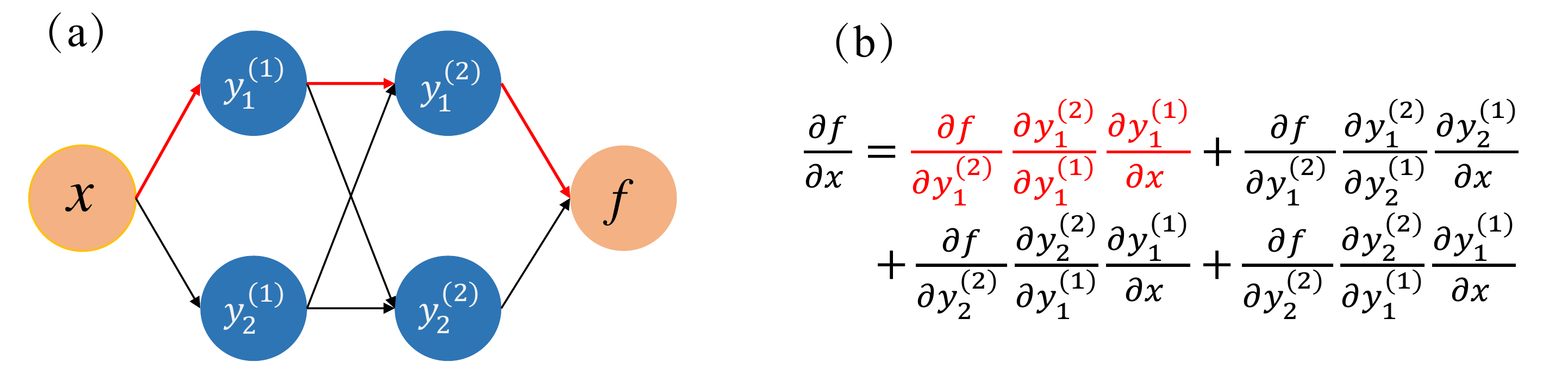}
	\caption{Diagram of automatic differentiation: (a) A simple feedforward neural network with two layers and two neurons in each layer; (b) The analytical expression of $\frac{\partial f}{\partial x}$ computed from the neural network. Each term corresponds to a path from $x$ to $f$.}
	\label{sec2-auto-diff}
\end{figure}

To construct the loss function of the concerned PDEs in PINN, we need to evaluate the partial derivatives of the physical field with regard to the spatiotemporal variables, such as $\nabla\cdot\bm{\sigma}$ and $\mathbf{u}_{tt}$. The automatic differentiation \cite{baydin2017automatic} is able do this job perfectly. A simple feedforward neural network, as shown in Fig.\ref{sec2-auto-diff}(a), is provided as an example to illustrate how it works. The network takes the variable $x$ as input and outputs $f$ after the nonlinear transformation defined in Eq. \eref{feedforwardnn}. Here, $y_{j}^{(i)}$ denotes the output of the $j$th neuron in the $i$th hidden layer. Since the output $f$ can be represented as a nested function of $x$, we can apply the chain rule to calculate the derivative of $f$ with respect to $x$, as depicted in Fig. \ref{sec2-auto-diff}(b) where the highlighted term corresponds the red path on the graph Fig. \ref{sec2-auto-diff}(a). The network can be implemented on the platform of TensorFlow  \cite{abadi2016tensorflow} which supports the definition of the partial derivatives in a symbolic way. Hence, unlike the numerical differentiation techniques that suffers from the approximation error, the automatic differentiation produces the exact derivatives (except round-off error) from the computational graph.

\subsection{PINN as PDE solver}\label{pinnassolver}
In this subsection, we focus on the discretization strategy and the formulation of PINN for solving general PDEs. Let us consider a nonhomogeneous second-order one-dimensional PDE subject to initial/boundary condition as an example, given by
\begin{equation} 
\label{pde} 
\mathcal{L}(u)=f(x,t),\,x\in \Omega,\,t\in [0, T]
\end{equation}
\begin{equation} 
\label{boundarycd} 
\begin{aligned}
u(x, t)&={\mathcal{B}_D}(t),~x\in\partial\Omega_D,~t\in [0, T] \\
\mathbf{n}\cdot\nabla u(x, t)&=\mathcal{B}_N(t),~x\in\partial\Omega_N,~t\in [0, T]
\end{aligned}
\end{equation}
\begin{equation} 
\label{initialcd} 
\begin{aligned}
u(x, 0)&=\mathcal{I}_0(x),~x\in\Omega \\
u_t(x, 0)&=\mathcal{I}_1(x),~x\in\Omega
\end{aligned}
\end{equation}
\noindent where $u(x, t)$ is the solution to this PDE, $\mathcal{L}(\cdot)$ is a differential operator and $\Omega$ is the domain of interests, $\partial\Omega=\partial\Omega_D\cup\partial\Omega_N $ is the boundary of the domain composed exclusively by Dirichlet boundary $\partial\Omega_D$ and Neumann boundary $\partial\Omega_N$, $\mathbf{n}$ is the unit outer normal vector of the boundary, $\mathcal{B}_D(t)$ and $\mathcal{B}_N(t)$ are the prescribed functions for two types of boundaries (i.e., Dirichlet and Neumann), $\mathcal{I}_0(x)$ and $\mathcal{I}_1(x)$ are initial states of the domain.

\begin{figure}[t!]
	\centering
	\includegraphics[width=0.85\textwidth]{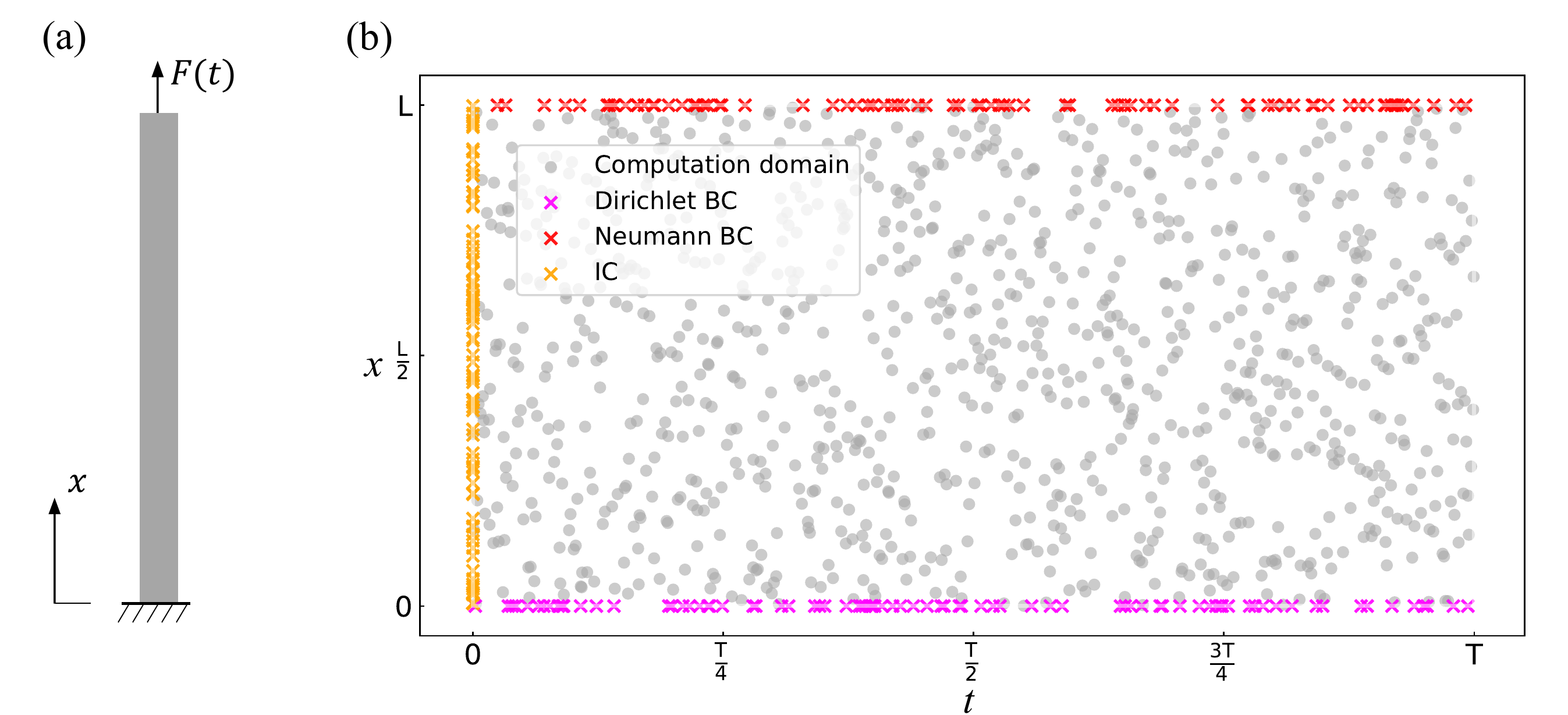}
	\caption{Diagram of discretization by collocation points in PINN  (a) A simple 1D bar subject to Dirichlet and Neumann boundary conditions at two ends; (b) The set of collocation points generated for the evaluation of initial/boundary value and equation residual.}
	\label{sec2-col-pt}
\end{figure}

As mentioned previously, the solution to the PDE can be approximated by a DNN, i.e., $\hat{u}(x,t)\equiv\mathcal{N}(x,t|\mathbf{W},\mathbf{b})$. We could view the initial and boundary value problem (IBVP) of the PDE being solved if we can find a set of DNN parameters that make the residual of the equation and I/BCs equal (or close) to zero. To this end, a finite set of spatiotemporal collocation points in the computation domain are introduced for evaluation of the residuals (see Fig. \ref{sec2-col-pt}). The whole set of collocation points can be denoted by $S=\left \{ S_\text{D}, S_\text{DBC}, S_\text{NBC}, S_\text{IC} \right\}$ consisting of spatiotemporal coordinates for the entire spatiotemporal domain ($S_\text{D}$), Dirichlet boundary condition ($S_\text{DBC}$), Neumann boundary condition ($S_\text{NBC}$) and the initial condition ($S_\text{IC}$), given by
\begin{equation} 
\label{collocationpt} 
\begin{aligned}
S_\text{D}&=\left\{(x,t)|x\in\Omega \text{ and } t\in[0,T]\right\}\\
S_\text{DBC}&=\left\{(x,t)|x\in\partial\Omega_D \text{ and } t\in[0,T]\right\}\\
S_\text{NBC}&=\left\{(x,t)|x\in\partial\Omega_N \text{ and } t\in[0,T]\right\}\\
S_\text{IC}&=\left\{(x,t=0)|x\in\Omega\right\}
\end{aligned}
\end{equation}
Note that the collocation points can be generated via the Latin hypercube sampling (LHS) strategy \cite{mckay1979comparison}. Plugging the DNN-approximated solution $\hat{u}(x,t; \mathbf{W,b})$ (for simplicity, we denote it by $\hat{u}(x,t)$) into the PDE and I/BCs renders us a constrained optimization problem as follows
\begin{equation} 
\label{maximization} 
\begin{aligned}
(\mathbf{W^*,b^*}) = &\arg\underset{(\mathbf{W,b})}{\min}\Bigg\{ \sum_{(x_i,t_i)\in S_\text{D}}\left[\mathcal{L}(\hat{u}(x_i,t_i))-f(x_i,t_i)\right]^2 \Bigg\} \\
\text{s.t.}~~~~ &\hat{u}(x_i, t_i)={\mathcal{B}_D}(t_i), ~~~(x_i,t_i)\in S_\text{DBC}\\
&\mathbf{n}\cdot\nabla \hat{u}(x_i, t_i)=\mathcal{B}_N(t_i), ~~~(x_i,t_i)\in S_\text{NBC}\\
&\hat{u}(x_i, t_i=0)=\mathcal{I}_0(x_i), ~~~(x_i,t_i)\in S_\text{IC}\\
&\hat{u}_t(x_i, t_i=0)=\mathcal{I}_1(x_i), ~~~(x_i,t_i)\in S_\text{IC}
\end{aligned}
\end{equation}
where $(\mathbf{W^*,b^*})$ are the target parameters that make the DNN best approximate the solution to the PDE, which can be obtained by DNN training using optimization strategies such as the family of gradient descent methods. However, how to enforce I/BC equality constraints is non-trivial since they guarantee the uniqueness of the solution. In our work, a composite scheme of PINN is proposed to ensure that the I/BCs are satisfied exactly for the solution. The detailed discussion on hard imposition OF I/BCs will be presented in the following subsection. 

\subsection{Elasticity theory and mixed-variable formulation}\label{elas_theory}
In this subsection, we present the PINN in the context of solving elastodynamics problems. The governing equations for elastodynamics can be written in the tensor form, as follows:
\begin{equation} 
\label{elstody1} 
\bm{\nabla \cdot \sigma} + \mathbf{F}=\rho \mathbf{u}_{tt}
\end{equation}
\begin{equation} 
\label{elstody2} 
\bm{\epsilon } =\frac{1}{2}\left [ \mathbf{\nabla u} + \mathbf{{(\nabla u)}}^T \right ]
\end{equation}
\begin{equation} 
\label{elstody3} 
\bm{\sigma}=\mathbf{C}:\bm{\epsilon}
\end{equation}
\noindent where $\bm{\sigma}$ is the Cauchy stress tensor, $\bm{\epsilon}$ is the strain tensor, $\mathbf{u}$ is the displacement vector, $\mathbf{F}$ is the body force vector, $\mathbf{C}$ is the fourth-order constitutive tensor and $\bm{\nabla}$ is the Nabla operator. This set of equations are subject to certain I/BCs, e.g., defined in Eqs. \eref{boundarycd} and \eref{initialcd}. 

To establish a PINN framework for solving elastodynamics problems, the input and output of the network must be specified first. A straightforward means might be using the spatiotemporal coordinates $\mathbf{X}=(\mathbf{x},~t)$ as the input while the displacement field $\mathbf{Y}=\mathbf{u}$ as the output. However, inspired by the hybrid finite element, we propose the mixed-variable formulation, i.e., displacement and stress fields $\mathbf{Y}=(\mathbf{u},~\bm{\sigma})$ as the DNN output in this work. This formulation is found to be superior to the displacement formulation with regard the trainability of the network. A detailed comparison between these two types of formulation is presented in the \ref{mix_vs_disp}.

\begin{figure}[t!]
\centering
\subfigure[]{
\begin{minipage}[h]{0.90\linewidth}
\centering
\includegraphics[width=\linewidth]{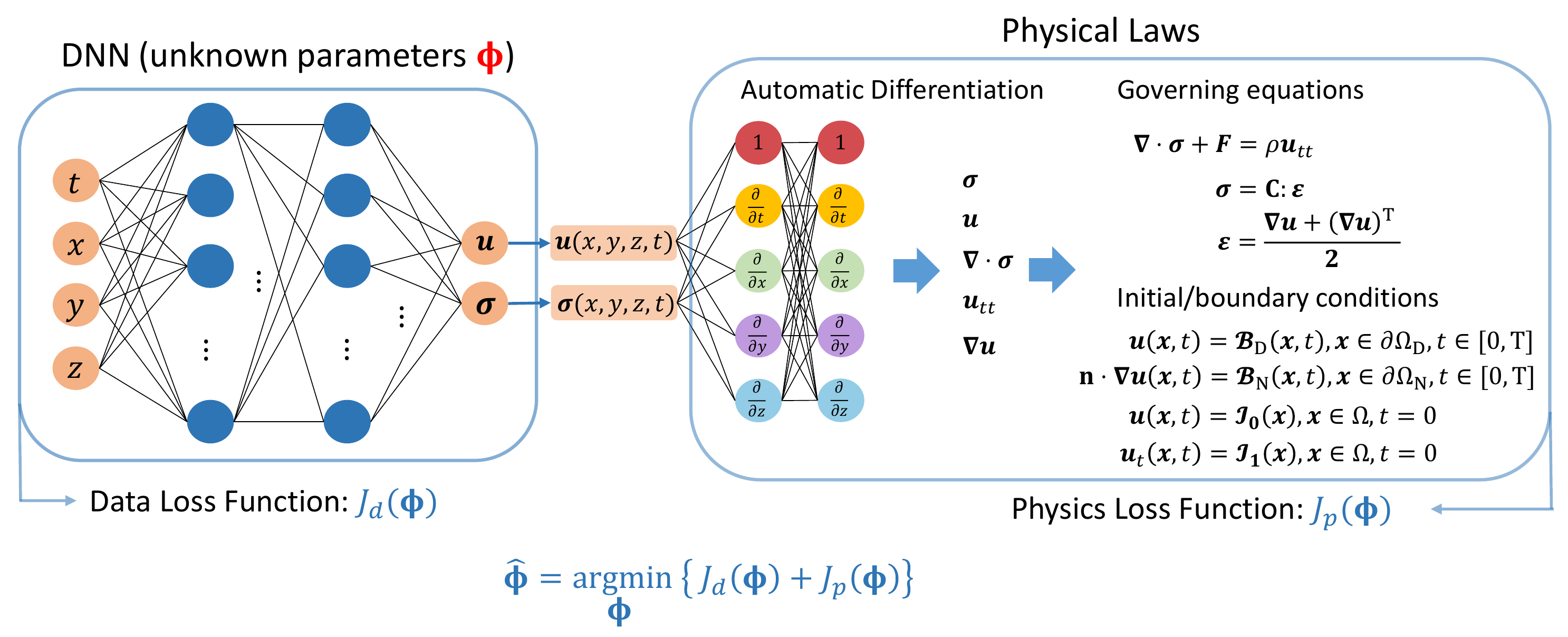}
\end{minipage}%
}%
\hfill
\subfigure[]{
\begin{minipage}[h]{0.99\linewidth}
\centering
\includegraphics[width=\linewidth]{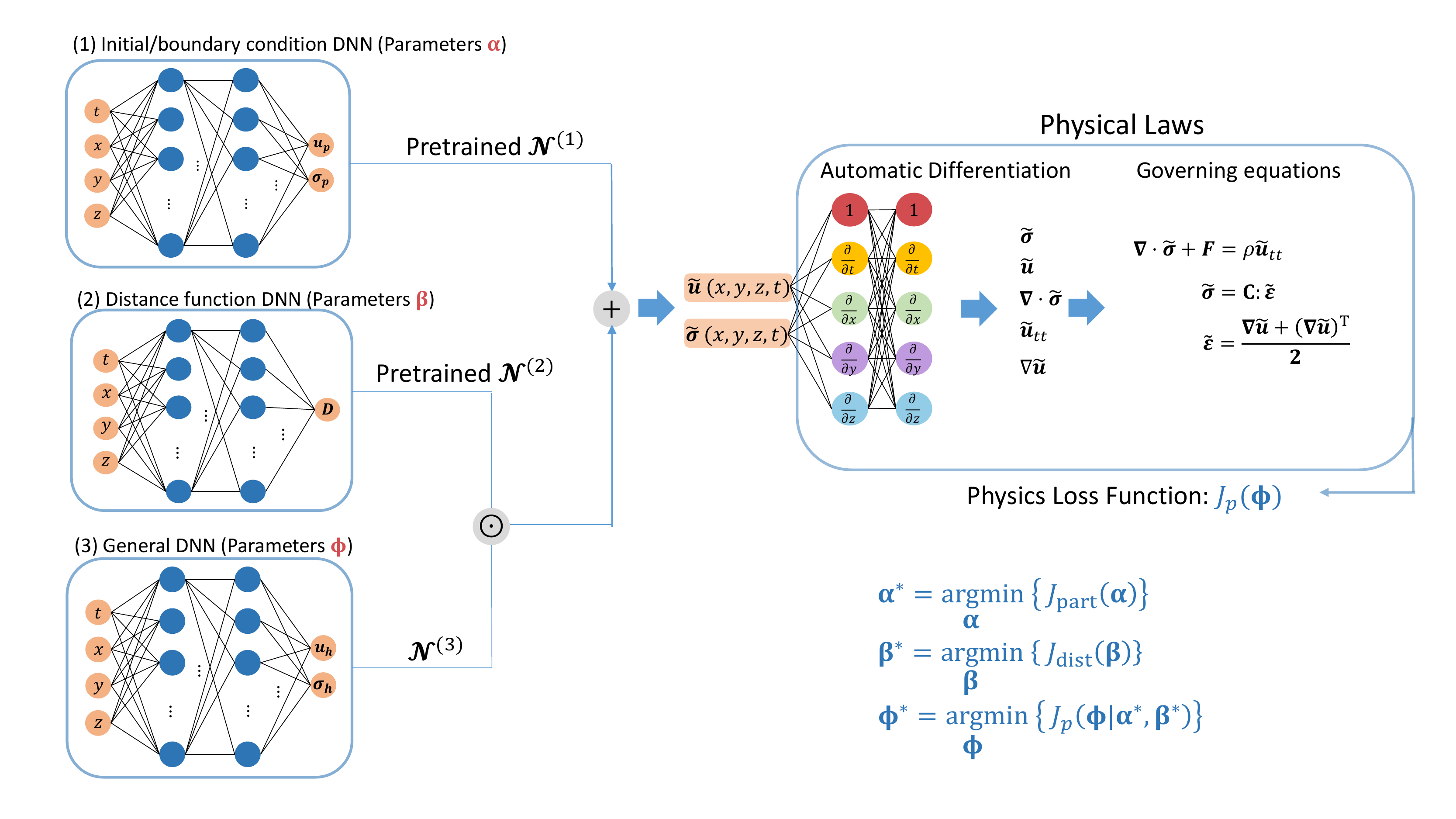}
\end{minipage}%
}%
\caption{Schematics  of PINN: (a) classical scheme with softly enforced I/BCs, and (b) composite scheme with hardly enforced I/BCs.}
\label{sec2_PINN_diagram}
\end{figure}

Another problem to be addressed is how to enforce the I/BCs. The imposition of I/BCs is crucial for solving the PDEs since it allows a unique solution. Considering the optimization nature of the PINN, the primitive way of applying I/BCs is to penalize the residual loss function of the PDE by the residuals of initial and boundary values via Lagrange multipliers in a ``soft'' manner. This strategy has been widely used in existing PINN methods \cite{weinan2018deep, sirignano2018dgm, raissi2019physics, jin2020nsfnets}. In this case, the physics loss function $J_p$ can be constructed with three components (i.e., equation loss $J_g$, initial value loss $J_{ic}$ and boundary value loss $J_{bc}$), given by
\begin{equation} 
\label{softloss} 
J_p=J_\text{g}+\lambda_1J_\text{bc}+\lambda_2J_\text{ic}
\end{equation}
where $\lambda_1>$0 and $\lambda_2>$0 are the relative weighting coefficients. To be more specific, each component of the total physics loss in Eq. \eref{softloss} is defined as
\begin{equation} 
\label{msef} 
J_\text{g}=||\bm{\nabla \cdot {\sigma}} +\mathbf{F}-\rho\mathbf{u}_{tt}||^2_{\Omega\times \left [ 0,T \right ]}+||\bm{{\sigma}-\mathbf{C}:{\epsilon}}||^2_{\Omega \times \left [ 0,T \right ]}
\end{equation}
\begin{equation} 
\label{msebc} 
J_\text{bc}=||\mathbf{u}-\mathcal{B}_D||^2_{\partial \Omega_D\times \left [ 0,T \right ]}+||\mathbf{n}\cdot\bm{\sigma}-\mathcal{B}_N||^2_{\partial \Omega_N\times \left [ 0,T \right ]}
\end{equation}
\begin{equation} 
\label{mseic} 
J_\text{ic}=||\mathbf{u}-\mathcal{I}_0||^2_{\Omega\times \left \{ t=0 \right \}}+||\mathbf{u}_t-\mathcal{I}_1||^2_{ \Omega\times\left \{ t=0 \right \}}
\end{equation}
where $||\cdot||^2$ denotes the mean square error (MSE) on the set of collocation points annotated by the corresponding subscript. For instance, the term $||\mathbf{u}-\mathcal{B}_D||^2_{\partial \Omega_D\times \left [ 0,T \right ]}$ represents the MSE evaluated at the collocation points on Dirichlet boundary (i.e. ${\partial \Omega_D\times \left [ 0,T \right ]}$). The physical quantities $\mathbf{u}$ and $\bm{\sigma}$ are obtained from the output of the DNN. This type of I/BC imposition is also called ``soft'' enforcement \cite{marquez2017imposing, sun2020surrogate} because the initial and boundary values may not be enforced accurately, due to the pathology issue of the gradient. A detailed study on this issue can be found in \cite{wang2020understanding}. To mitigate the issue of inaccurate I/BC enforcement, a trail-and-error procedure is usually involved to find suitable weighting coefficients $\lambda$'s. The framework of PINN developed with softly enforced I/BCs is depicted in Fig. \ref{sec2_PINN_diagram}(a). Noteworthy, having the measurement data makes the PINN modeling data-driven, which is however not a prerequisite.

To address this problem, a counterpart of ``hard'' I/BC enforcement, as shown in Fig. \ref{sec2_PINN_diagram}(b), is proposed by using a composite scheme, which consists of three single DNNs that represent the I/BC (or particular) network $\mathcal{N}^{(1)}$, the distance function network $\mathcal{N}^{(2)}$ and the general solution network $\mathcal{N}^{(3)}$. The final solution to the elastodynamics problem is thus constructed as follows
\begin{equation} 
\label{sec2_hard_eqn} 
\left\{\begin{matrix}
u(\mathbf{x},t)=\mathcal{N}_{u}^{(1)}(\mathbf{x},t)+\mathcal{N}_{u}^{(2)}(\mathbf{x},t)\cdot\mathcal{N}_{u}^{(3)}(\mathbf{x},t)\\ 
v(\mathbf{x},t)=\mathcal{N}_{v}^{(1)}(\mathbf{x},t)+\mathcal{N}_{v}^{(2)}(\mathbf{x},t)\cdot\mathcal{N}_{v}^{(3)}(\mathbf{x},t)\\ 
\sigma_{11}(\mathbf{x},t)=\mathcal{N}_{\sigma_{11}}^{(1)}(\mathbf{x},t)+\mathcal{N}_{\sigma_{11}}^{(2)}(\mathbf{x},t)\cdot\mathcal{N}_{\sigma_{11}}^{(3)}(\mathbf{x},t)\\ 
\sigma_{22}(\mathbf{x},t)=\mathcal{N}_{\sigma_{22}}^{(1)}(\mathbf{x},t)+\mathcal{N}_{\sigma_{22}}^{(2)}(\mathbf{x},t)\cdot\mathcal{N}_{\sigma_{22}}^{(3)}(\mathbf{x},t)\\ 
\sigma_{12}(\mathbf{x},t)=\mathcal{N}_{\sigma_{12}}^{(1)}(\mathbf{x},t)+\mathcal{N}_{\sigma_{12}}^{(2)}(\mathbf{x},t)\cdot\mathcal{N}_{\sigma_{12}}^{(3)}(\mathbf{x},t)
\end{matrix}\right.
\end{equation}
The goal of this practice is to guarantee that the constructed solution satisfies the I/BC automatically. To explain how it is achieved, let us discuss each component within the synergy solution in Eq. \eref{sec2_hard_eqn} . The first part $\mathcal{N}^{(1)}$ represents the ``particular'' solution network pre-trained only with the I/BC values, which are known beforehand. The network $\mathcal{N}^{(2)}$ represents the distance from a given point to the initial time or boundaries in the spatiotemporal space. It equals to zero at the initial time and boundaries, and nonzero within the domain. We can pre-train it with the generated distance function value (see Eq. \eref{d_component}) within the spatiotemporal domain. The general solution network $\mathcal{N}^{(3)}$ is the only adjustable part after pre-training of $\mathcal{N}^{(1)}$ and $\mathcal{N}^{(2)}$. We should note that it is called ``general'' because the only constraint for this network is to satisfy the governing equations. The beauty of the constructed solution is that, at the initial time or boundaries where the distance function evaluates to zero, the solution degrades to the particular solution so that the I/BCs will be imposed forcibly. This type of ``hard'' I/BC enforcement strategy can be generalized as follows 
\begin{equation} 
\label{hard} 
\mathcal{U}(\mathbf{x},t)=\mathcal{U}_p(\mathbf{x},t)+ \mathcal{D}(\mathbf{x},t)\odot\mathcal{U}_h(\mathbf{x},t)
\end{equation}
where $\mathcal{U}_{(\cdot)}$ denotes the physical quantity of interest , $\mathcal{D}$ the distance function and $\odot$ denotes the element-wise multiplication which will be omitted for simplicity in the remainder of this paper. 

When the geometry of the computation domain is relatively simple, both of the ``particular'' solution $\mathcal{U}_p(\mathbf{x},t)$ and distance function $\mathcal{D}(\mathbf{x},t)$ can be expressed analytically (see \cite{sun2020surrogate, samaniego2020energy} for example). However, to render PINN the capability of handling complex geometries, two low-capacity auxiliary neural networks are used to approximate the $\mathcal{U}_p(\mathbf{x},t)$ and $\mathcal{D}(\mathbf{x},t)$ in this paper, as shown by Fig. \ref{sec2_PINN_diagram}(b). After $\mathcal{U}_p(\mathbf{x},t)$ and $\mathcal{D}(\mathbf{x},t)$ are pre-trained and fixed, the composite DNNs will be finally trained as a whole to make sure the remaining constraint, residual of governing PDEs, is satisfied. In the final training of the composite PINN, only the weights and biases of $\mathcal{U}_h(\mathbf{x},t)$ will be the trainable variables exposed to the optimizer. 

To train the distance function $\mathcal{D}(\mathbf{x},t)$, we can sample the points $\{(\mathbf{x}_i,~t_i)\}_{i=1}^n$ in the domain $\Omega\times \left [ 0,T \right ]$ and compute the distance $\hat{\mathcal{D}}$ to the spatiotemporal boundaries, shown as follows
\begin{equation} 
\label{d_component} 
\left\{\begin{matrix}
\hat{\mathcal{D}}_u(\mathbf{x},t)=\text{min(distance to the spatiotemporal boundary of $u$)}\\ 
\hat{\mathcal{D}}_v(\mathbf{x},t)=\text{min(distance to the spatiotemporal boundary of $v$)}\\ 
\hat{\mathcal{D}}_{\sigma_{11}}(\mathbf{x},t)=\text{min(distance to the spatiotemporal boundary of $\sigma_{11}$)}\\ 
\hat{\mathcal{D}}_{\sigma_{22}}(\mathbf{x},t)=\text{min(distance to the spatiotemporal boundary of $\sigma_{22}$)}\\ 
\hat{\mathcal{D}}_{\sigma_{12}}(\mathbf{x},t)=\text{min(distance to the spatiotemporal boundary of $\sigma_{12}$)}
\end{matrix}\right.
\end{equation}
where the spatiotemporal boundary is defined as the combination of initial and boundary conditions for a specific solution variable. Here, for the two dimensional problem, the distance function $\mathcal{D}$ has five components which correspond to the output $\mathbf{Y}=(u,v,\sigma_{11},\sigma_{22},\sigma_{12})$. 

With the above composite network, we are able to ensure the satisfaction of the boundary conditions ($\mathcal{B}_D$ and $\mathcal{B}_N$) and the initial displacement condition ($\mathcal{I}_0$). To enforce the initial velocity condition ($\mathcal{I}_1$), extra constraints must be imposed on the $\mathcal{U}_p$ and $\mathcal{D}$. To illustrate this concept, let us assume the initial displacement and velocity of the domain to be zero. Enforcement of $\dot{\mathcal{U}}_p(\mathbf{x},t=0)$ and $\mathcal{D}(\mathbf{x},t=0)$ equal to zero will not guarantee the $\dot{\mathcal{U}}(\mathbf{x},t)$ to be zero since it contains the $\dot{\mathcal{D}}(\mathbf{x},t)$ term, as indicated in Eq. \eref{timederivative}. Therefore, we need to constrain $\dot{\mathcal{D}}(\mathbf{x},t=0)$ to be zero.
\begin{equation} 
\label{timederivative} 
\dot{\mathcal{U}}(\mathbf{x},t)=\dot{\mathcal{U}}_p(\mathbf{x},t)+ \dot{\mathcal{D}}(\mathbf{x},t)\mathcal{U}_h(\mathbf{x},t)+\mathcal{D}(\mathbf{x},t)\dot{\mathcal{U}}_h(\mathbf{x},t)
\end{equation}

Let us summarize the constraints we need to enforce on the $\mathcal{D}(\mathbf{x},t)$ and $\mathcal{U}_p(\mathbf{x},t)$ for a IBVP problem defined in Eqs. \eref{pde}--\eref{initialcd} so that we can formulate the corresponding loss functions. To train the $\mathcal{D}(\mathbf{x},t)$ network, the following conditions should be enforced
\begin{equation} 
\label{d_constraint_1} 
\mathcal{D}(\mathbf{x},t)=\left\{\begin{matrix}
\text{zero},  &\text{for }(\mathbf{x},t)\in(\partial \Omega\times \left [ 0,T \right ])\cup (\Omega\times \left \{ t=0 \right \})\\ 
\text{nonzero}, &\text{otherwise}
\end{matrix}\right.
\end{equation}
\begin{equation} 
\label{d_constraint_2} 
\dot{\mathcal{D}}(\mathbf{x},t)=0,~\text{for } (\mathbf{x},t)\in(\Omega\times \left \{ t=0 \right \})
\end{equation}
Meanwhile, to train the $\mathcal{U}_p(\mathbf{x},t)$ network, we have the following constraints
\begin{equation} 
\label{up_constraint_1} 
\mathcal{U}_p(\mathbf{x},t)=\left\{\begin{matrix}
&\mathcal{B}_D,  &\text{for }(\mathbf{x},t)\in\partial \Omega_D\times \left [ 0,T \right ]\\ 
&\mathcal{I}_0, &\text{for }(\mathbf{x},t)\in \Omega\times \{ t=0 \}
\end{matrix}\right.
\end{equation}
\begin{equation} 
\label{up_constraint_2} 
\mathbf{n}\cdot\nabla \mathcal{U}_p(\mathbf{x},t)=\mathcal{B}_N(t),~ \text{for }(\mathbf{x},t)\in\partial \Omega_N\times \left [ 0,T \right ]
\end{equation}
\begin{equation} 
\label{up_constraint_3} 
\dot{\mathcal{U}}_p(\mathbf{x},t)=\mathcal{I}_1,~\text{for }(\mathbf{x},t)\in \Omega\times \{ t=0 \}
\end{equation}
Therefore, the loss functions for the $\mathcal{N}^{(1)}$ and $\mathcal{N}^{(2)}$ networks shown in Fig. \ref{sec2_PINN_diagram}(b) are written as
\begin{equation} 
\label{train_part} 
J_\text{part}=||\mathbf{n}\cdot\nabla \mathcal{U}_p - \mathcal{B}_N||^2_{\partial \Omega_N\times \left [ 0,T \right ]}
+||\mathcal{U}_p-\mathcal{B}_D||^2_{\partial \Omega_D\times \left [ 0,T \right ]} \\
+||\mathcal{U}_p-\mathcal{I}_0||^2_{\Omega\times \{ t=0 \}}
+||\dot{\mathcal{U}}_p-\mathcal{I}_1||^2_{\Omega\times \{ t=0 \}}
\end{equation}
\begin{equation} 
\label{train_dist} 
J_\text{dist}=||\mathcal{D} - \hat{\mathcal{D}}||^2_{\Omega\times \left [ 0,T \right ]}
+||\dot{\mathcal{D}}||^2_{\Omega\times \left \{ t=0 \right \}}
\end{equation}
while the only remaining term in the $J_p$ would be the governing equation loss $J_\text{g}$ as shown in Eq. \eref{msef}. 

A simple static example is provided in \ref{soft_vs_hard} to compare the performance between the proposed ··hard‘’ enforcement of I/BCs and the conventional ``soft'' enforcement. The accuracy of the boundary value shows an improvement by the proposed PINN scheme over the conventional approach. The source code for each numerical example in this paper can be found in {\color{blue}https://github.com/Raocp/PINN-elastodynamics} upon publication.


\section{Results}\label{results}
\subsection{Defected plate under periodic uni-axial tension}\label{notchplate}
\begin{figure}[b!]
	\centering
	\includegraphics[width=0.4\textwidth]{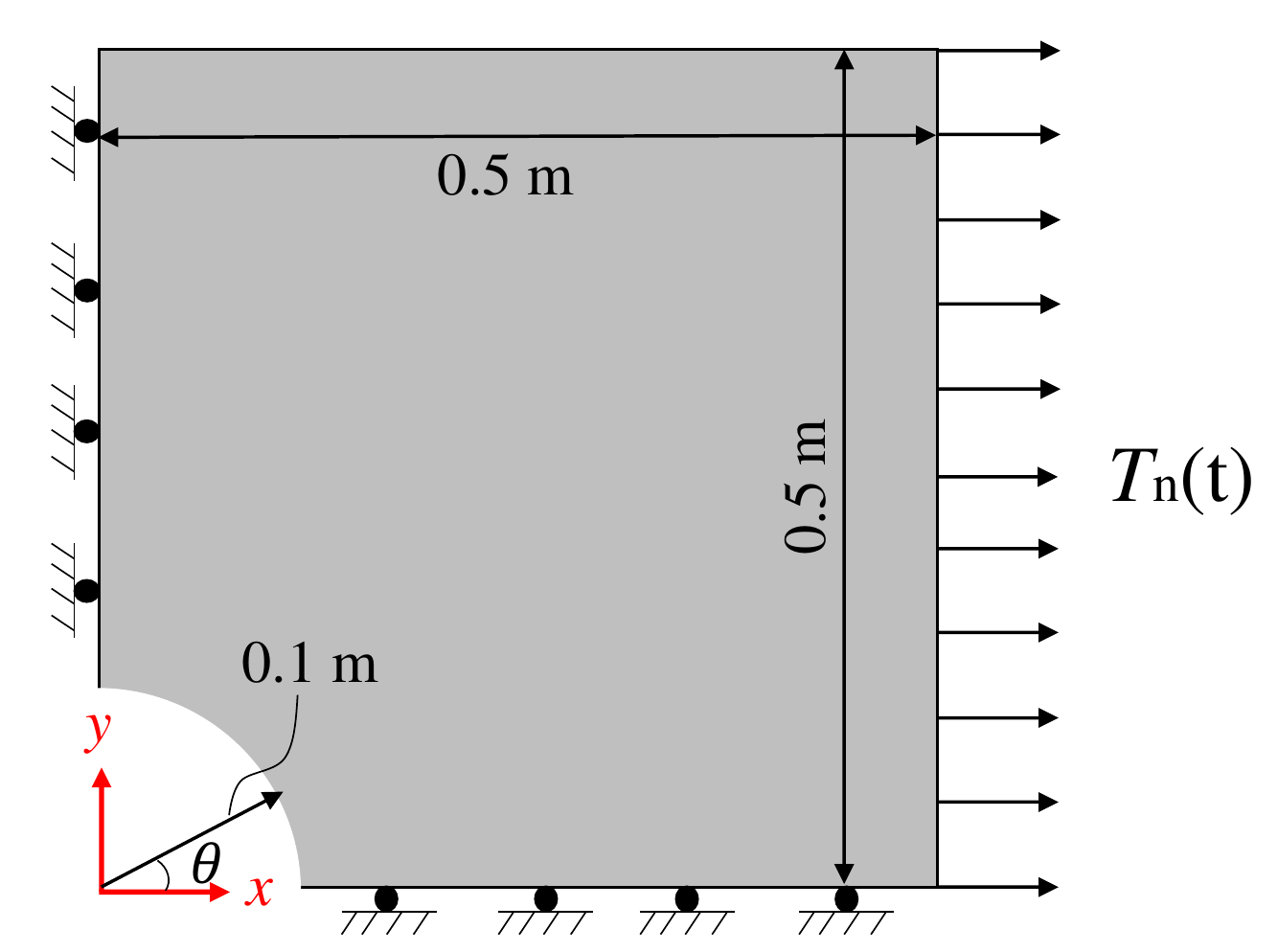}
	\caption{Diagram of defected plate under uni-axial load.}
	\label{sec3-2-sketch}
\end{figure}

A two-dimensional plane stress problem, i.e., a defected plate under uni-axial tension, is considered in this example. The total length of the square plate is 1.0 m while the radius of the circular defection located in the center is 0.1 m. Due to the symmetry of the problem, only a quarter plate is simulated (see Fig. \ref{sec3-2-sketch}). The Young's modulus and Poisson's ratio of the plate are 20 MPa and 0.25, respectively. A uni-axial normal traction $T_n(t)$ is applied on the right edge as shown in Fig. \ref{sec3-2-sketch}. The I/BCs are imposed in a ``hard'' manner, as introduced in the previous section. It is noted that we impose the traction free condition of the hole surface as an extra equation, altogether with the governing equations, which is represented in matrix form as 
\begin{equation}\label{tractionfree}
\begin{bmatrix}
T_x\\ T_y
\end{bmatrix}
=
\begin{bmatrix}
\sigma_{xx} &\sigma_{xy} \\ 
\sigma_{yx} &\sigma_{yy} 
\end{bmatrix} 
\begin{bmatrix}
n_x\\ n_y
\end{bmatrix} =
\begin{bmatrix}
0\\ 0
\end{bmatrix} 
\end{equation}
where $[n_x, n_y]^\text{T}$ is the unit normal vector of the surface. This system of equations standalone will not reveal any boundary value on the hole surface since the number of unknowns (3) is greater than that of equations (2).

\begin{table}[t!]
  \caption{Results of the convergence test with regard to the width and depth of the $\mathcal{U}_h$ network. Networks of $\mathcal{D}$ ($4\times10$) and $\mathcal{U}_p$ ($2\times5$) are fixed while training the final synergy solution $\mathcal{U}=\mathcal{U}_p+ \mathcal{D}\odot\mathcal{U}_h$. The FE solution is used as reference to compute the relative $\ell_2$ errors of von Mises stress and displacement vector. \vspace{6pt}}
  \label{convtest}
  \small
  \begin{tabularx}{\textwidth}{p{0.15\textwidth} p{0.15\textwidth} p{0.15\textwidth} p{0.15\textwidth} p{0.15\textwidth} p{0.15\textwidth}}
    \hline\hline
    {Depth$\times$width} & {$J_{\text{dist}}$ ($\mathcal{D}$)} & {$J_{\text{part}}$ ($\mathcal{U}_p$)} & {$J_p$ ($\mathcal{U}_h$)} & {$\mathcal{E}$($\mathbf{u}$)} & {$\mathcal{E}$($\sigma_v$)}\\
    \midrule
    4$\times$30 & \multirow{6}{*}{$1.1\times10^{-6}$} &  \multirow{6}{*}{$9.4\times10^{-8}$} &  $1.4\times10^{-4} $ & $4.1\times10^{-2}$ & $2.6\times10^{-2}$\\
    4$\times$40 &  &  &  $9.8\times10^{-5}$ & $4.0\times10^{-2}$ & $2.6\times10^{-2}$\\
    5$\times$40 &  &  &  $7.7\times10^{-5}$ & $3.7\times10^{-2}$ & $2.3\times10^{-2}$\\
    5$\times$50 &  &  & $4.3\times10^{-5}$ & $3.1\times10^{-2}$ & $1.5\times10^{-2}$\\
    6$\times$50 &  &  & $1.4\times10^{-5}$ & $2.0\times10^{-2}$ & $4.1\times10^{-3}$\\
    6$\times$60 &  &  & $1.3\times10^{-5}$ & $1.9\times10^{-2}$ & $3.4\times10^{-3}$\\
    \hline\hline
  \end{tabularx}
\end{table}
\normalsize

\begin{figure}[t!]
\centering
\subfigure[$\sigma_{11}$]{
\begin{minipage}[t]{0.32\linewidth}
\centering
\includegraphics[width=\linewidth]{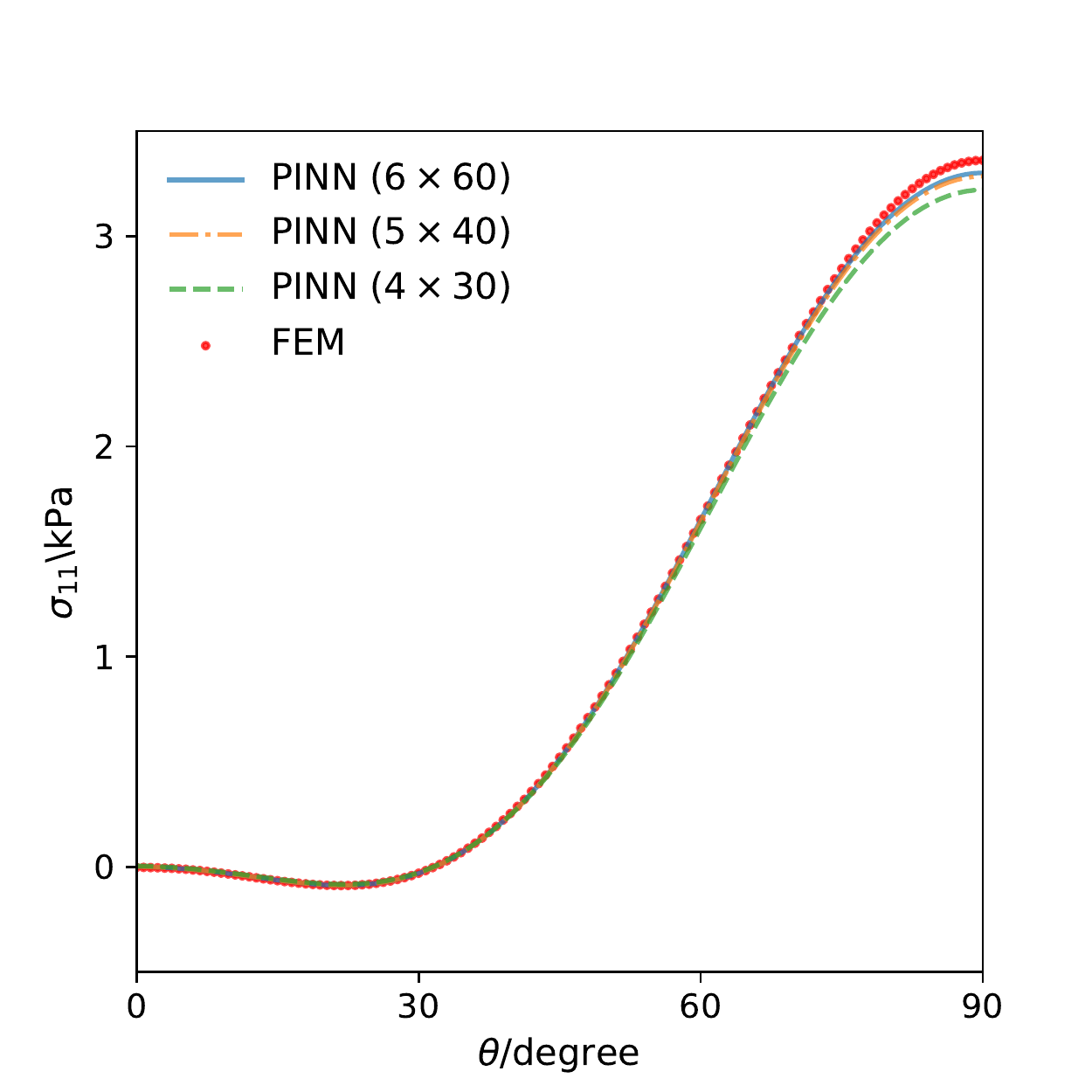}
\end{minipage}%
}%
\subfigure[$\sigma_{22}$]{
\begin{minipage}[t]{0.32\linewidth}
\centering
\includegraphics[width=\linewidth]{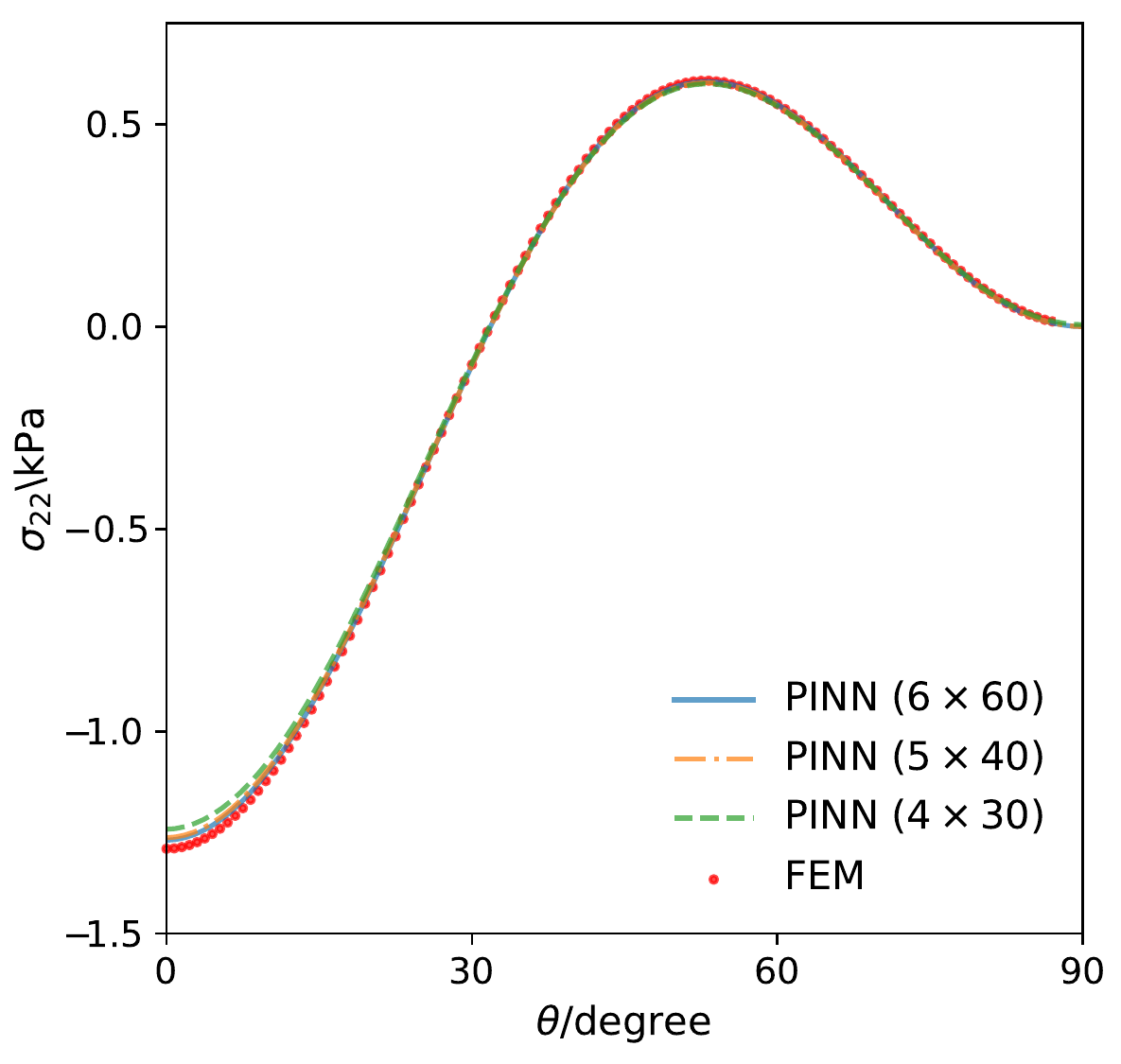}
\end{minipage}%
}%
\hfill
\centering
\subfigure[$\sigma_{12}$]{
\begin{minipage}[t]{0.32\linewidth}
\centering
\includegraphics[width=\linewidth]{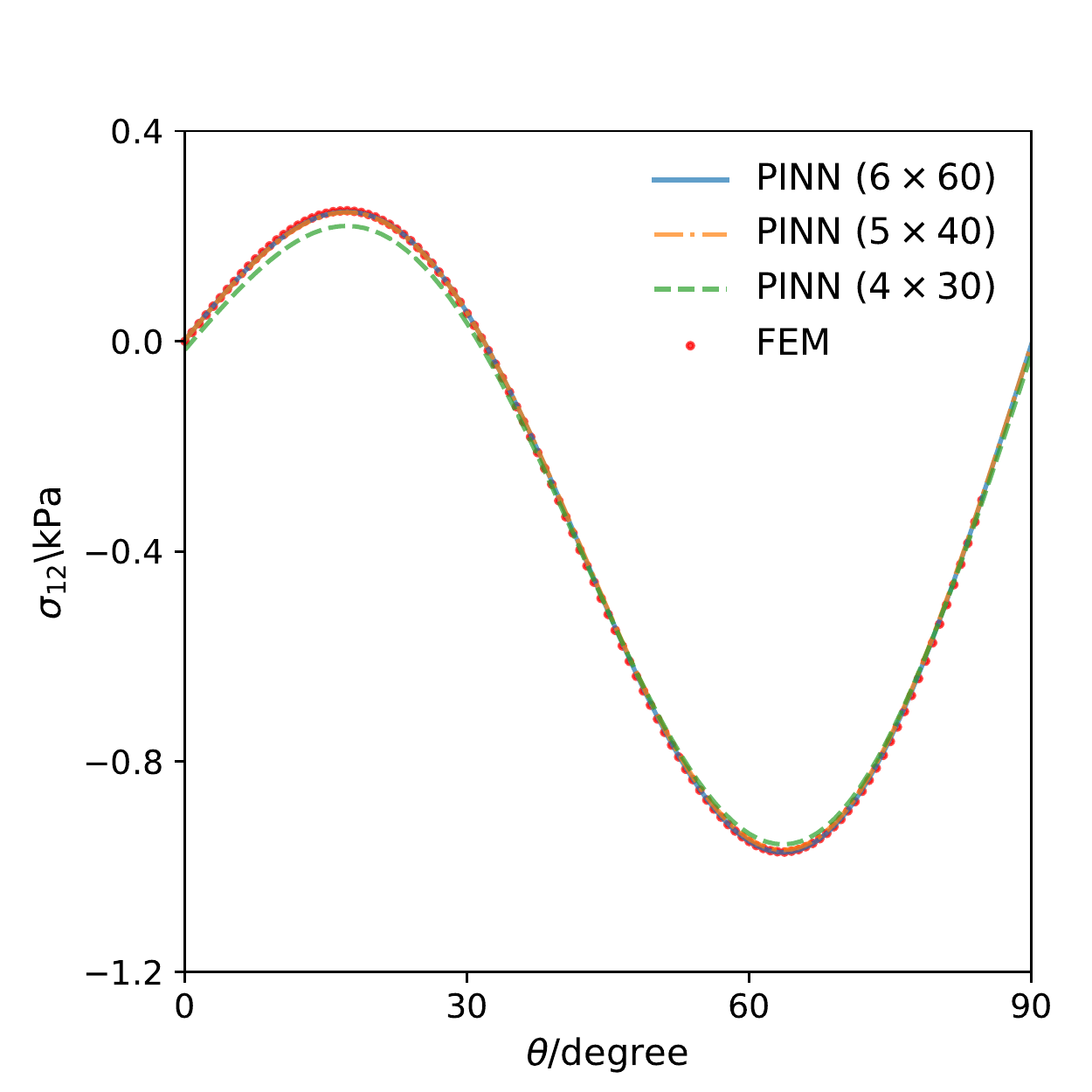}
\end{minipage}
}%
\caption{Convergence analysis of the stress distribution on notched surface.}\label{converanaly}
\end{figure}

\begin{figure}[b!]
	\centering
	\includegraphics[width=0.6\textwidth]{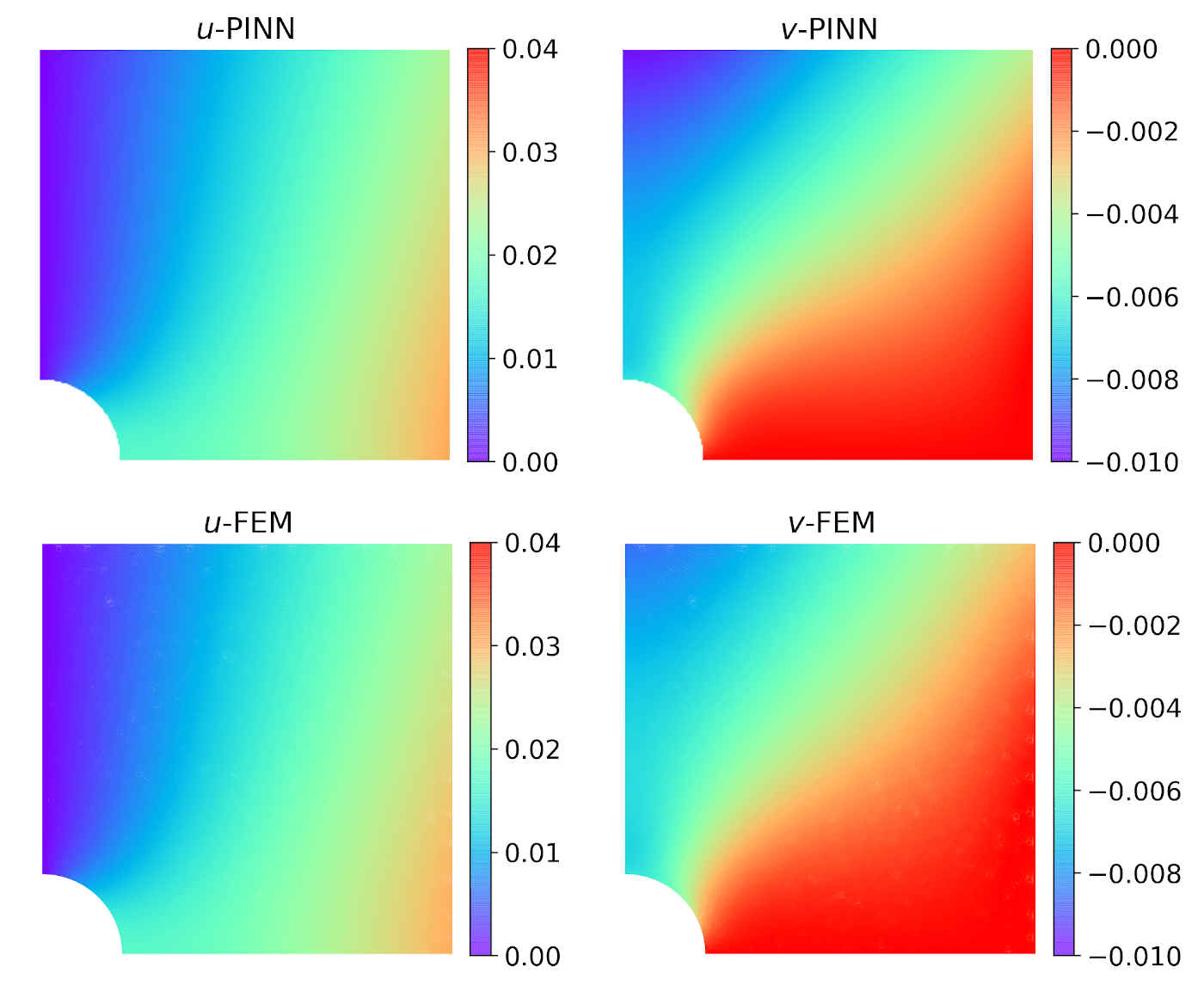}
	\caption{Comparison of the obtained displacement fields (top: current with $6\times60$ net; bottom: FEM).}
	\label{staticdispfield}
\end{figure}

\begin{figure}[t!]
	\centering
	\includegraphics[width=0.9\textwidth]{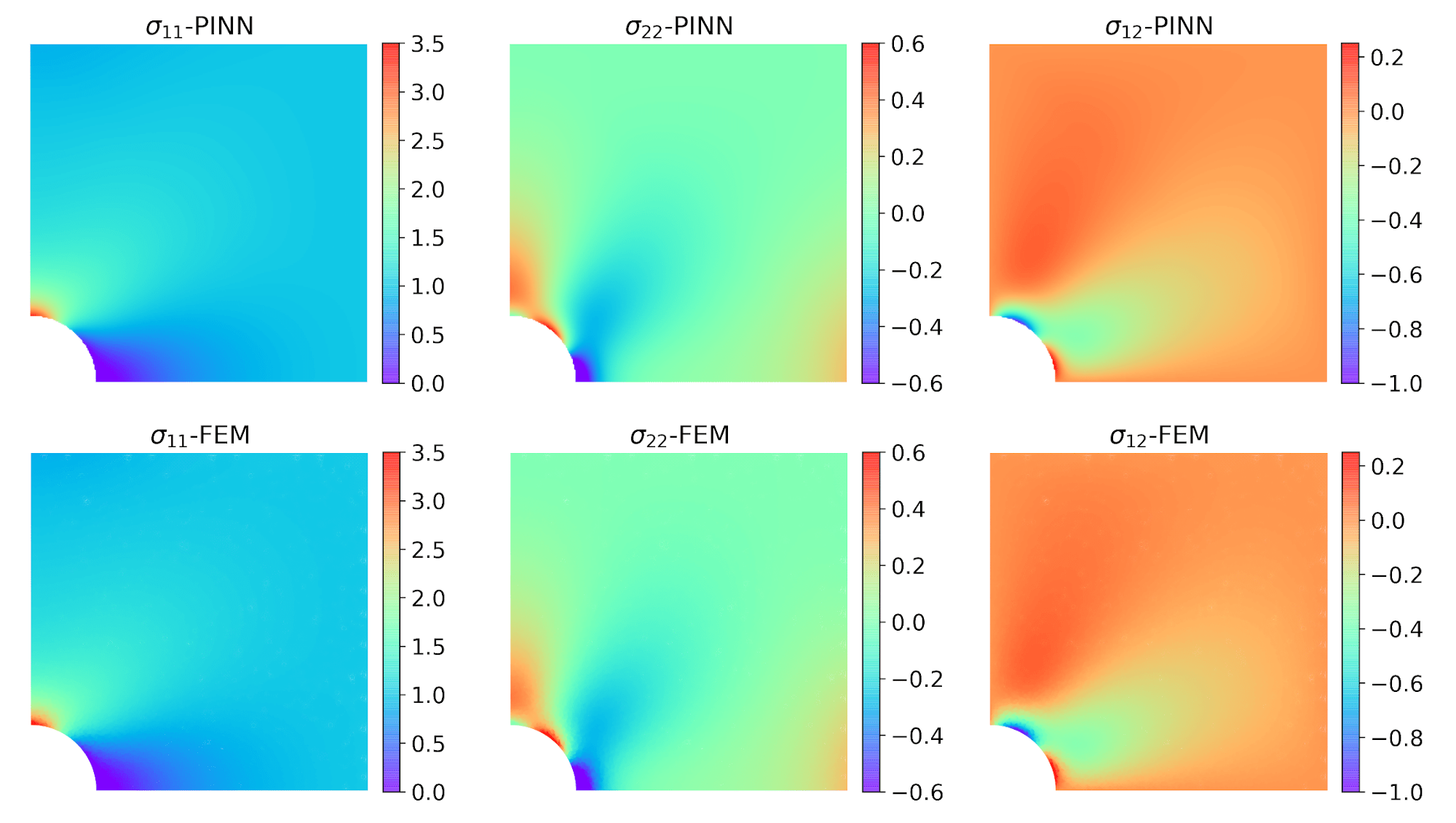}
	\caption{Comparison of the obtained stress fields (top: current with $6\times60$ net; bottom: FEM).}
	\label{staticstressfield}
\end{figure}

To make PINN a general approach for solving PDEs, the obtained results should converge along with the increase of network complexity, which is primarily controlled by the width (number of neurons) and depth (number of layers) of the DNN. This convergence property of PINN is guaranteed by the universal approximation theorem \cite{cybenko1989approximation}. To verify this, we consider a static case in which a constant $T_n(t)=1.0$ MPa is applied on the edge. The selection of width and depth is usually done in a grid-search way. However, to reduce the number of trails, a step-search way is adopted herein, e.g., the width and depth of the network increase alternately. It should be noted that the convergence test is only conducted on the general solution network of $\mathcal{U}_h$, since the $\mathcal{D}$ and $\mathcal{U}_p$ can be easily trained with two low-capacity networks. Therefore, the architecture for $\mathcal{D}$ and $\mathcal{U}_p$ are kept the same for each case, whose depth $\times$ width settings are $4\times10$ and $2\times5$ respectively. The number of collocation points for evaluating the equation residuals and the traction-free boundary condition of the notch surface is 25,000 and 160. The equation residual points are refined near the notch to capture the stress concentration. To train the distance function $\mathcal{D}$, 400 Cartesian grid points within the domain and 60 points on the notch are employed. In addition, 200 collocation points on the symmetry boundaries and 400 points on the traction boundaries are used for training the particular solution $\mathcal{U}_p$. The tanh activation function and the Xavier initialization \cite{glorot2010understanding} for trainable parameters are employed throughout the paper. The combined Adam \cite{kingma2014adam} and limited-memory Broyden-Fletcher-Goldfarb-Shanno algorithm with bound constraints (L-BFGS-B) \cite{zhu1997algorithm} are employed as the optimization algorithm to enhance both global search and local tuning. The whole training process consists of three stages, namely, 2,000 iterations of the Adam optimizer with $10^{-3}$ learning rate (LR), 2,000 Adam iterations with LR = $5\times10^{-4}$, and the L-BFGS-B optimizer until the loss converges to a small tolerance. Table \ref{convtest} summarizes the configurations of the general solution network ($\mathcal{N}^{(3)}$) in each case, as well as the final loss and the relative $\ell_2$ errors defined by 
\begin{equation} 
\label{relativeL2E} 
\begin{aligned}
\mathcal{E}(\mathbf{f}) = \frac{\sqrt{\sum_{i=1}^{M}||\mathbf{f}_{\text{pred}}^i-\mathbf{f}_{\text{ref}}^i||^2}}{\sqrt{\sum_{i=1}^{M}||\mathbf{f}_{\text{ref}}^i||^2}}
\end{aligned}
\end{equation}
where $\mathbf{f}$ is the physical quantity of interest, and $M$ is the total number of reference points. The physical quantities used for comparison with the finite element (FE) solution are the von Mises stress and the displacement vector. It can be seen that as the network becomes deeper and wider, both the achieved loss value and the relative $\ell_2$ errors become smaller. Also, the comparison of the stress distribution on the notch surface between the FE solution and the PINN result is plotted in Fig. \ref{converanaly}. It can be observed that the results produced by PINN converge as the total number of neurons increases. The predicted displacement and stress fields by the proposed PINN approach are presented in Fig. \ref{staticdispfield} and \ref{staticstressfield} which show satisfactory agreement with the FE reference solution. In the following examples, the architecture of the network is directly given after the convergence test with regard to the width and depth.

Next we consider a dynamic case with cyclic traction $\displaystyle T_n(t)=0.5\text{sin}\left(2\pi t/{T_0} + 1.5\pi\right) + 0.5$ where $T_0=5$ s denotes the period of the load. The total duration of the simulation is 10 s. The whole plate is initially at rest, i.e., $\displaystyle \mathbf{u}_t$ and $\mathbf{u}$ equal to zero on the entire domain. In addition to the boundary conditions, these two initial conditions are also enforced in a ``hard'' way  (see Section \ref{elas_theory}). The networks with the architecture of $4\times20$, $4\times20$ and $8\times80$ are used to represent the $\mathcal{U}_p$, $\mathcal{D}$ and $\mathcal{U}_h$, respectively. A total number of 10,000 collocation points at initial time, 16,000 points on the symmetry boundaries and 16,000 on the traction boundaries are used for training $\mathcal{U}_p$. Meanwhile, 8,651 Cartesian grid points and 840 notch points are used for training the distance function $\mathcal{D}$. The number of collocation points for evaluating the equation residuals and the notch traction-free boundary condition are 120,000 and 9,600. The training setting and procedure are the same as those in the static case. 

\begin{figure}[t!]
\centering
\subfigure[$\sigma_{11}$]{
\begin{minipage}[t]{0.32\linewidth}
\centering
\includegraphics[width=\linewidth]{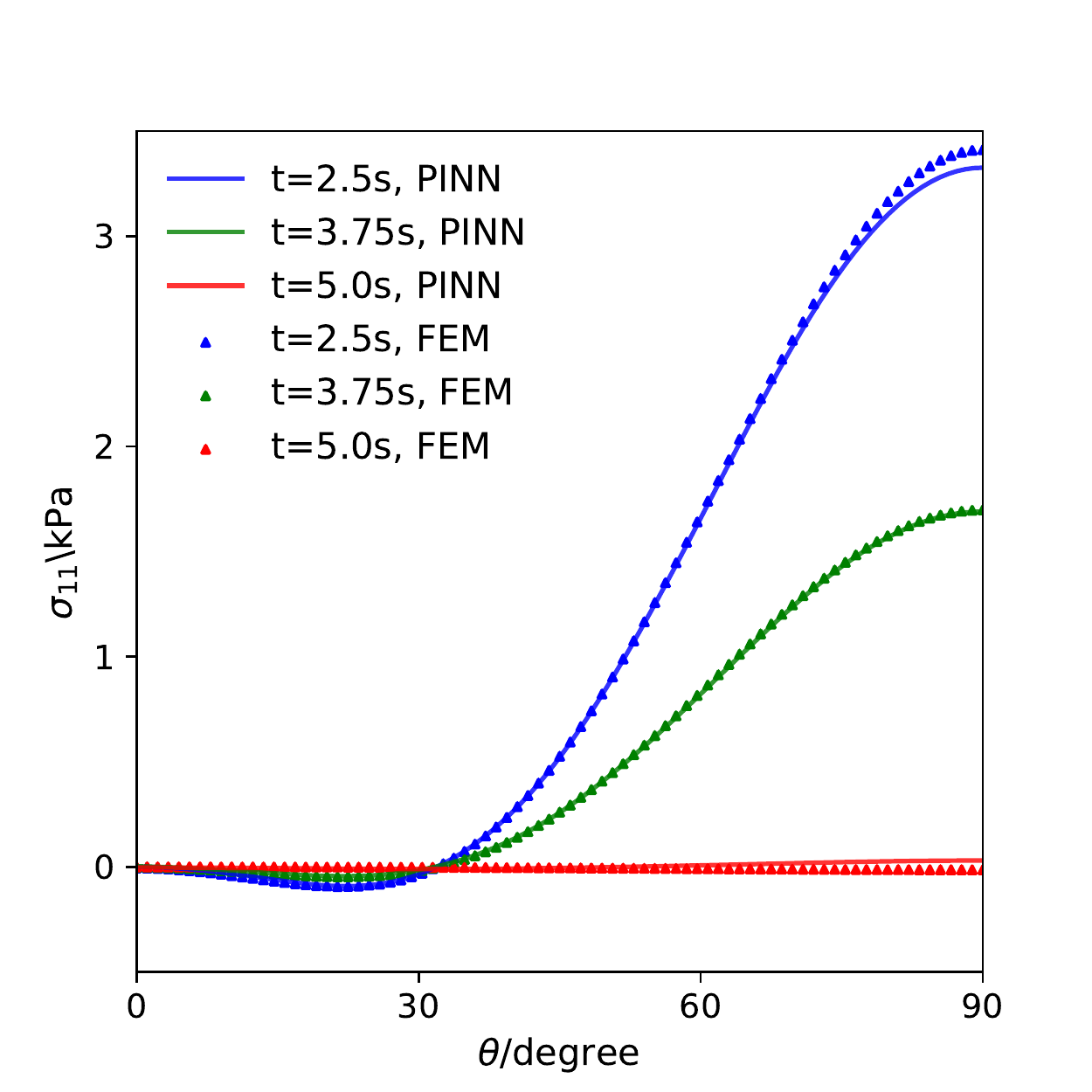}
\end{minipage}%
}%
\subfigure[$\sigma_{22}$]{
\begin{minipage}[t]{0.32\linewidth}
\centering
\includegraphics[width=\linewidth]{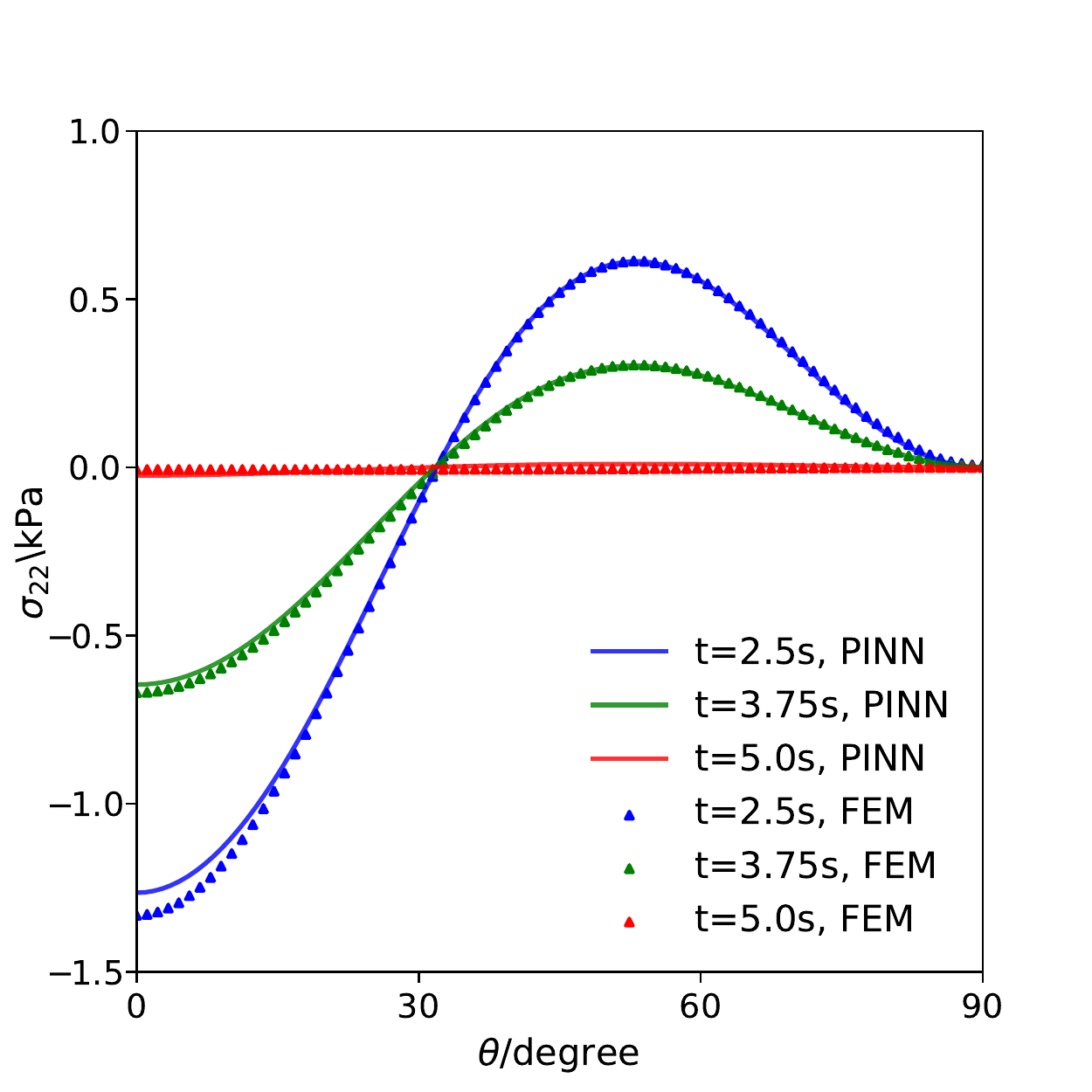}
\end{minipage}%
}%
\hfill
\centering
\subfigure[$\sigma_{12}$]{
\begin{minipage}[t]{0.32\linewidth}
\centering
\includegraphics[width=\linewidth]{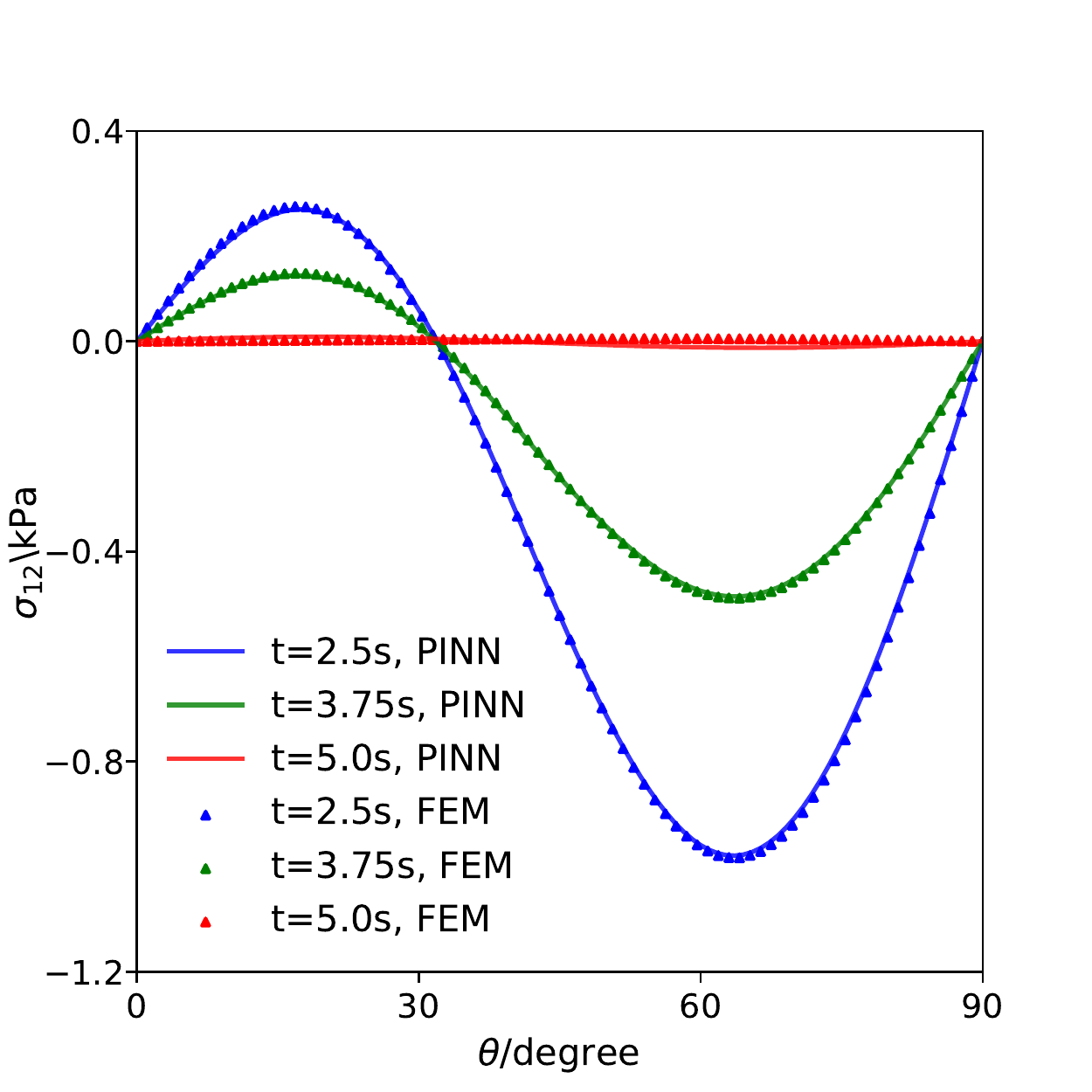}
\end{minipage}
}%
\caption{Comparison of the stress distribution on notched surface.}\label{dynstressdistr}
\end{figure}

\begin{figure}[t!]
	\centering
	\includegraphics[width=0.55\textwidth]{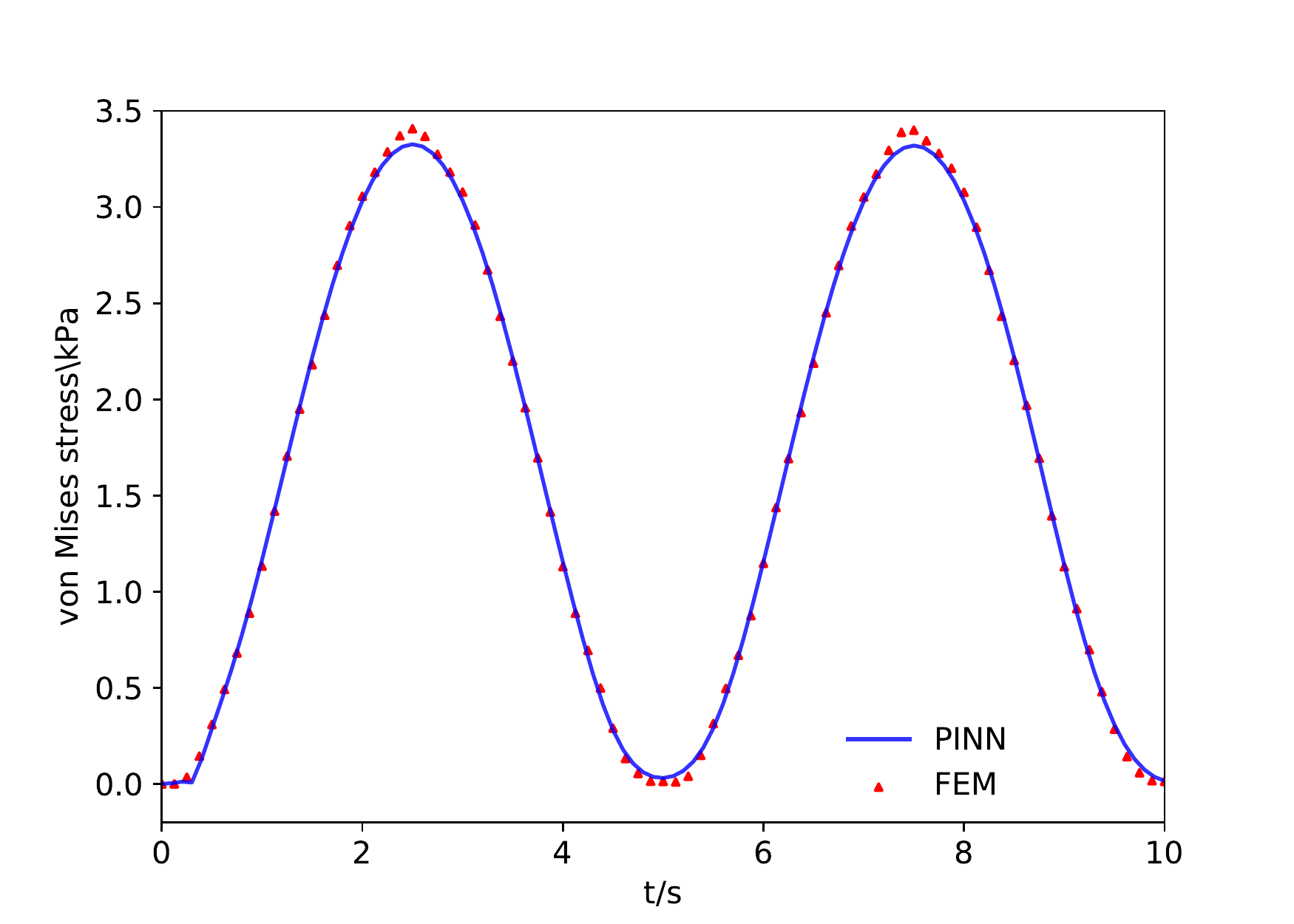}
	\caption{History of von Mises stress at point $(x,y)=(0,0.1)$.}
	\label{sec3-2-vonmises}
\end{figure}

The stress distributions, predicted by the proposed PINN, at various moments on the notch are shown in Fig. \ref{dynstressdistr} which are consistent with the reference solution obtained by the implicit finite element method (FEM). The von Mises stress on the point most vulnerable to the yielding at ($x, y$)=(0, 0.1) m is plotted in Fig. \ref{sec3-2-vonmises}, which demonstrates the capability of PINN for capturing the evolution of stress over time.

\subsection{Elastic wave propagation}\label{ela_wave_sec}
The performance of the proposed PINN approach is tested on the elastic wave propagation in a 2D homogeneous media in this section. As introduced in Section \ref{method}, PINN deals with the strong form of PDEs, while many of the other numerical methods, such as the Galerkin method, handle the weak form. As a result of the energy nature of the functional in the Galerkin method, the boundary without prescribed displacement or traction will be treated as free surface since it has no contribution to the total potential energy. As the PINN deals with the strong form governing equations directly, we are allowed to apply boundary conditions on only part of the domain boundaries. This owes to the nature of hyperbolic PDEs as further explained in Remark \ref{domain_dep}, which has been elaborated in \cite{evans2010partial}. 

\begin{remark}\label{domain_dep}
The domain of dependence of a hyperbolic PDE for a given point is the portion of the problem domain that influences the value of this point. That is to say, the solution at a generic point is unique as long as the solution in the domain of dependence is determined. The domain of dependence of a point is bounded by a surface (or hypersurface) in the spatiotemporal space.
\end{remark}

The problems considered herein are featured with only one circular wave source in the middle of the domain. Hence, the solution at a given point can be computed as long as the solution in the upstream domain is determined. This feature would be extremely useful for full wave-inversion problems in geophysics that involves infinite (or semi-infinite) domain. To reduce the infinite/semi-infinite domain to a truncated computation domain, researchers usually resort to the artificial absorbing layers (e.g., perfectly matched layer (PML) \cite{berenger1994perfectly}) or the enlargement of the computation domain to avoid the wave reflection issue. These treatments would result in complicated numerical implementation or unnecessary computational burden for uninterested regions. 

\begin{figure}[t!]
	\centering
	\includegraphics[width=0.32\textwidth]{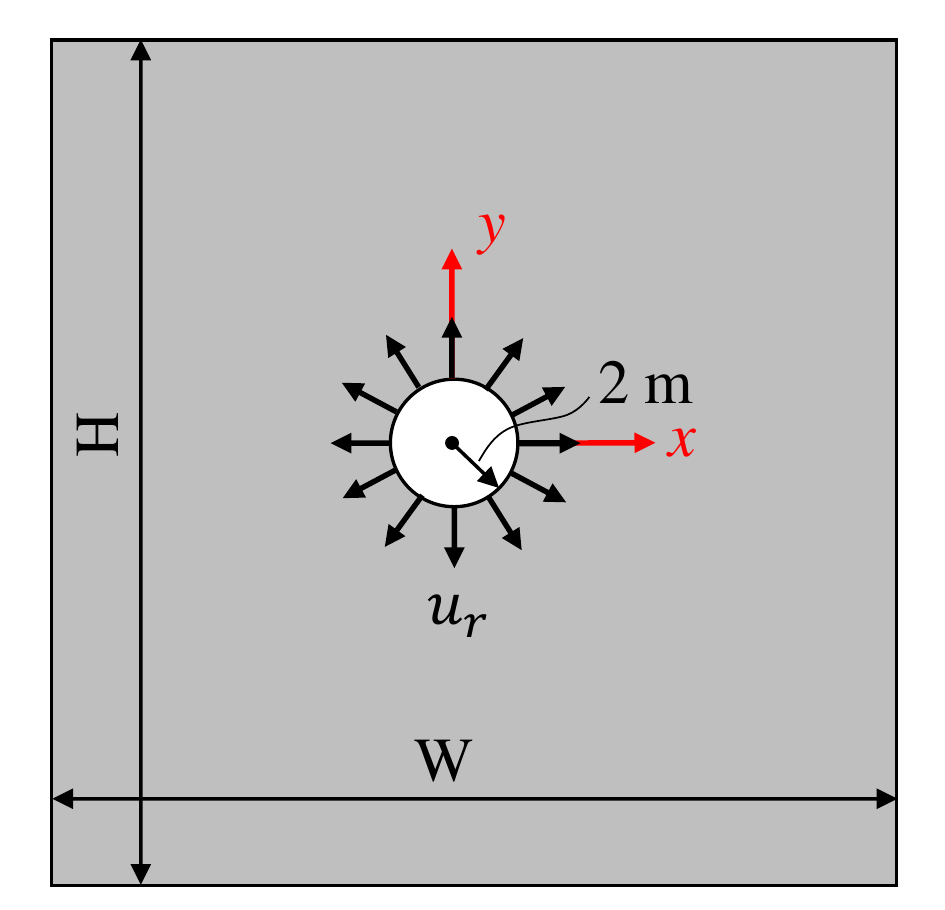}
	\caption{Diagram of the computation domain for elastic wave propagation.}
	\label{sec3-3-infi-domain}
\end{figure}

\begin{figure}[t!]
\centering
\subfigure[$t=$ 5 s]{
\begin{minipage}[t]{0.48\linewidth}
\centering
\includegraphics[width=\linewidth]{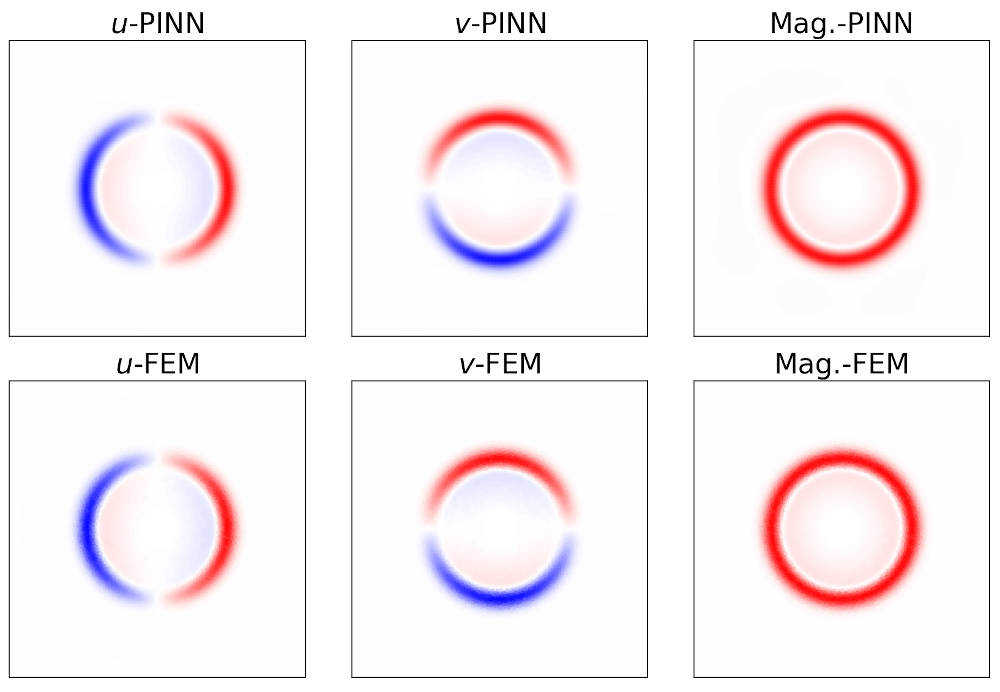}
\end{minipage}%
}%
\subfigure[$t=$ 8 s]{
\begin{minipage}[t]{0.48\linewidth}
\centering
\includegraphics[width=\linewidth]{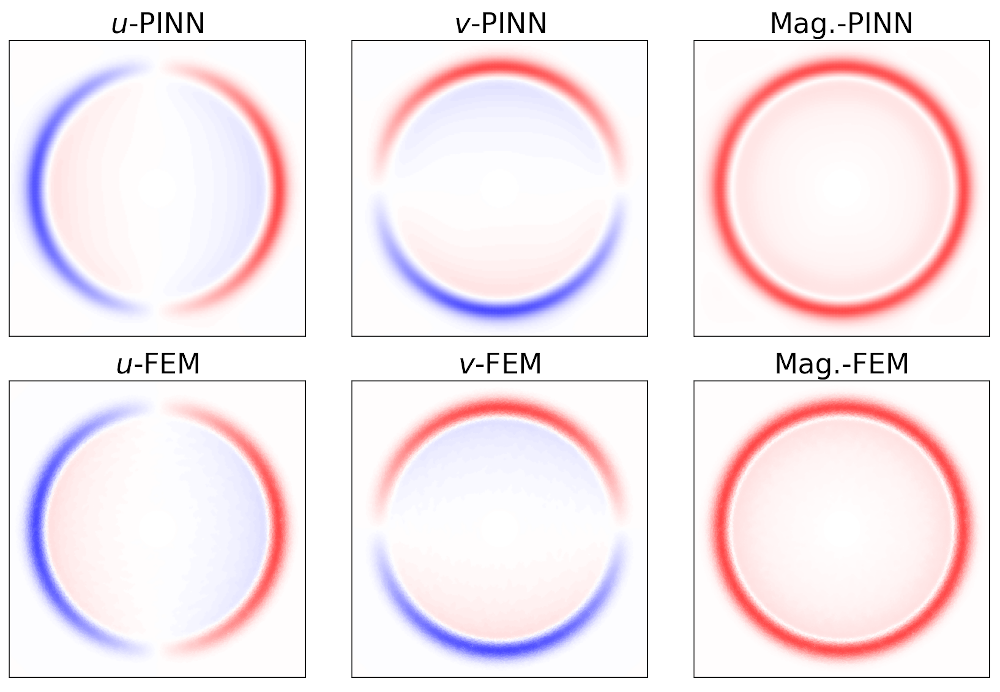}
\end{minipage}%
}%
\newline
\centering
\subfigure[$t=$ 11 s]{
\begin{minipage}[t]{0.48\linewidth}
\centering
\includegraphics[width=\linewidth]{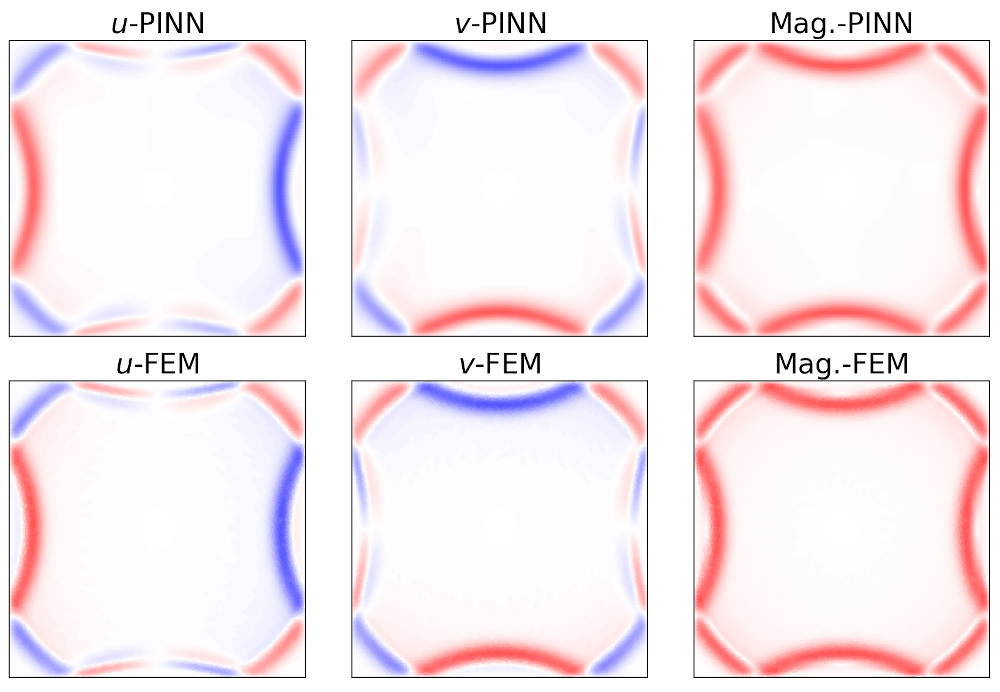}
\end{minipage}
}%
\centering
\subfigure[$t=$ 14 s]{
\begin{minipage}[t]{0.48\linewidth}
\centering
\includegraphics[width=\linewidth]{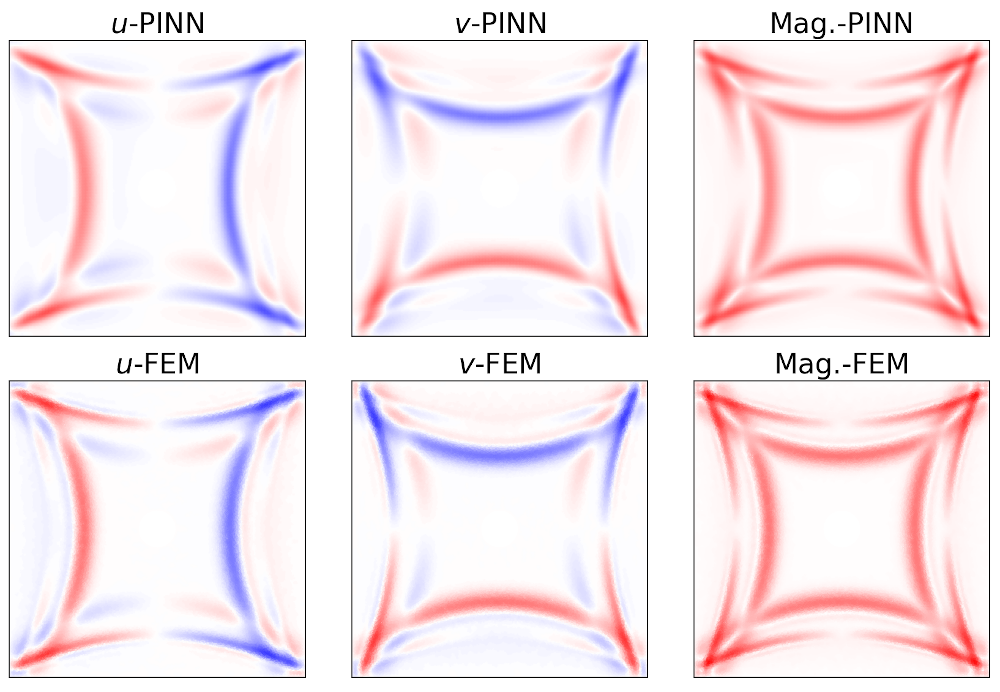}
\end{minipage}%
}%
\newline
\centering
\subfigure{
\begin{minipage}[t]{0.2\linewidth}
\centering
\includegraphics[width=\linewidth]{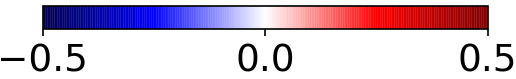}
\end{minipage}%
}%
\caption{Predicted displacement fields at various moments in confined domain. (upper) The implicit finite element solution with 41,696 linear quadrilateral elements. The time step is 0.01 s. (lower) The PINN solution. The architectures of $6\times140$, $3\times 30$ and $3\times 20$ are used in the networks for approximating $\mathcal{D}$, $\mathcal{U}_p$ and $\mathcal{U}_h$, respectively. I/BCs are enforced in a ``hard'' way. }\label{confined_disp}
\end{figure}

\begin{figure}[t!]
\centering
\subfigure[$t=$ 5 s]{
\begin{minipage}[t]{0.48\linewidth}
\centering
\includegraphics[width=\linewidth]{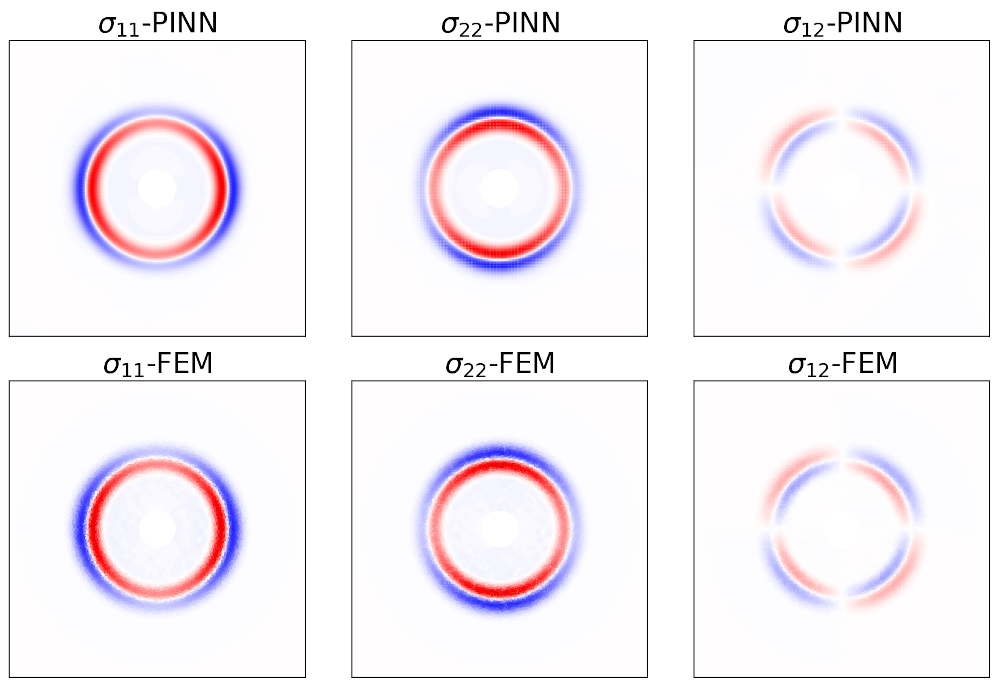}
\end{minipage}%
}%
\subfigure[$t=$ 8 s]{
\begin{minipage}[t]{0.48\linewidth}
\centering
\includegraphics[width=\linewidth]{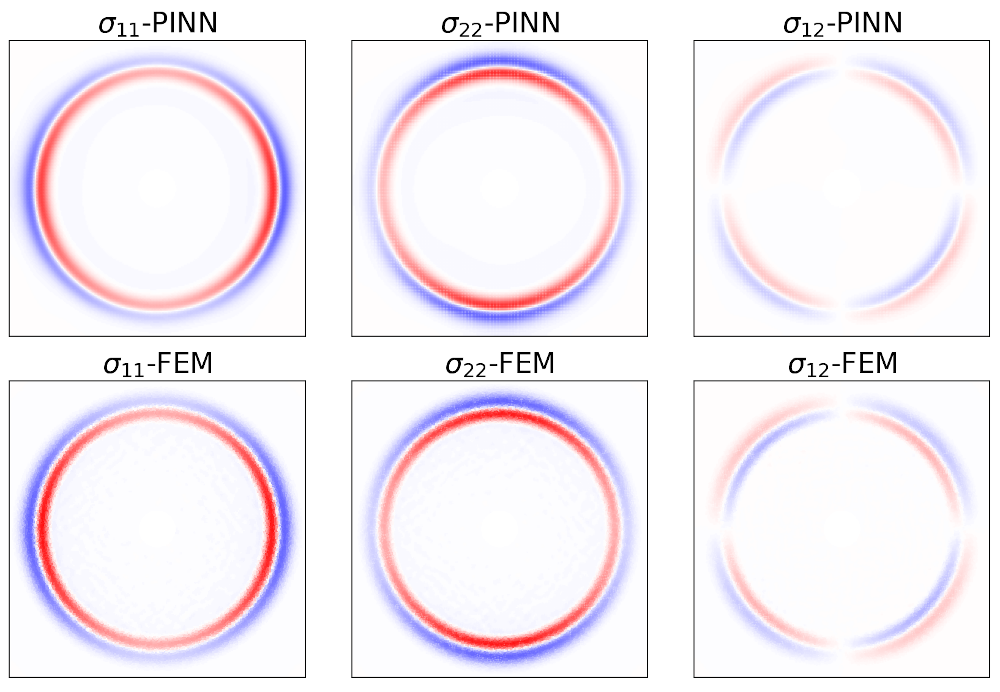}
\end{minipage}%
}%
\newline
\centering
\subfigure[$t=$ 11 s]{
\begin{minipage}[t]{0.48\linewidth}
\centering
\includegraphics[width=\linewidth]{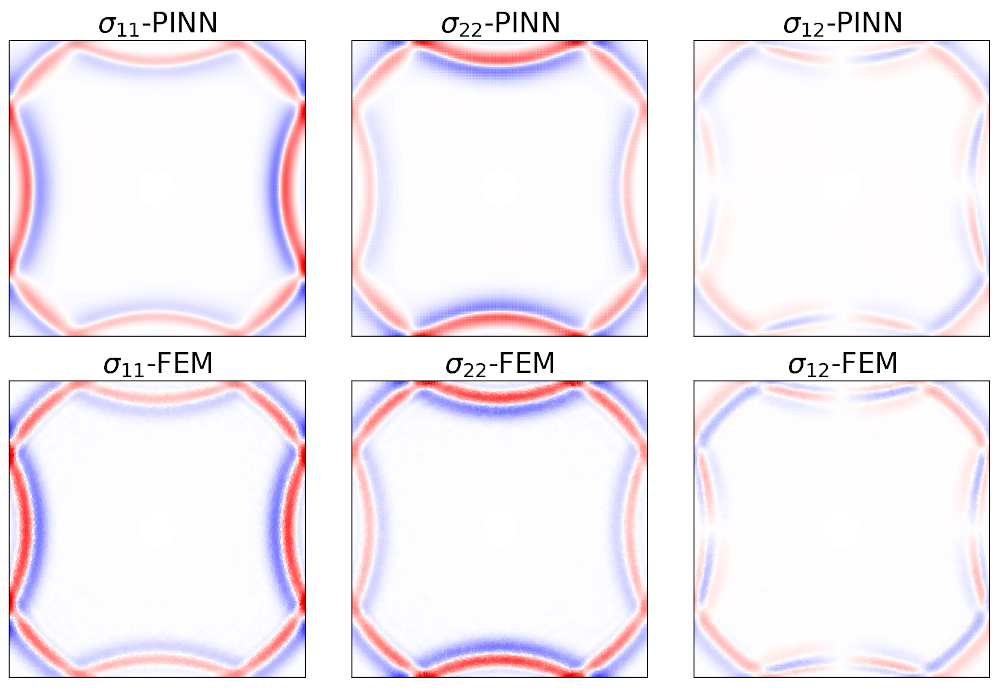}
\end{minipage}
}%
\centering
\subfigure[$t=$ 14 s]{
\begin{minipage}[t]{0.48\linewidth}
\centering
\includegraphics[width=\linewidth]{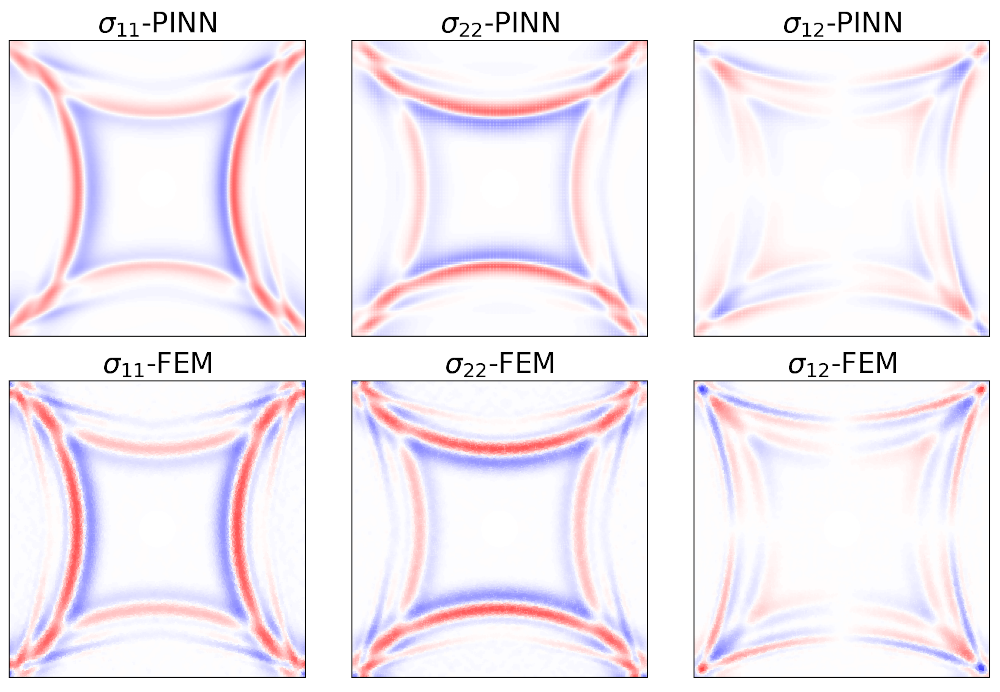}
\end{minipage}%
}%
\newline
\centering
\subfigure{
\begin{minipage}[t]{0.2\linewidth}
\centering
\includegraphics[width=\linewidth]{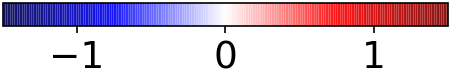}
\end{minipage}%
}%
\caption{Predicted stress fields at various moments in confined domain.}\label{confined_stress}
\end{figure}

We will employ the proposed PINN to simulate the wave propagation in homogeneous elastic media, within the region of interest considered (see Fig. \ref{sec3-3-infi-domain}). Three cases, with the same computation domain but different types of boundary conditions (i.e., completely confined, infinite and semi-infinite domains), are considered to examine the capability of the proposed approach. Besides, the domain is initially at rest for all the cases. The Young's modulus and Poisson's ratio of the media are 2.5 MPa and 0.25, respectively. The plane strain constitutive equations are used in the formulation of the loss function.

\subsubsection{Fully confined domain}
In the first case, we consider the wave problem in a confined domain (i.e., four fixed edges) whose boundary conditions will be hardly enforced. The wave source is prescribed by the radial displacement in the form of a Gaussian-like pulse, defined as
\begin{equation} 
\label{gausswave} 
u_r=U_0\exp\left[-\left(\frac{t-t_s}{t_p}\right)^2\right]
\end{equation}
where $U_0=0.5$ m represents the amplitude, $t_s=2.0$ s is the offset and $t_p=0.5$ s controls the width of the pulse. The distance function $\mathcal{D}$ represented by a $3\times 30$ network is trained with 3,840 Cartesian grid points and 1,000 notch surface points, while the particular solution $\mathcal{U}_p$ approximated by the $3\times 20$ network is trained with 6,000 collocation points at the initial state, 28,000 points on the fixed edges and 38,000 points on the wave source. For the general solution $\mathcal{U}_h$ approximated by a $6\times 140$ network, the total number of collocation points for evaluating the PDE residuals is 150,000, with a denser distribution near the wave source and four edges. The time duration of the simulation is 14 seconds. While training the general solution network, the parameters of $\mathcal{D}$ and $\mathcal{U}_p$ networks are fixed. The training consists of two stages: 10,000 epochs by the Adam optimizer with $LR=5\times10^{-4}$ followed by the L-BFGS-B optimizer for solution fine tuning. This training procedure is also adopted in the following two cases. The displacement and stress fields predicted by the proposed PINN, in comparision with the reference numerical solution, are presented in Fig. \ref{confined_disp} and \ref{confined_stress}, respectively. It can be seen that the reflection and interaction of the waves are accurately reproduced by the proposed PINN.

\begin{figure}[t!]
\centering
\subfigure[$t=$ 3 s]{
\begin{minipage}[t]{0.48\linewidth}
\centering
\includegraphics[width=\linewidth]{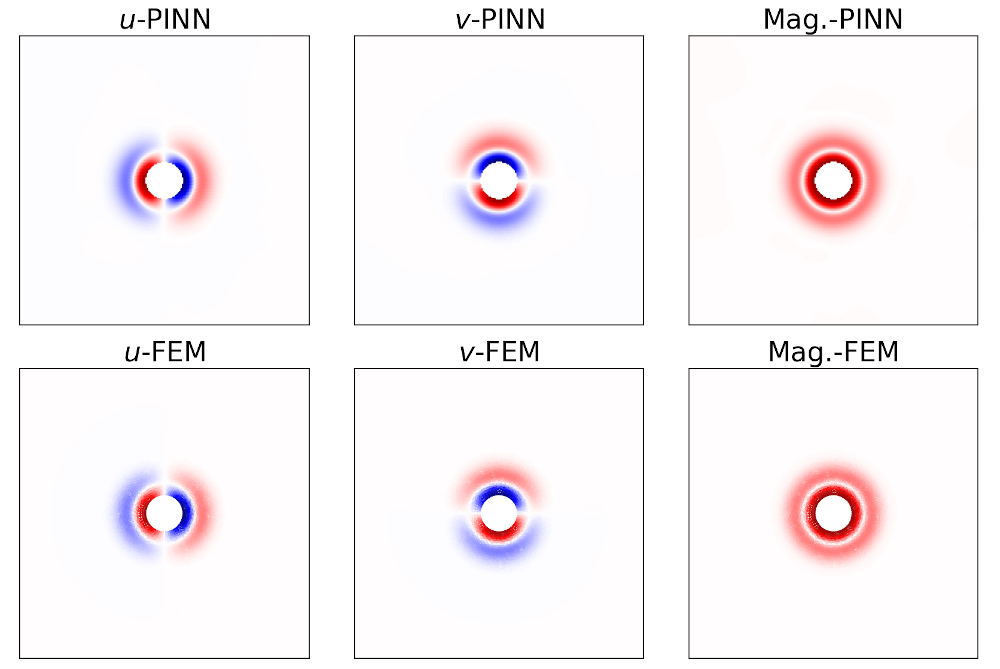}
\end{minipage}%
}%
\subfigure[$t=$ 6 s]{
\begin{minipage}[t]{0.48\linewidth}
\centering
\includegraphics[width=\linewidth]{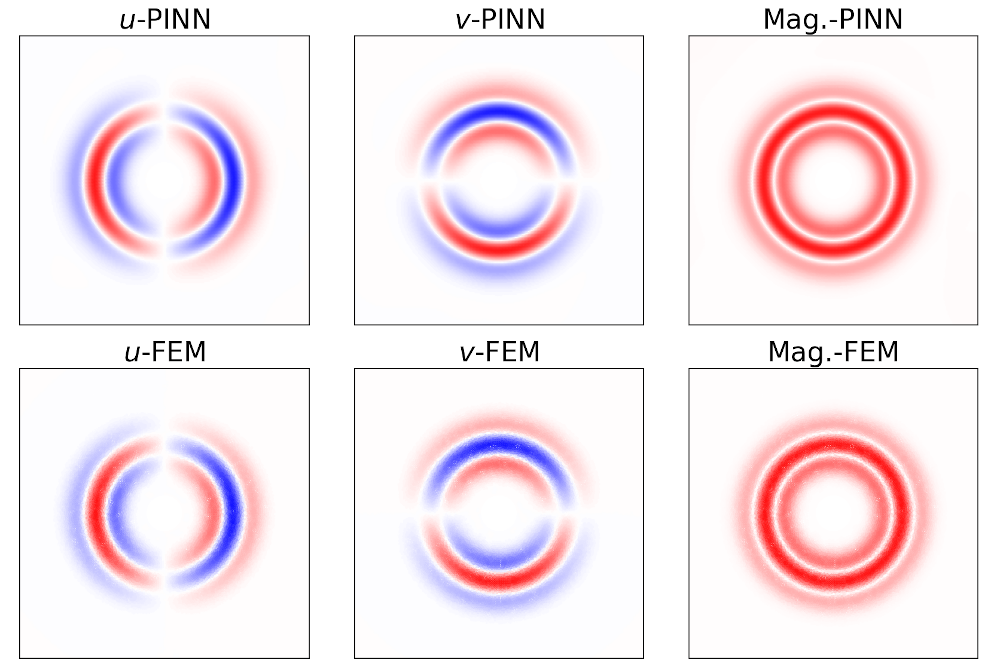}
\end{minipage}%
}%
\newline
\centering
\subfigure[$t=$ 9 s]{
\begin{minipage}[t]{0.48\linewidth}
\centering
\includegraphics[width=\linewidth]{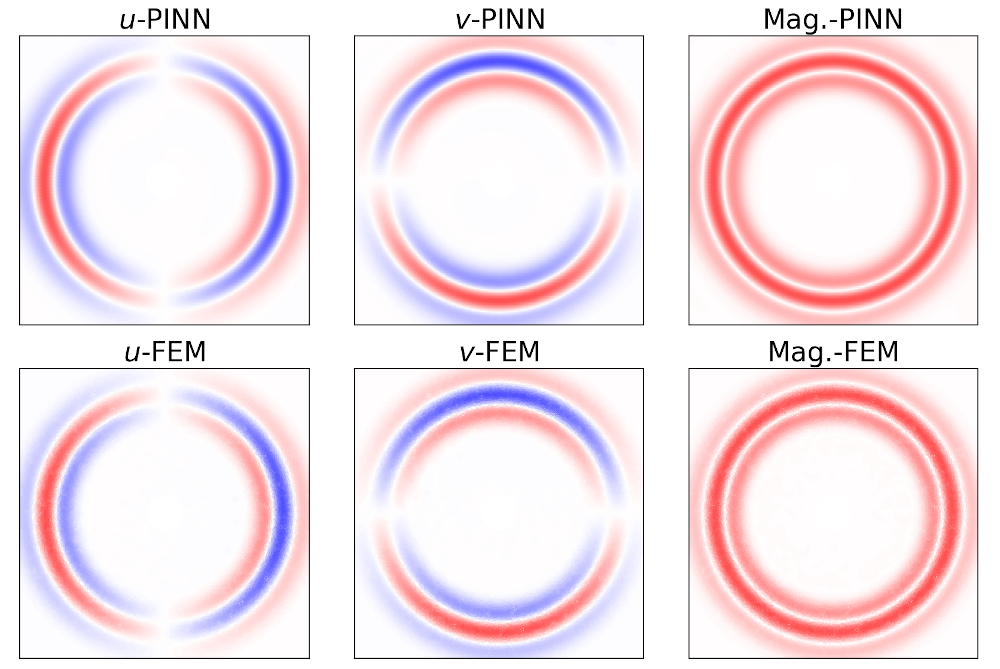}
\end{minipage}
}%
\centering
\subfigure[$t=$ 12 s]{
\begin{minipage}[t]{0.48\linewidth}
\centering
\includegraphics[width=\linewidth]{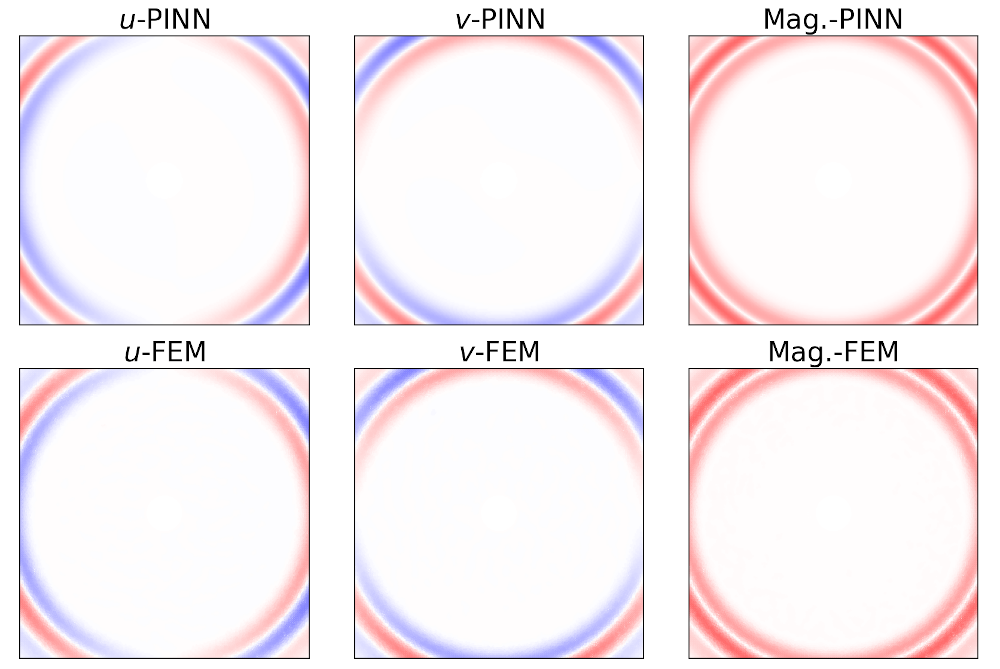}
\end{minipage}%
}%
\newline
\centering
\subfigure{
\begin{minipage}[t]{0.2\linewidth}
\centering
\includegraphics[width=\linewidth]{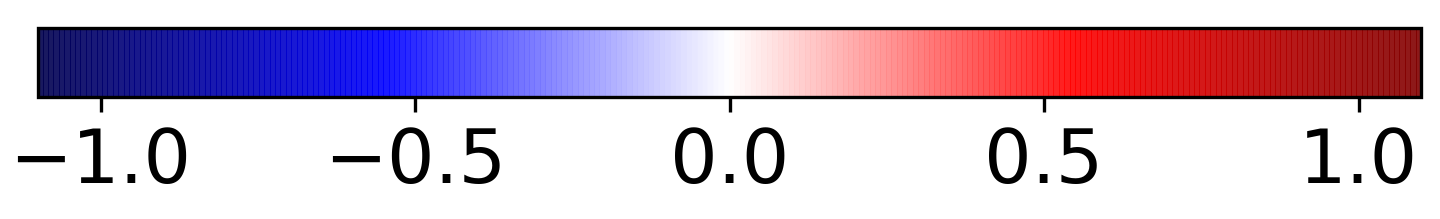}
\end{minipage}%
}%
\caption{Predicted displacement fields at various moments in the infinite domain. (upper) The finite element solution is computed from a enlarged domain (90 m by 90 m) with 61,294 linear quadrilateral elements. Time step equals to 0.01 s. (lower) The PINN solution. The $8\times80$ network with softly enforced I/BCs are employed. The first 10,000 epochs are trained by the Adam optimizer followed by the L-BFGS-B optimizer.}\label{inf_disp}
\end{figure}

\subsubsection{Infinite domain}
\begin{figure}[t!]
\centering
\subfigure[$t=$ 3 s]{
\begin{minipage}[t]{0.48\linewidth}
\centering
\includegraphics[width=\linewidth]{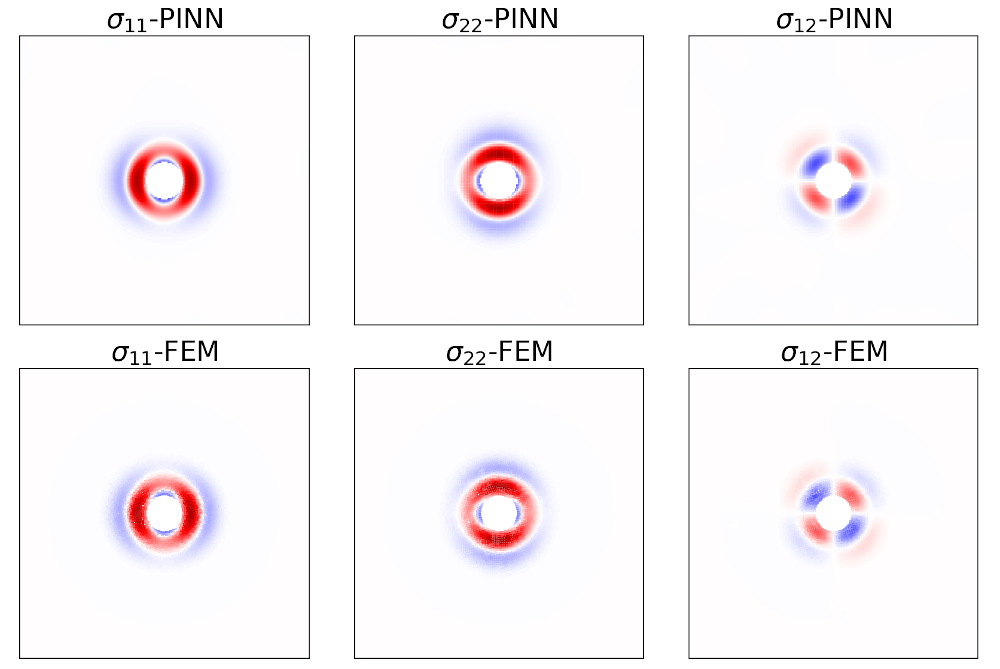}
\end{minipage}%
}%
\subfigure[$t=$ 6 s]{
\begin{minipage}[t]{0.48\linewidth}
\centering
\includegraphics[width=\linewidth]{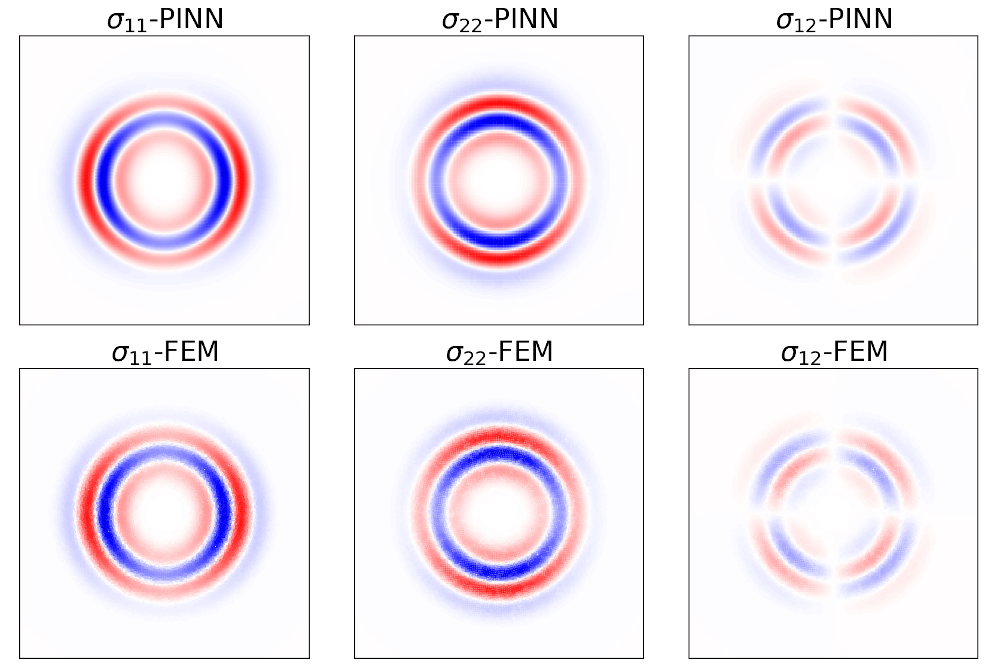}
\end{minipage}%
}%
\newline
\centering
\subfigure[$t=$ 9 s]{
\begin{minipage}[t]{0.48\linewidth}
\centering
\includegraphics[width=\linewidth]{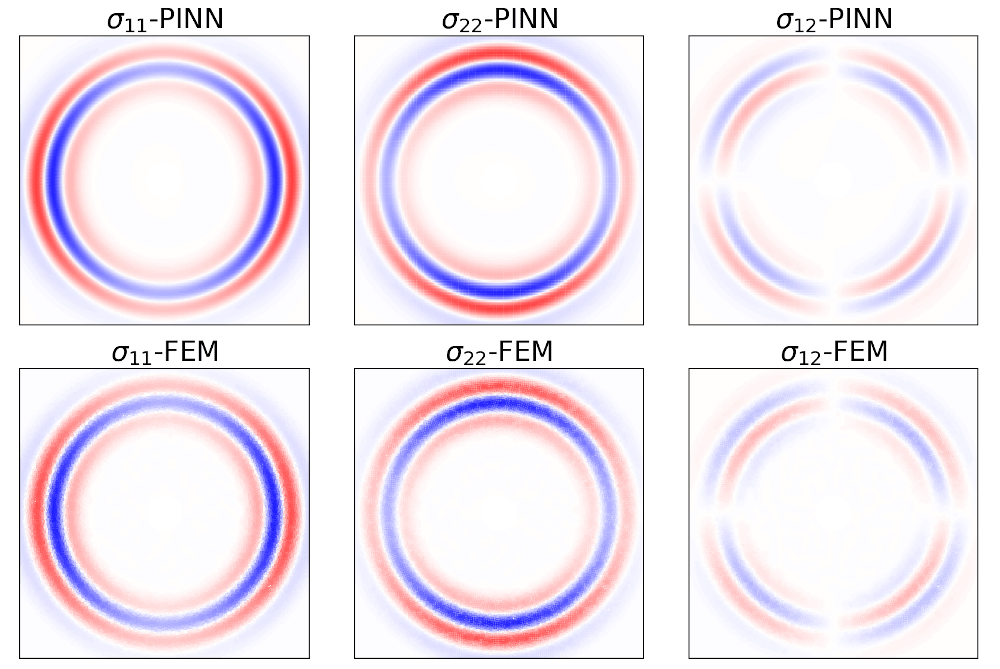}
\end{minipage}
}%
\centering
\subfigure[$t=$ 12 s]{
\begin{minipage}[t]{0.48\linewidth}
\centering
\includegraphics[width=\linewidth]{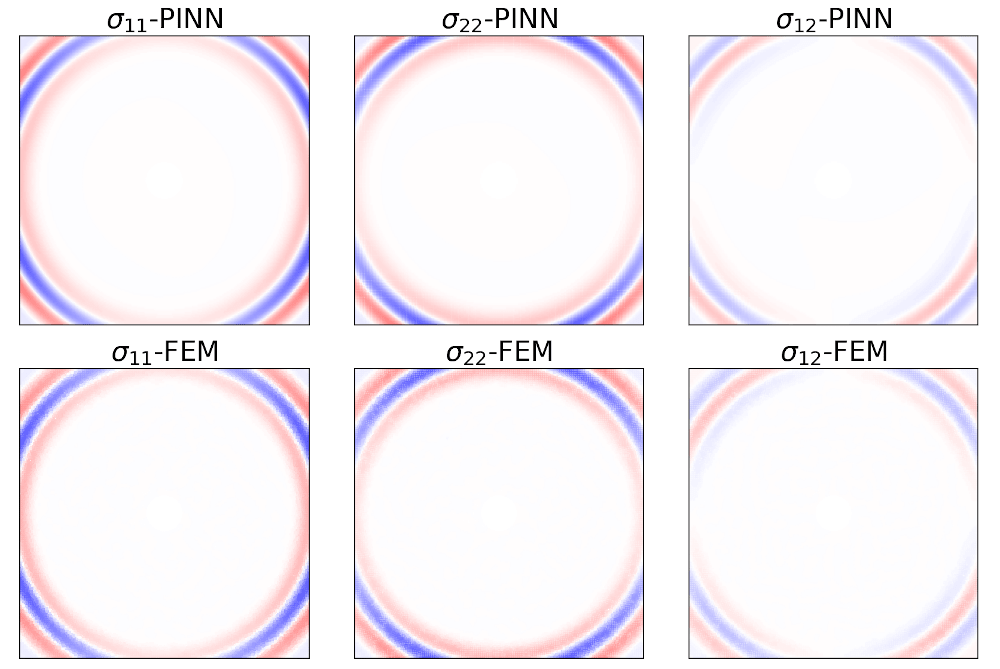}
\end{minipage}%
}%
\newline
\centering
\subfigure{
\begin{minipage}[t]{0.2\linewidth}
\centering
\includegraphics[width=\linewidth]{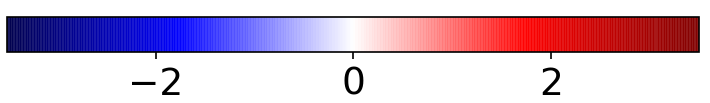}
\end{minipage}%
}%
\caption{Predicted stress fields at various moments in the infinite domain}\label{inf_stress}
\end{figure}

\begin{figure}[t!]
\centering
\includegraphics[width=0.52\textwidth]{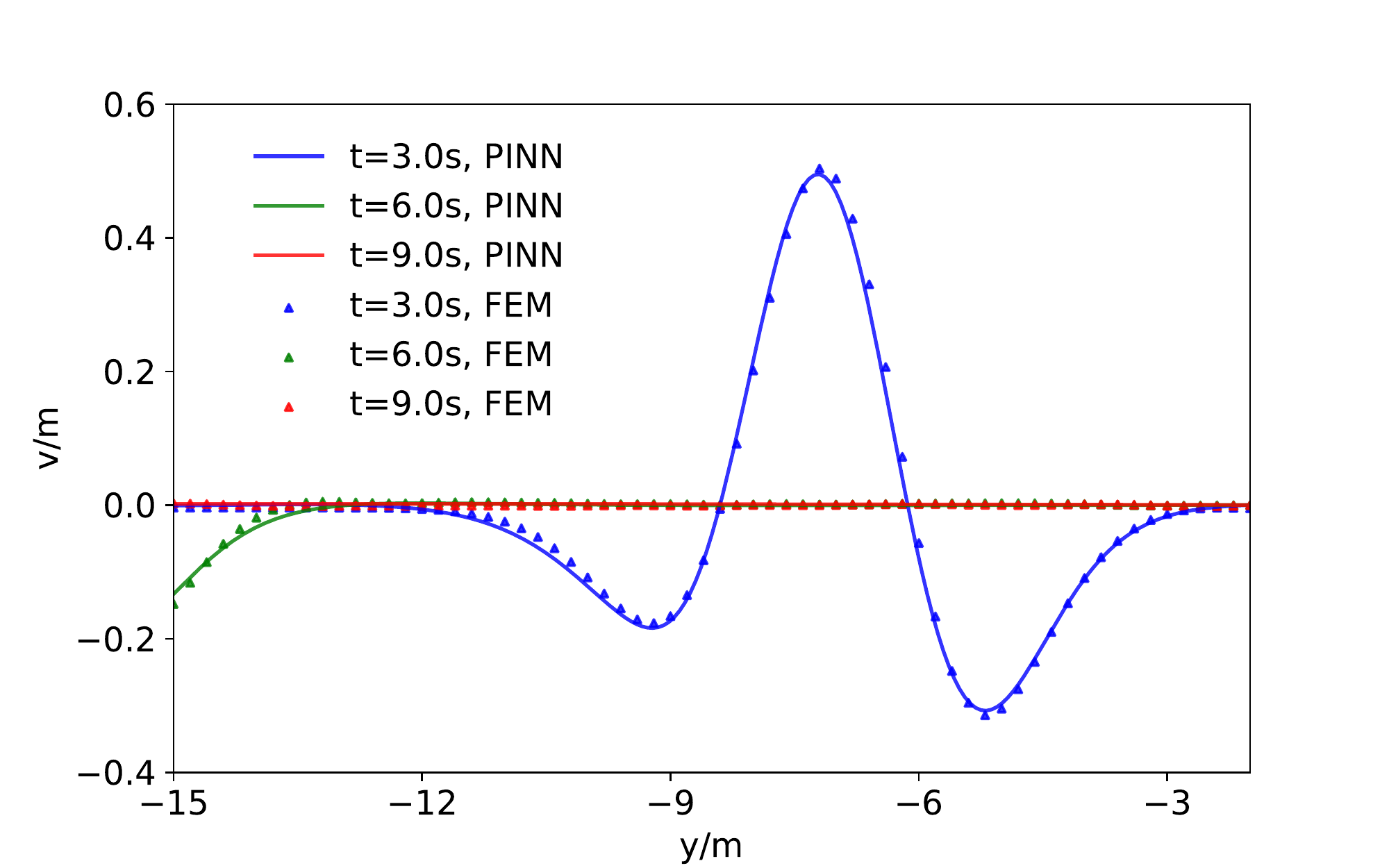}
\caption{Vertical displacements on the mid-line ($x=0$, $y\in[-15,-2]$ m) at various moments in the infinite domain}\label{inf_wave_height}
\end{figure}

In the second case, the wave propagation in an infinite domain is considered. The truncated domain for modeling is given in Fig. \ref{sec3-3-infi-domain}. This case differs from the previous one in that no any boundary conditions are applied except the wave source. That is to say, we do not impose any constraints on the four edges. Since the wave source is the only boundary, we employ the soft enforcement approach (see Fig. \ref{sec2_PINN_diagram}(a)) for the sake of simplicity. The radial displacement $u_r$ in the form of Ricker wavelet is prescribed on the circular wave source, which reads
\begin{equation} 
\label{rickerwave} 
u_r=U_0\left (2\pi^2\left (\frac{t-t_s}{t_p}\right)^2-1 \right)\exp\left[-\pi^2\left (\frac{t-t_s}{t_p}\right)^2\right]
\end{equation}
where $U_0=1.0$ m, $t_s=3$ s and $t_p=3.0$ s. The time duration of the simulation is 16 seconds. A total of 120,000 collocation points are used to evaluate the equation residuals, while 25,600 wave source points and 10,000 initial condition points are employed to evaluate the I/BC constraints generated via LHS sampling. The weighting coefficients $\lambda_1$ and $\lambda_2$ are set to be 1. The configuration of the DNN is $8\times 80$. Figures \ref{inf_disp} and \ref{inf_stress} show the displacement and stress fields of the truncated infinite domain at different moments. The reference result is obtained through the implicit FE solver \cite{abqsolver} for which an enlarged domain ($90\times90$ m) is considered to avoid the wave reflection. It can be observed from the PINN result that the wave shape is not affected by the edge of computation domain ($30\times30$ m) in the context of PINN. The vertical wave (i.e. $v$) distribution on the mid-line ($x=0,~y=[-15,-2]$ m) obtained by PINN is compraed with the reference solution as illustrated in Fig. \ref{inf_wave_height}, which shows excellent agreement.

\begin{figure}[t!]
\centering
\subfigure[$t=$ 6 s]{
\begin{minipage}[t]{0.48\linewidth}
\centering
\includegraphics[width=\linewidth]{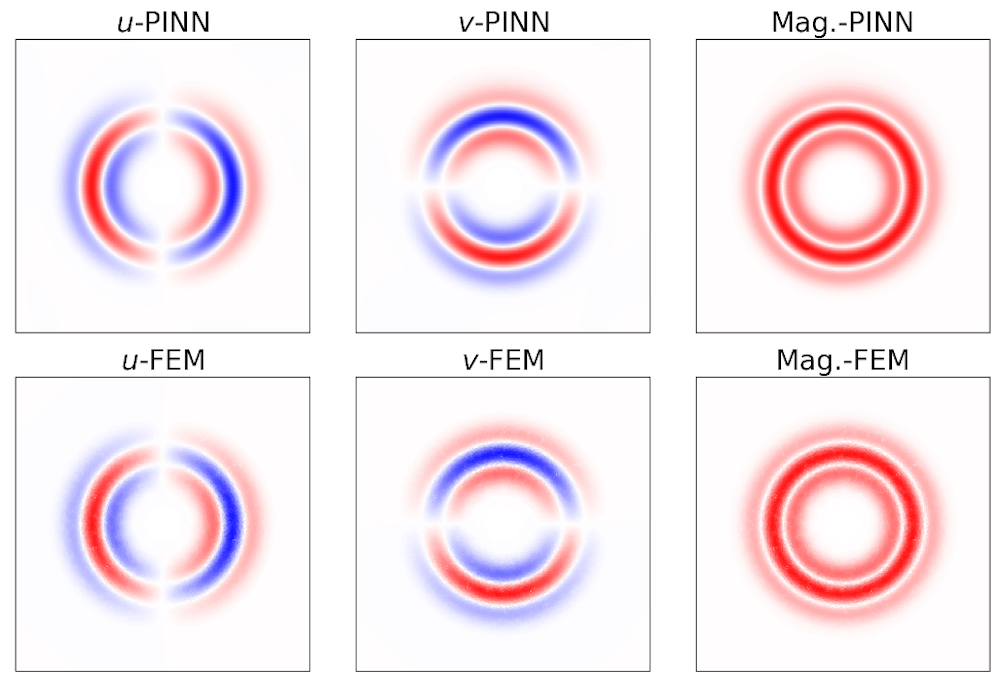}
\end{minipage}%
}%
\subfigure[$t=$ 9 s]{
\begin{minipage}[t]{0.48\linewidth}
\centering
\includegraphics[width=\linewidth]{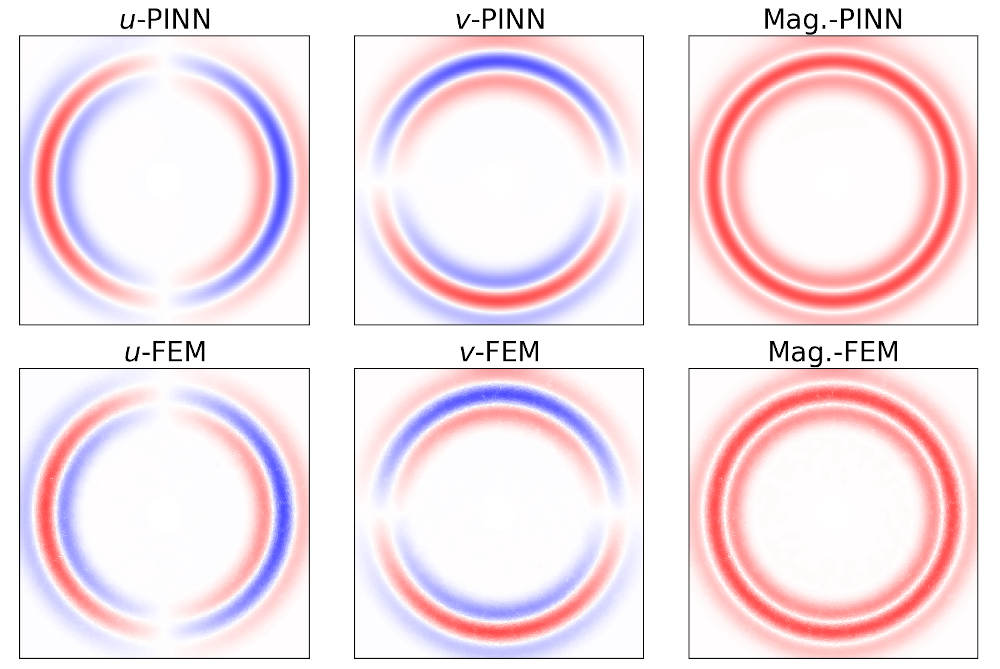}
\end{minipage}%
}%
\newline
\centering
\subfigure[$t=$ 12 s]{
\begin{minipage}[t]{0.48\linewidth}
\centering
\includegraphics[width=\linewidth]{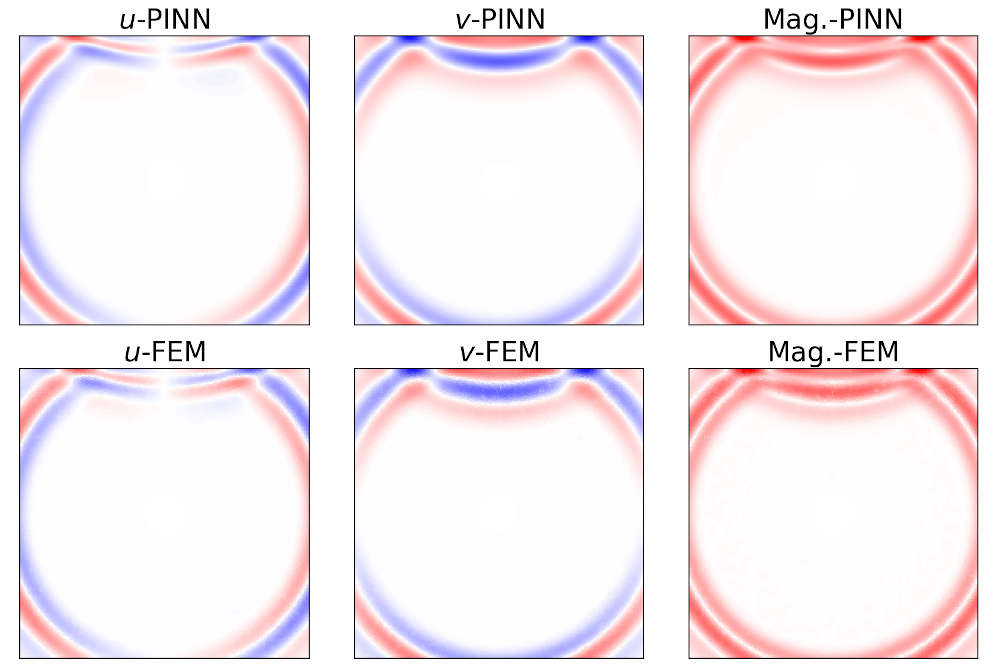}
\end{minipage}
}%
\centering
\subfigure[$t=$ 15 s]{
\begin{minipage}[t]{0.48\linewidth}
\centering
\includegraphics[width=\linewidth]{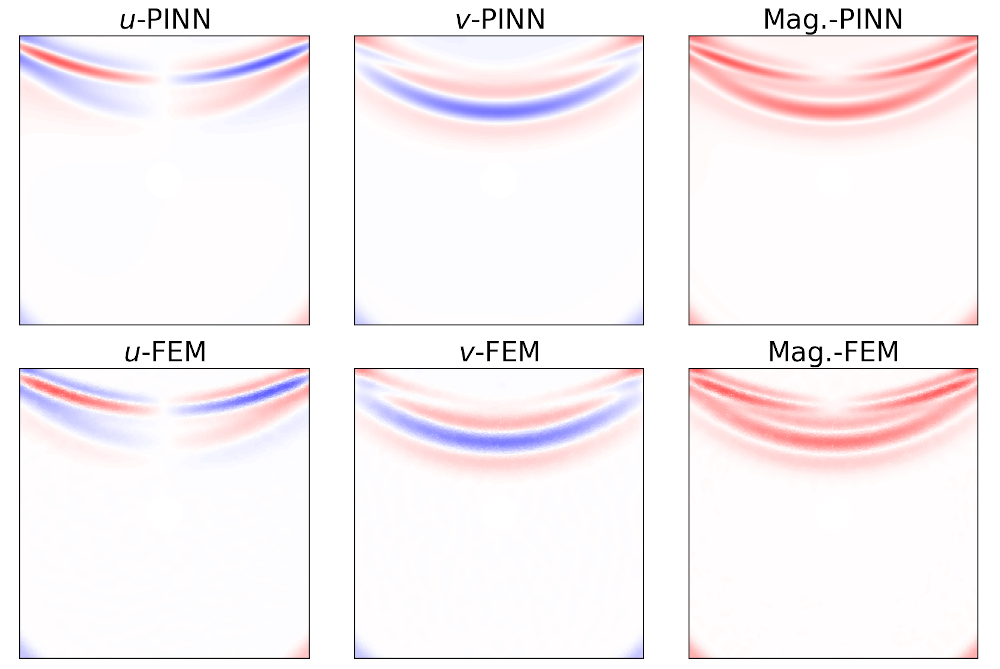}
\end{minipage}%
}%
\newline
\centering
\subfigure{
\begin{minipage}[t]{0.2\linewidth}
\centering
\includegraphics[width=\linewidth]{sec3-3_inf_disp_cb.png}
\end{minipage}%
}%
\caption{Snapshots of the predicted displacement fields in the semi-infinite domain.}\label{semi_inf_disp}
\end{figure}

\subsubsection{Semi-infinite domain}
\begin{figure}[t!]
\centering
\subfigure[$t=$ 6 s]{
\begin{minipage}[t]{0.48\linewidth}
\centering
\includegraphics[width=\linewidth]{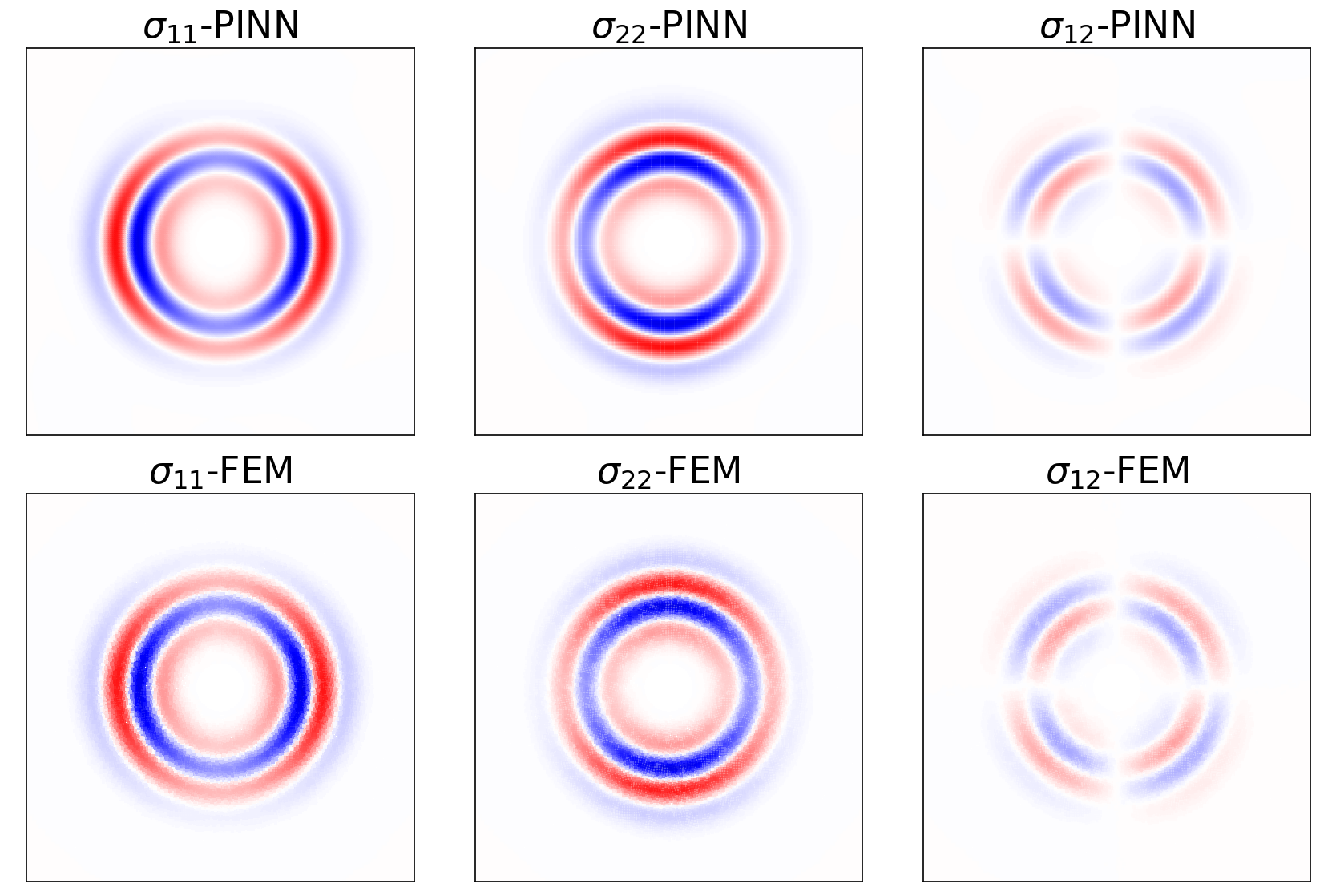}
\end{minipage}%
}%
\subfigure[$t=$ 9 s]{
\begin{minipage}[t]{0.48\linewidth}
\centering
\includegraphics[width=\linewidth]{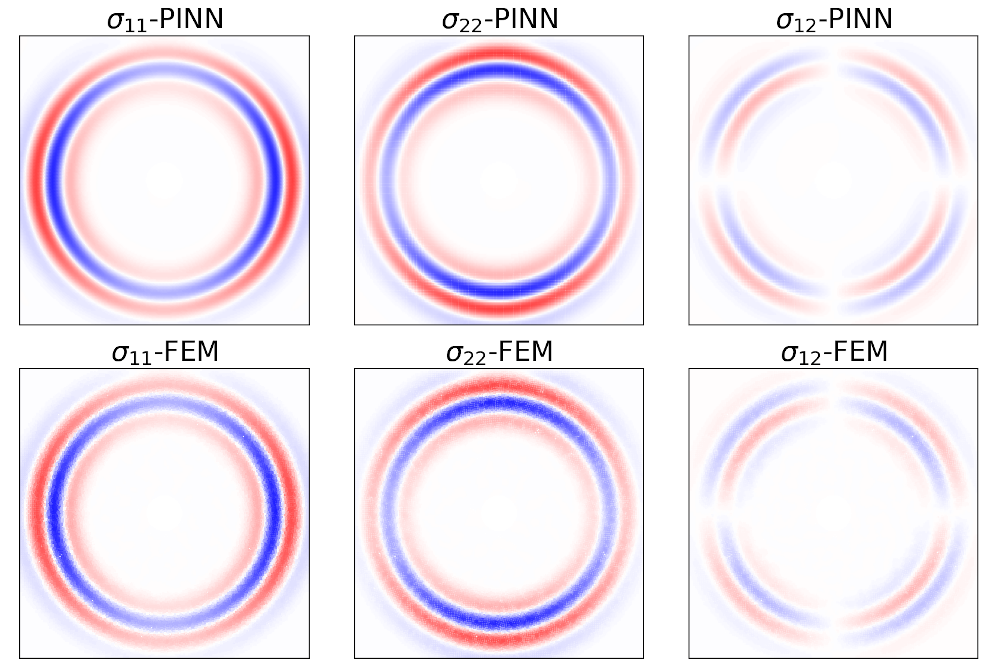}
\end{minipage}%
}%
\newline
\centering
\subfigure[$t=$ 12 s]{
\begin{minipage}[t]{0.48\linewidth}
\centering
\includegraphics[width=\linewidth]{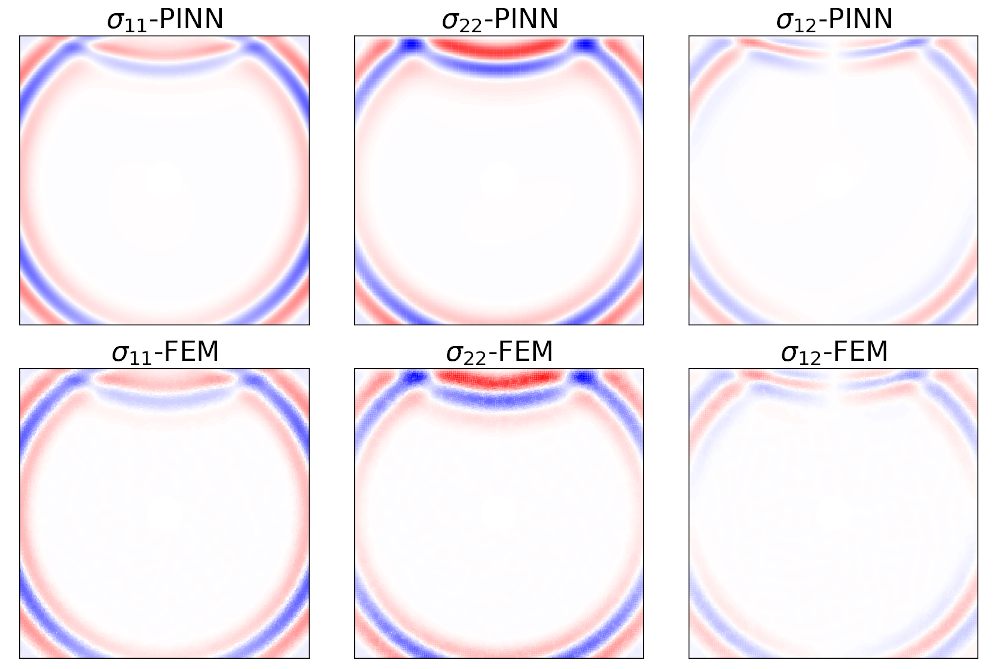}
\end{minipage}
}%
\centering
\subfigure[$t=$ 15 s]{
\begin{minipage}[t]{0.48\linewidth}
\centering
\includegraphics[width=\linewidth]{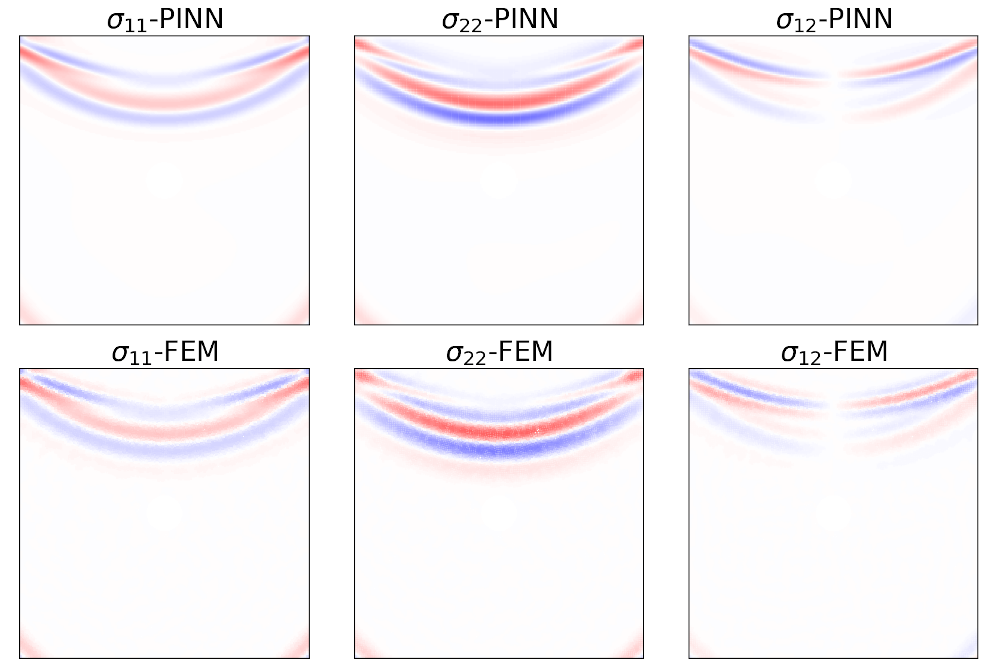}
\end{minipage}%
}%
\newline
\centering
\subfigure{
\begin{minipage}[t]{0.2\linewidth}
\centering
\includegraphics[width=\linewidth]{sec3-3_inf_stress_cb.png}
\end{minipage}%
}%
\caption{Snapshots of predicted stress fields in the semi-infinite domain.}\label{semi_inf_stress}
\end{figure}

\begin{figure}[t!]
\centering
\includegraphics[width=0.5\textwidth]{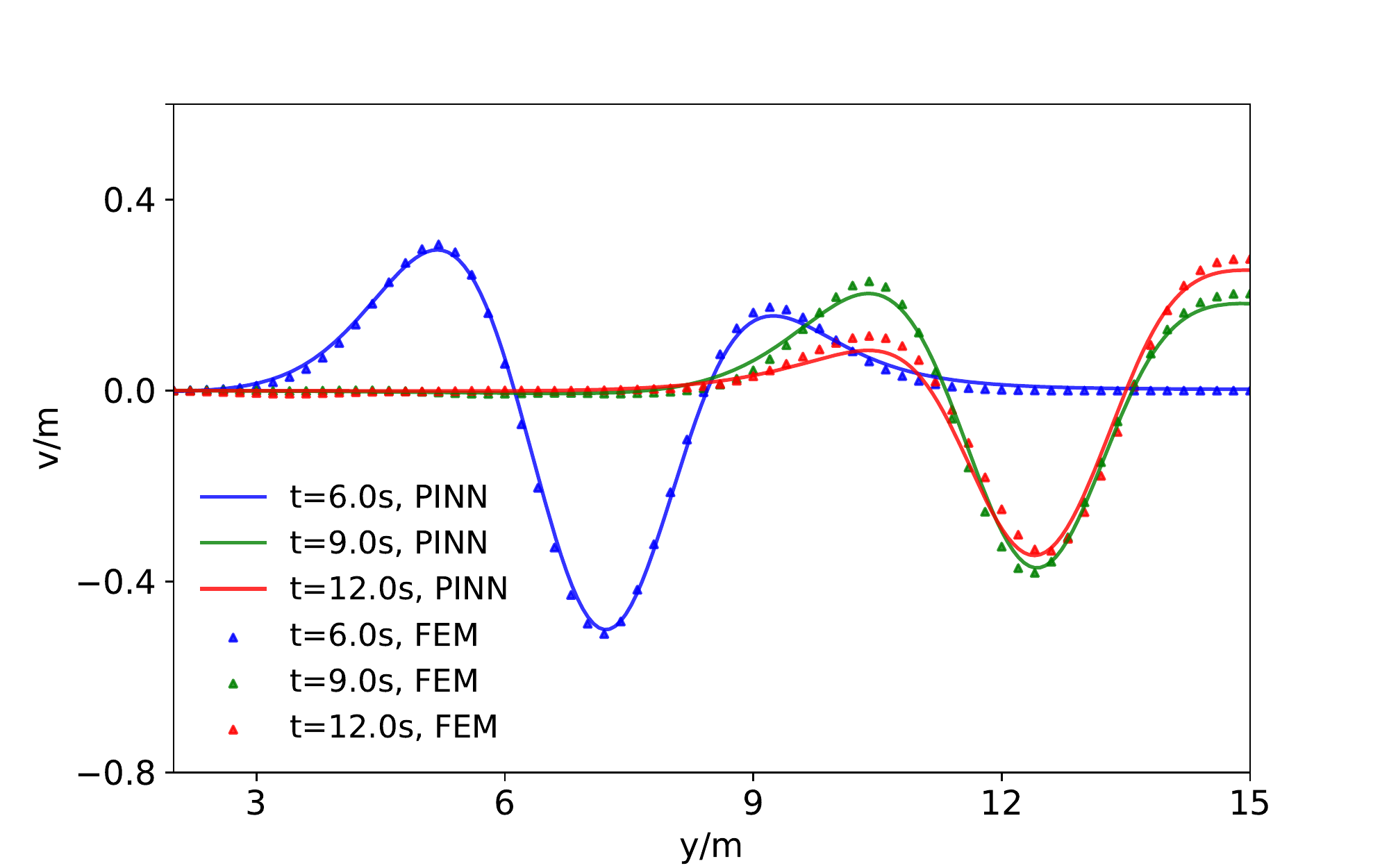}
\caption{Vertical displacement distribution on the mid-line ($x=0$, $y\in[2,15]$ m) at various moments. The truncated semi-infinite domain is considered with free surface.}\label{semi_inf_wave_height}
\end{figure}

In the last case, we consider a semi-infinite domain, i.e., the top edge is modeled as a traction-free boundary condition. This scenario is commonly seen in the modeling of earthquake or underground explosion. The wave source is the same as the second case, defined by Eq. \eref{rickerwave}. The I/BCs are enforced in a ``soft'' way (see Fig. \ref{sec2_PINN_diagram}(a)) due to the simple boundary conditions considered herein. The network with $8\times 100$ layers/neurons is trained with 12,0000 collocation points for equation residuals, 10,000 points for the initial state, 30,000 points on the wave source and 15,000 points on the top edge. It is noted that the equation residual points are refined near the wave source and the free surface to better capture the details of the wave. The results of the wave propagation are presented in Fig. \ref{semi_inf_disp} and \ref{semi_inf_stress}, in comparision with the implicit FE solution obtained from an enlarged domain ($90\times 60$ m) whose other three edges are fixed. It can be seen that the PINN is able to predict the surface reflection while the other three edges have no interfere on course of the wave propagation. To quantitatively examine the accuracy of the PINN result, we compare the vertical displacement distribution on the mid-line ($x=0$, $y\in [2,15]$) at various moments with that of the FE solution, as shown in Fig. \ref{semi_inf_wave_height}, which match well with each other.

\section{Conclusions}\label{conc}
In this work, we present a PINN framework for solving elastodynamics problems without using labeled data (although measurement data can also be incorporated when available). In particular, we propose a mixed-variable scheme of PINN where both the displacement and stress components are taken as the output of the neural network. This scheme is found to significantly improve the accuracy and trainability of the network compared with the pure displacement-based PINN (see \ref{mix_vs_disp}). We also propose a composite scheme of PINN to construct a synergy solution to elastodynamics problems. The basic concept is to enforce the I/BCs in a ``hard'' manner through decomposition of the solution to a PDE system into general and particular solutions. The major benefit of the constructed solution is that the I/BCs are imposed forcibly and satisfied in nature. The parametric study (see \ref{soft_vs_hard}) shows that the inaccuracy near the boundaries encountered by the conventional PINN with ``soft'' I/BC enforcement can be mitigated significantly. A series of dynamics problems, including the defected plate under cyclic uni-axial tension and the elastic wave propagation in confined/semi-infinite/infinite domains, are studied to illustrate the capability of the proposed PINN. Since the proposed PINN deals with the strong-form PDEs for wave propagation modeling, the wave field is not affected by the unconstrained boundaries when it propagates out of the truncated domain, hence avoiding the wave reflection issue commonly seen in Galerkin-based methods. The proposed method shows great potential in full wave-inversion problems which will be studied in the future. In a broader sense, the proposed PINN framework can solve general PDE systems, despite linear or nonlinear, determined or stochastic.

In summary, this paper has discussed some basic issues of PINN for solving solid mechanics problems, such as I/BC enforcement, the numerical scheme (governing PDEs and unknown variables) and the formulation of loss functions. In addition to the properties discussed above, PINN has been characterized with many other advantages which fall out of the scope of this paper, including (1) the capability to achieve the data-driven solution when measurements are available \cite{raissi2019physics}, (2) the convenience to formulate the inverse problem in mechanics (e.g., crack detection \cite{jung2013identification,jung2014modeling,sun2016sweeping} and material damage model calibration \cite{liu2014regularized}) due to the optimization nature of training PINN \cite{raissi2019physics}, and (3) the capability to conduct uncertainty quantification (UQ) on physical systems \cite{yang2019adversarial, zhu2019physics}. The proposed PINN framework can also be potentially extended to take into account the aforementioned capabilities.

\bibliographystyle{elsarticle-num}
\bibliography{refs}

\appendix
\setcounter{figure}{0}

\begin{figure}[b!]
\centering
\includegraphics[width=0.6\textwidth]{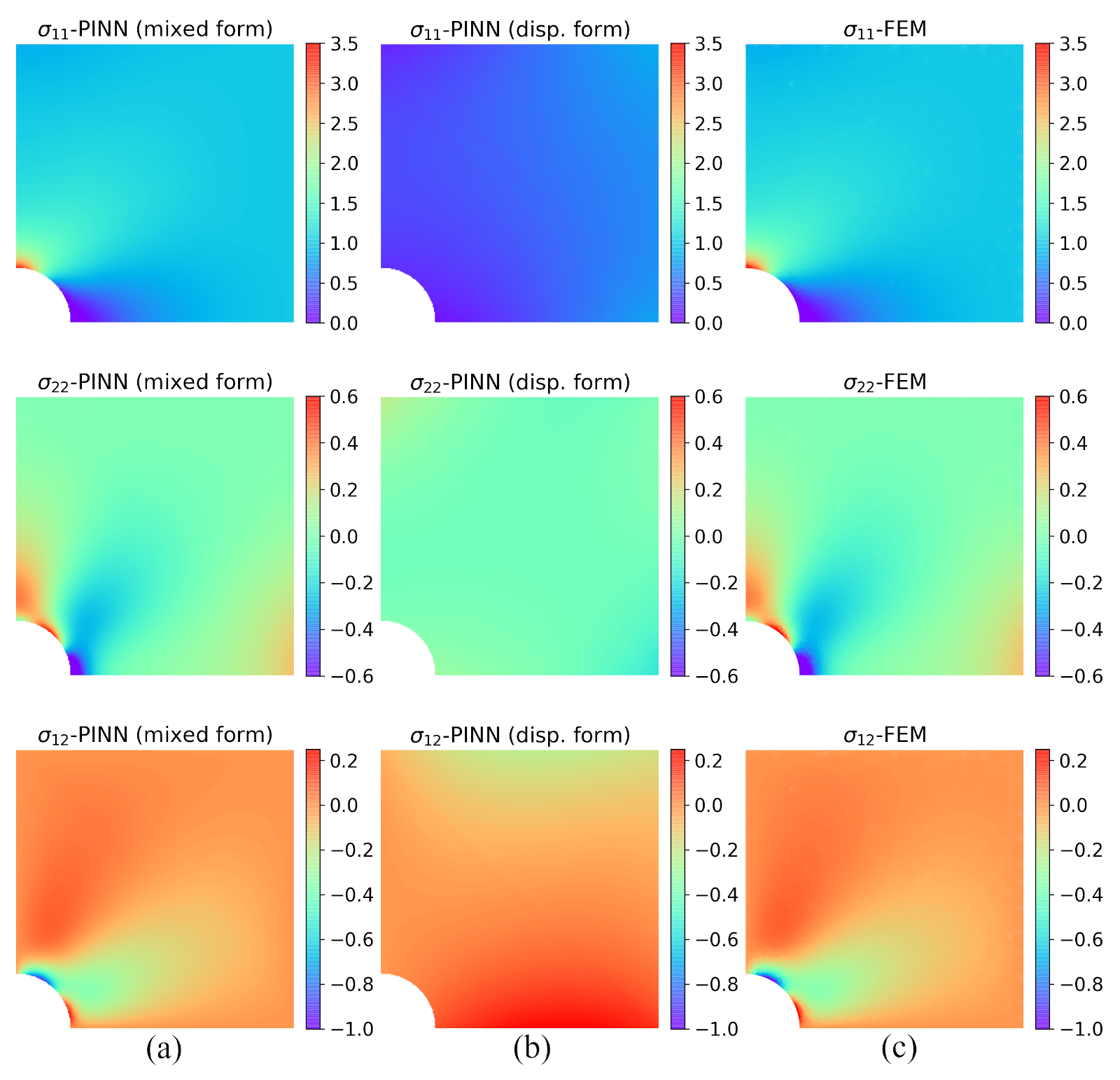}
\caption{Comparison of the stress fields produced by (a) the mixed-variable PINN, (b) the pure-displacement PINN, and (c) the finite element solution. Both networks have 6 hidden layers and 60 neurons in each layer.} 
\label{mixed_disp_compare}
\end{figure}

\section{Comparison between mixed-variable/pure-displacement formulation}\label{mix_vs_disp}

\begin{figure}[t!]
\centering
\includegraphics[width=0.7\textwidth]{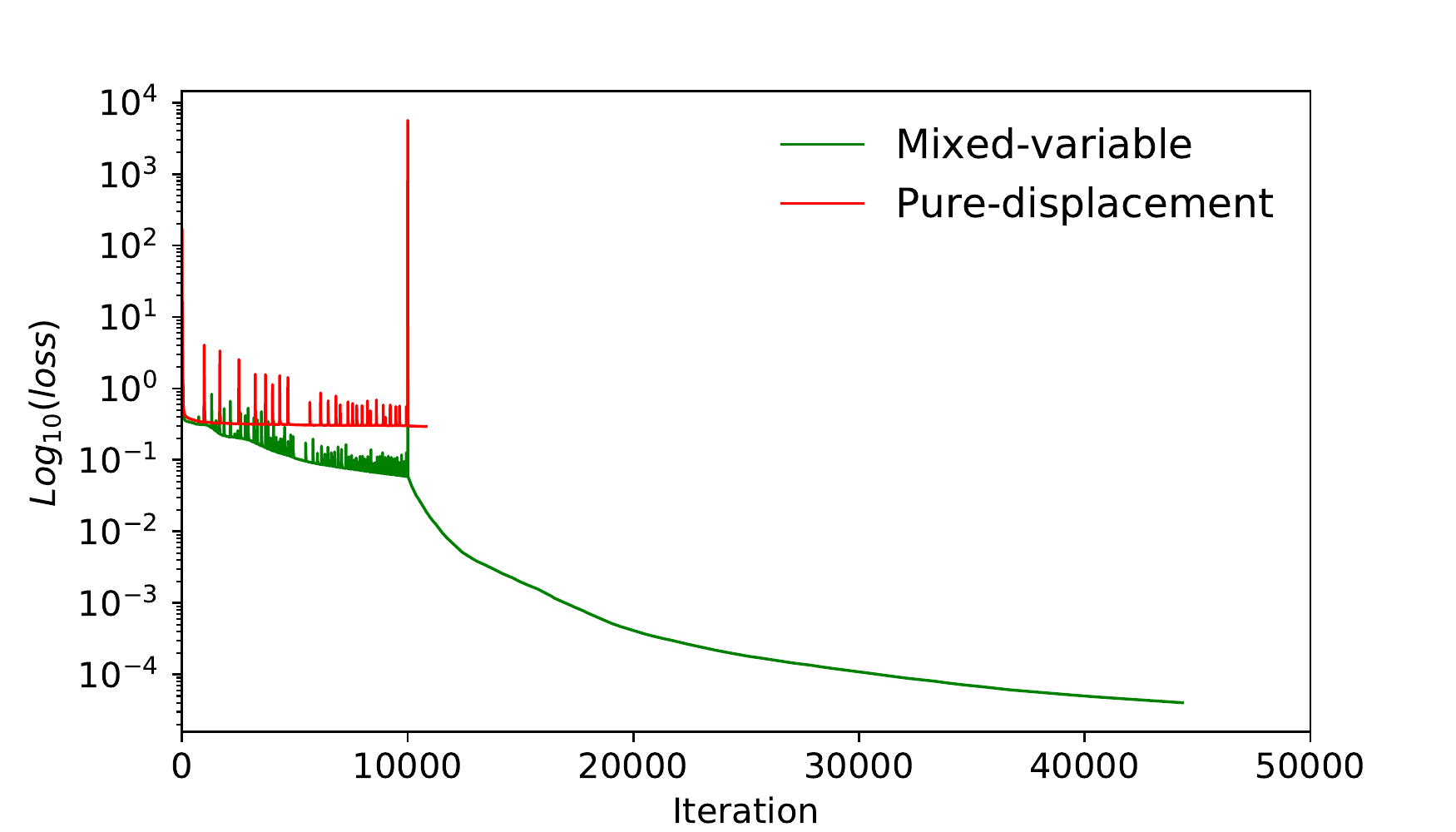}
\caption{Convergence of the loss functions. Both PINNs (mixed-variable vs. pure-displacement) are trained with the Adam optimizer for 10,000 iterations (learning rate $10^{-3}$), followed by the L-BFGS-B optimizer. The final loss values are $2.89\times10^{-1}$ and $2.93\times10^{-5}$ for mixed-variable and pure-displacement PINN, respectively}\label{loss_curve_compare}
\end{figure}

In this part, we will show the comparison of the results produced by mixed-variable scheme and traditional displacement-based formulation for PINN. The mixed-variable network maps the spatiotemporal location $(x,y)$ to the displacement and stress variables $(u,v,\sigma_{11},\sigma_{22},\sigma_{12})$, which include the dependent variables $\bm{\sigma}$ compared with the pure displacement output $(u,v)$. This major difference in output will also affect the construction of the loss function. For the current static case, the loss function in the mixed-variable PINN is given in Eqs. \eref{msef} and \eref{msebc} while the counterparts in the pure-displacement form read 
\begin{equation} 
\label{msef_disp} 
J_\text{g}=||\bm{\nabla} \cdot (\bm{\mathbf{C}:\bm{\epsilon}}) +\mathbf{F}-\rho\mathbf{u}_{tt}||^2_{\Omega\times \left [ 0,T \right ]}
\end{equation}
\begin{equation} 
\label{msebc_disp} 
J_\text{bc}=||\mathbf{u}-\mathcal{B}_D||^2_{\partial \Omega_D\times \left [ 0,T \right ]}+||\mathbf{n\cdot} (\mathbf{C}:\bm{\epsilon})-\mathcal{B}_N||^2_{\partial \Omega_N\times \left [ 0,T \right ]}
\end{equation}
It can be seen that the mixed-variable PINN reduces the highest order of spatial derivatives from two to one considering the displacement-strain relationship $\displaystyle \bm{\epsilon } =\big [ \bm{\nabla}\mathbf{u} + {(\bm{\nabla} \mathbf{u})}^T \big ]\big/2$, which we believe is the major reason for its improved trainability and accuracy. Another benefit of the mixed-form output is that, in addition to the Dirichlet boundary, we are able to impose the Neumann boundary condition in a ``hard'' way (see Section \ref{elas_theory}) which is infeasible for the pure-displacement PINN \cite{sun2020surrogate}. The static case described in Fig. \ref{sec3-2-sketch} is considered herein with traction $T_n(t)=1.0$ MPa. We keep the same hyperparameters and collocation points, as described in Section \ref{notchplate}, for the two networks. The comparison of the produced stress fields is shown in Fig. \ref{loss_curve_compare}, in comparison with the FE reference result. The pure-displacement PINN fails to model the problem while the mixed-variable PINN gives satisfactory prediction. The convergence of the training loss is presented in Fig. \ref{loss_curve_compare}. It can be seen that the loss reaches a stagnation soon after the training is launched for the pure-displacement PINN, while the mixed-variable PINN demonstrates excellent trainability as indicated by the convergence curve.

\section{Comparison between soft/hard enforcement of boundary conditions}\label{soft_vs_hard}

\begin{figure}[t!]
	\centering
	\includegraphics[width=0.4\textwidth]{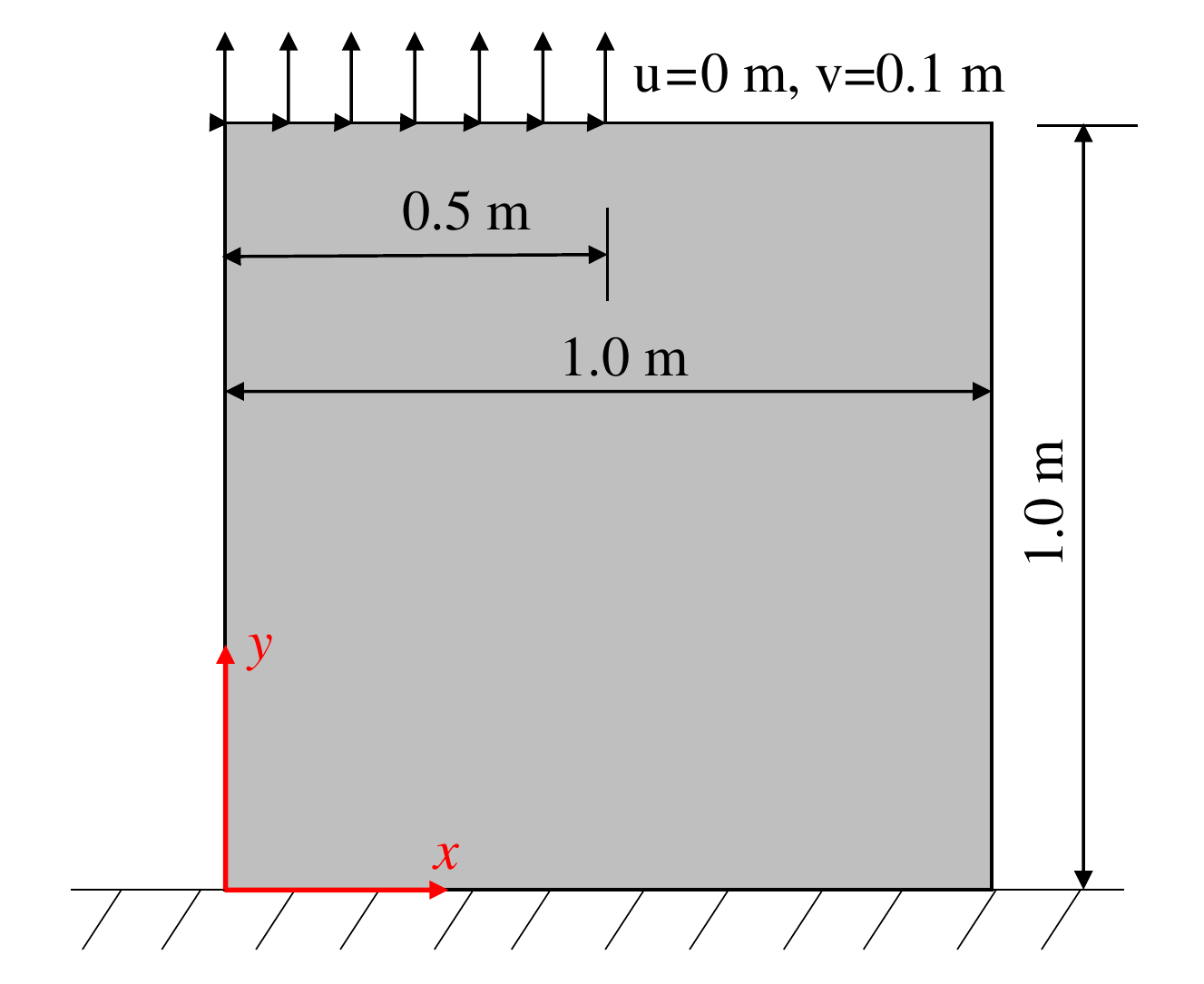}
	\caption{Diagram of the plane stress problem of a square plate with half edge loaded.}
	\label{ap_a_sketch}
\end{figure}

\begin{figure}[t!]
\centering
\subfigure[DNN approximation]{
\begin{minipage}[t]{0.9\linewidth}
\centering
\includegraphics[width=\linewidth]{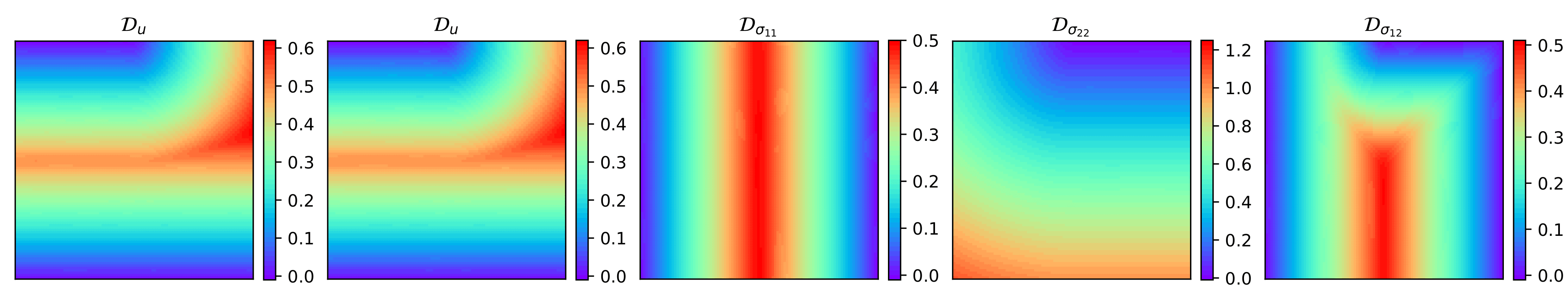}
\end{minipage}%
}%
\hfill
\subfigure[Analytical]{
\begin{minipage}[t]{0.9\linewidth}
\centering
\includegraphics[width=\linewidth]{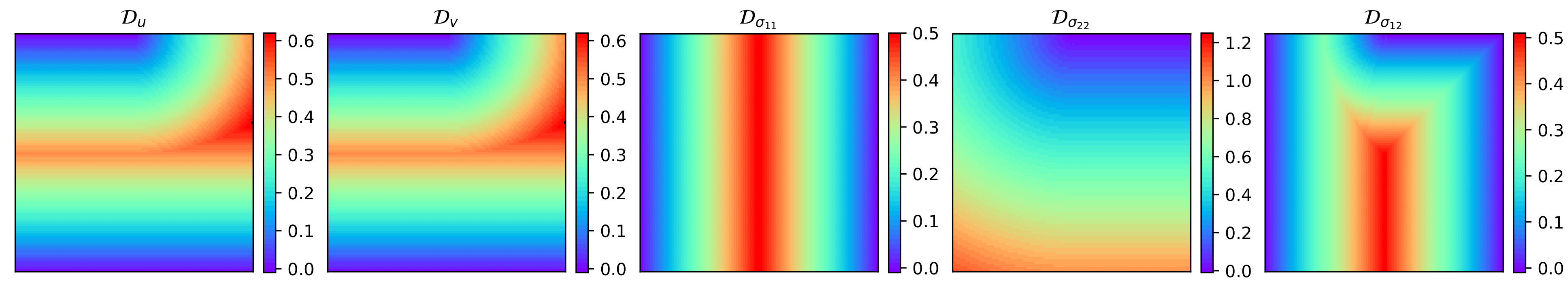}
\end{minipage}%
}%
\caption{Distance function $\mathcal{D}(x,y)$: (a) DNN approximation, (b) analytical value.}
\label{ap_a_dist_func}
\end{figure}

\begin{figure}[t!]
\centering
\subfigure[$\lambda=1$]{
\begin{minipage}[t]{0.18\linewidth}
\centering
\includegraphics[width=\linewidth]{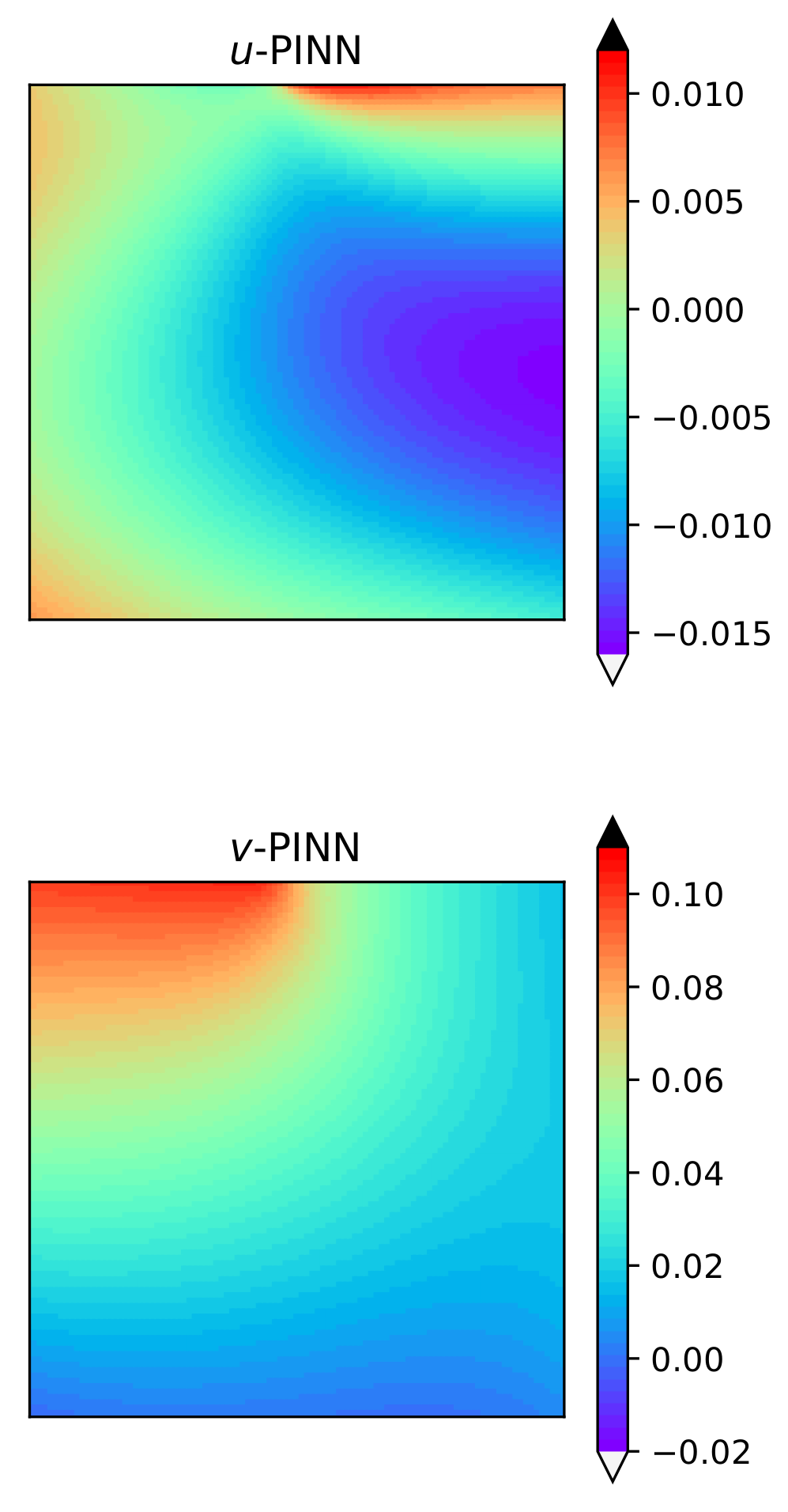}
\end{minipage}%
}%
\subfigure[$\lambda=10$]{
\begin{minipage}[t]{0.18\linewidth}
\centering
\includegraphics[width=\linewidth]{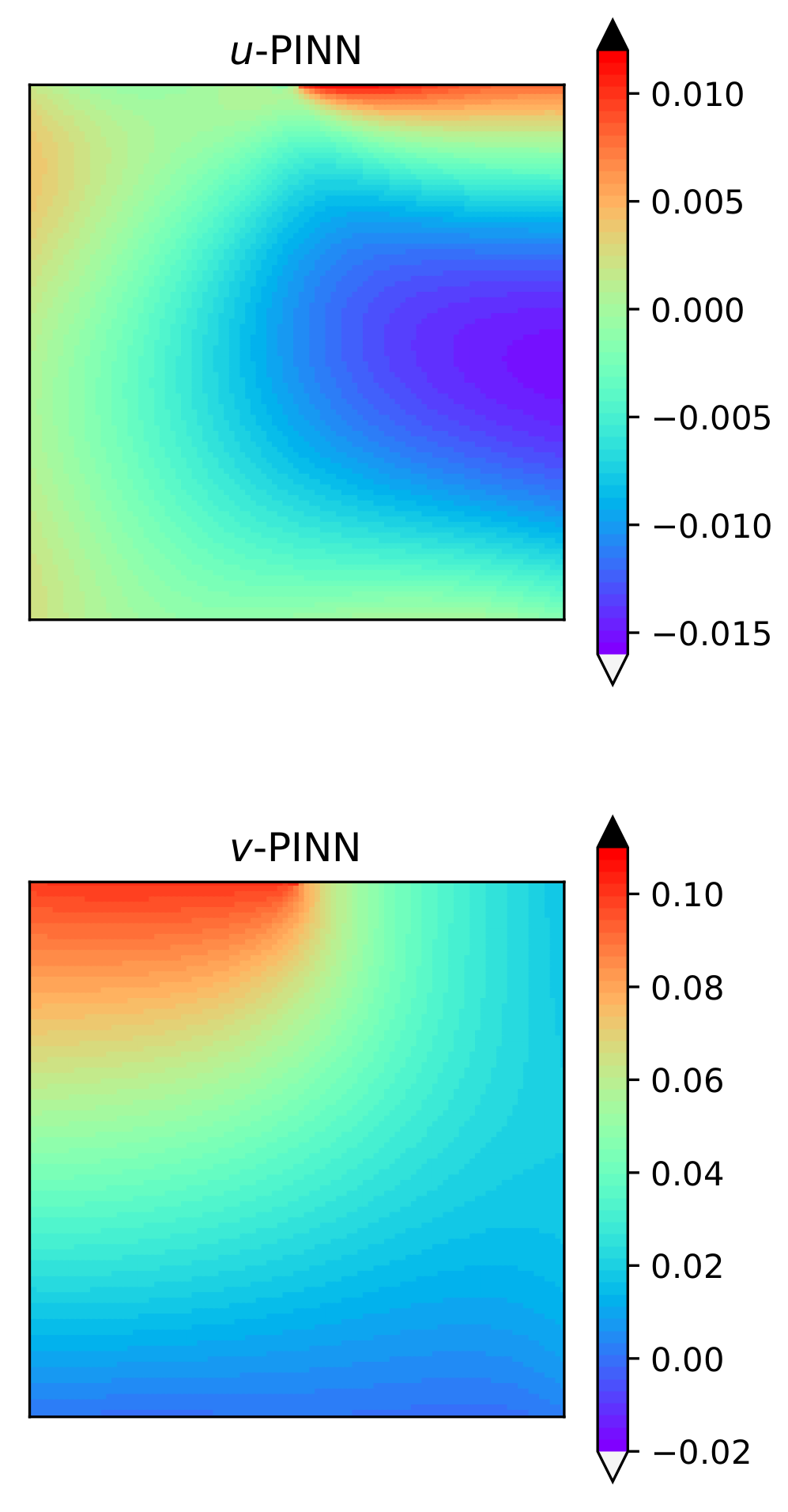}
\end{minipage}%
}%
\subfigure[$\lambda=100$]{
\begin{minipage}[t]{0.18\linewidth}
\centering
\includegraphics[width=\linewidth]{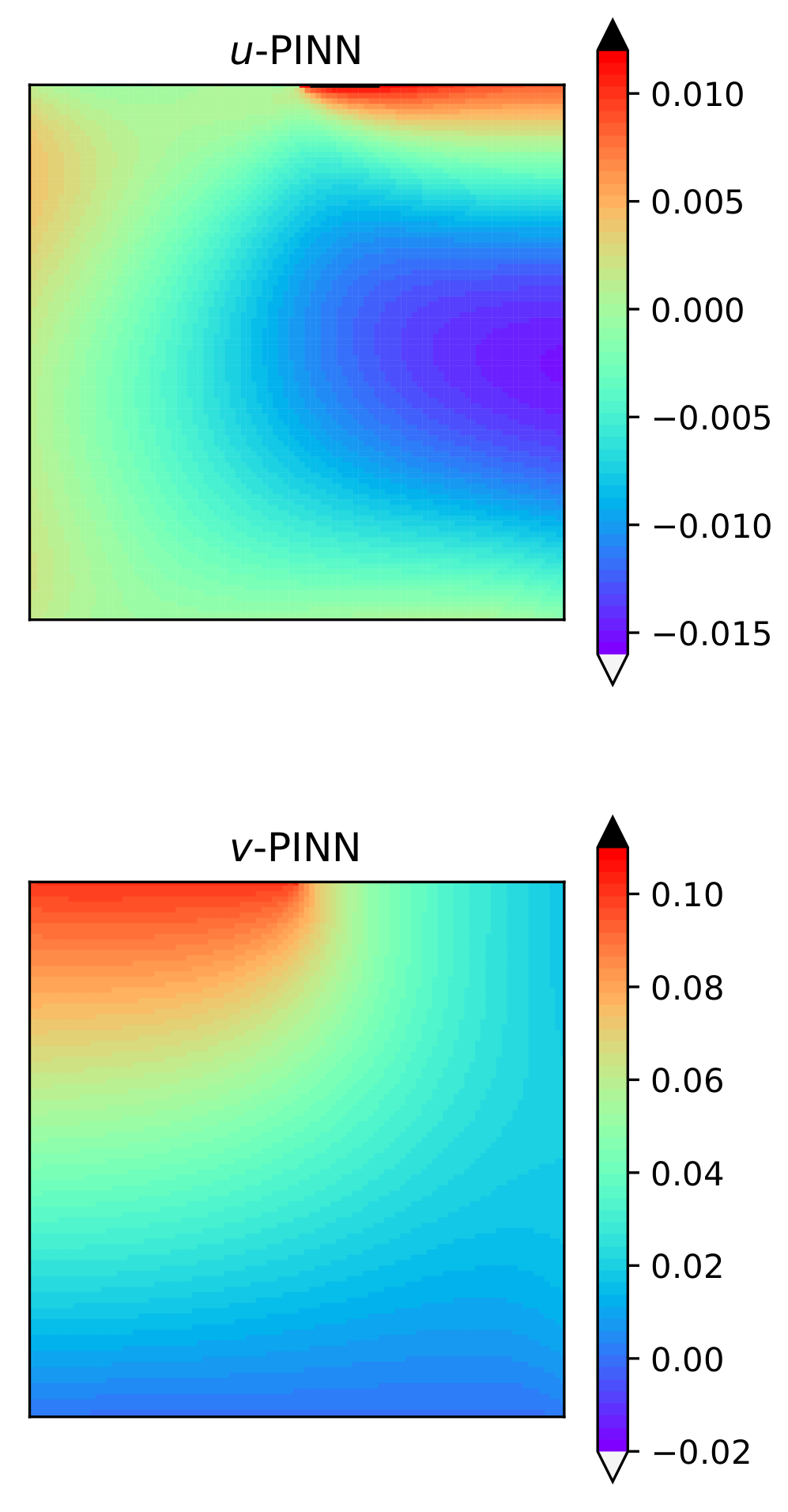}
\end{minipage}
}%
\subfigure[Hard]{
\begin{minipage}[t]{0.18\linewidth}
\centering
\includegraphics[width=\linewidth]{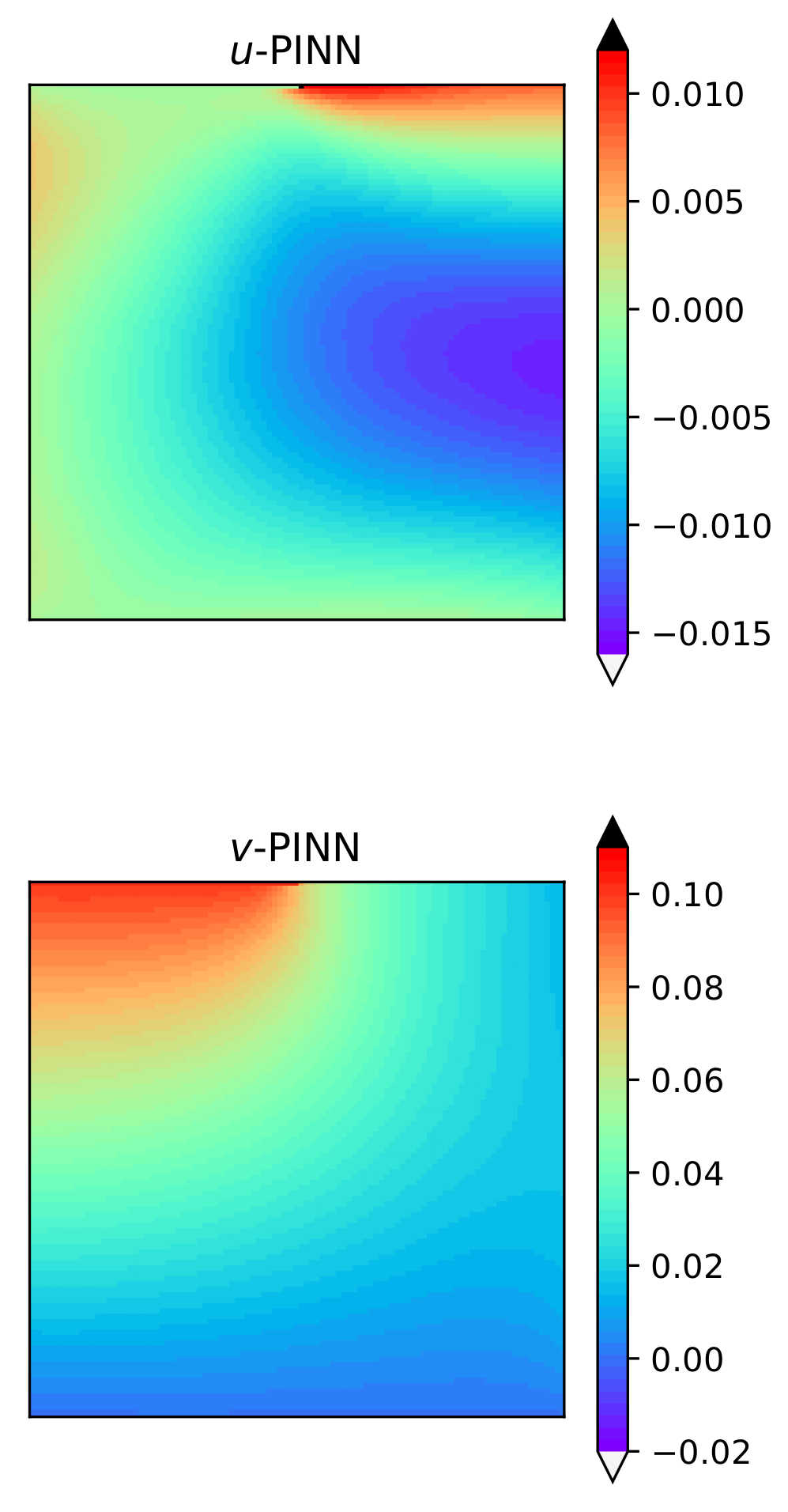}
\end{minipage}
}%
\subfigure[FEM]{
\begin{minipage}[t]{0.18\linewidth}
\centering
\includegraphics[width=\linewidth]{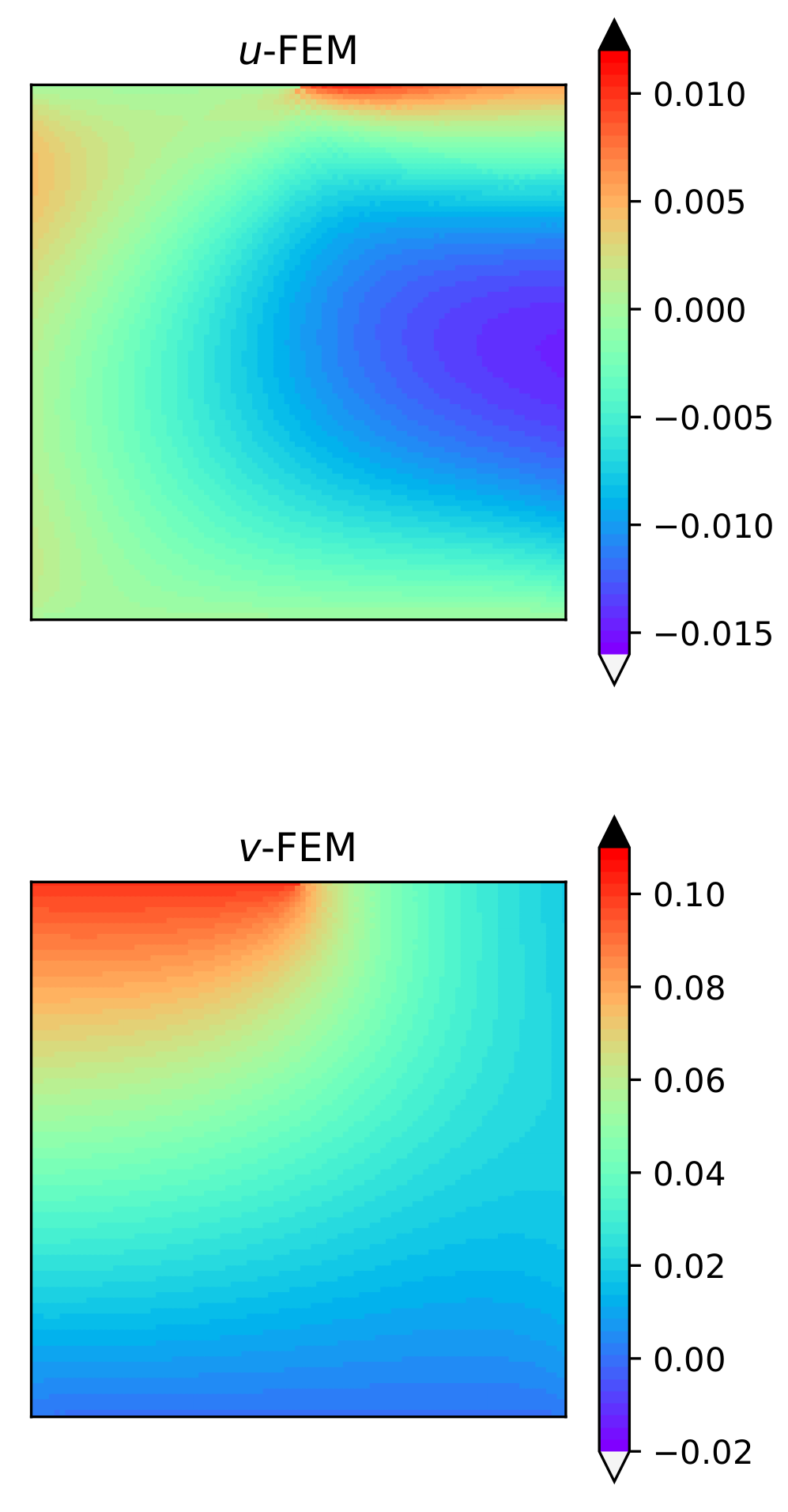}
\end{minipage}
}%
\caption{Comparison of the displacement fields for various configurations: (a), (b) and (c) are predictions from soft-enforced PINN with various weighing coefficients, while (d) uses ``hard'' enforcement, and (e) is obtained from the FE solver with 10,000 linear quadrilateral elements.}
\label{ap_a_disp_field}
\end{figure}

The performance of the proposed approach for boundary condition (BC) enforcement, i.e., the ``hard'' enforcement, is compared with its counterpart in the conventional PINN, soft enforcement. A two-dimensional plane stress problem, as shown in Fig. \ref{ap_a_sketch}, is considered for the sake of simplicity. The square plate has its lower edge fixed while the half of its top edge is applied with the forced displacement $(u,v)$=(0, 0.1) m. The Young's modulus and the Poisson's ratio of the plate are 10 MPa and 0.2 respectively. A total number of 10,000 collocation points, which includes 150 Dirichlet boundary (lower and upper-left edge) points and 250 Neumann boundary (left, right and upper-right edge) points, are generated using LHS sampling for training the network. The Adam and L-BFGS-B optimizers are employed subsequently to train the networks. As for the network architecture, we adopt $6\times30$ for the conventional PINN, while, for the proposed composite PINN, the networks for $\mathcal{U}_p$, $\mathcal{D}$ and $\mathcal{U}_h$ are configured to be $3\times20$, $3\times20$ and $6\times30$, respectively, with tanh($\cdot$) as the activation function. Simulations are also conducted with various values of the weighting coefficient $\lambda$ on the BC residual, for the soft-enforced PINN. The training data for the distance function, which comes from the analytical value, as well as the trained DNN's approximation, is shown in Fig. \ref{ap_a_dist_func}. 

In Fig. \ref{ap_a_disp_field}, we compare the displacement fields for different cases. The FE solution is provided as the reference. As can be seen from Fig. \ref{ap_a_disp_field}(a), the BC at the lower-left corner, i.e., $u(x,1)=0$, is not enforced accurately for the soft enforcement approach with $\lambda=1$. However, increasing the $\lambda$ can mitigate the inaccuracy as shown in Fig. \ref{ap_a_disp_field}(b-c). It is due to the unbalance between each term within the loss function, resulting in the gradient pathology issue during the training \cite{wang2020understanding}. Therefore, a trail-and-error procedure is usually involved in training the soft-enforced PINN to find a suitable coefficient, which is computationally costly. For the proposed PINN, this issue does not exist since the solution is constructed in a way that the BC is imposed forcibly. 


\end{document}